\documentclass[a4paper,10pt]{article}
\usepackage[utf8]{inputenc}
\DeclareUnicodeCharacter{00A0}{}
\title{Espaces de Berkovich, polytopes, squelettes et théorie des modèles}
\author{\sc Antoine Ducros\thanks{L'auteur est membre du projet ANR {\em Espaces de Berkovich}}\\ \small Institut de mathématiques de Jussieu \\ \small Projet {\em Topologie et géométrie algébriques}\\ \small 4 place Jussieu 75005 Paris ~~FRANCE}
\date{}

\usepackage{mathrsfs}

\usepackage[francais]{babel}
\usepackage{amssymb}
\usepackage{hyperref} 

\date{}
\input xy
\xyoption{all} 
\input{xypic}

\newcommand{\pl}[3]{\mathsf{PL}_{#1}(#2,#3)}
\newcommand{\wid}{\widehat}
\newcommand{\lgk}[2]{\Lambda_{#2}(#1)}
\newcommand{\rst}[1]{(\RR^*_+)^{#1}}

\newcommand{\inv}{^{-1}}

\renewcommand{\leq}{\leqslant}

\renewcommand{\geq}{\geqslant}
\newcommand{\grot}{_{\rm G}}

\newcommand{\hotimes}{\wid{\otimes}}
\renewcommand{\phi}{\varphi}
\renewcommand{\dim}[1]{\mbox{\rm dim}_{#1}\;}

\newcommand{\got}[1]{{\mathfrak #1}}
\renewcommand{\Bbb}{\mathbb}

\newcommand{\gm}{{\Bbb G}_{m}}
\newcommand{\gmk}{{\Bbb G}_{m,k}}
\newcommand{\gman}{\gmk^{n,\rm an}}

\newcommand{\RR}{{\Bbb R}}
\newcommand{\GG}{{\Bbb G}}

\newcommand{\sk}{{\bf S}}
\newcommand{\ZZ}{{\Bbb Z}}

\newcommand{\FF}{{\Bbb F}}

\newcommand{\PP}{{\Bbb P}}
\newcommand{\Aff}{{\Bbb A}}

\newcommand{\sch}[1]{\mathscr #1}

\newcommand{\NN}{{\Bbb N}}

\newcommand{\QQ}{{\Bbb Q}}
\renewcommand{\epsilon}{\varepsilon}

\newcommand{\an}{^{\rm an}}
\newcommand{\zero}{^{\mbox{\tiny o}}}

\newcommand{\spec}{\mathsf{Spec}\;}

\newcommand{\red}{\widetilde}
\newcommand{\hres}{{\sch H}}
\newcommand{\cf}{{\em cf}}

\newcommand{\val}[1]{|{#1}|}
\newcommand{\deux}[1]{\refstepcounter{subsection}\label{#1}\medskip\noindent {\bf (\thesubsection)}\hspace{.1cm}}
\newcommand{\trois}[1]{\refstepcounter{subsubsection}\label{#1}\medskip\noindent {\bf
    (\thesubsubsection)}\hspace{.1cm}}

\begin{document}
\maketitle

{\small
\noindent
{\bf Abstract}. Let~$X$ be an analytic space over a non-Archimedean, complete field~$k$ and let~${\bf f}=(f_1,\ldots, f_n)$ be a family of invertible functions on~$X$. Let us recall two results, both of which were proven using de Jong's alterations (these alterations could have been avoided for 1), which could have been deduced quite formally from a former result by Bieri and Groves, based upon explicit computations on Newton polygons). 

\medskip
1) The compact set~$|{\bf f}|(X)$ is a polytope of the~$\RR$-vector space~$\rst n$ (we use the multiplicative notation) ; this is due to Berkovich in the locally algebraic case, and has been extended to the general case
by the author.

2) If moreover~$X$ is Hausdorff and~$n$-dimensional, and if~$\phi$ denotes the morphism~$X\to \gman$ induced by~$\bf f$, then the pre-image of the skeleton~$S_n$ of~$\gman$ under~$\phi$ has a piecewise-linear structure making~$\phi\inv(S_n)\to S_n$ a piecewise immersion ; this  is due to the author.

\medskip
In this article, we improve 1) and 2), and give a new proofs of both of them. Our proofs are based upon the {\em model theory of algebraically closed, non-trivially valued fields} and don't involve de Jong's alterations. 

\medskip
Let us quickly explain what we mean by improving 1) and 2). 

\medskip
$\bullet$ Concerning 1), we also prove that if~$x\in X$, there exists a compact analytic neighborhood~$U$ of~$x$, such that for every compact analytic neighborhood~$V$ of~$x$ in~$X$, the germs of polytopes~$(|{\bf f}|(V),|{\bf f}|(x))$ and~$(|{\bf f}|(U),|{\bf f}|(x))$ coincide. 

$\bullet$ Concerning 2), we prove that the piecewise linear structure on~$\phi\inv (S_n)$ is canonical, that is, doesn't depend on the map we choose to write it as a pre-image of the skeleton; we thus answer a question which was asked to us by Temkin. 

Moreover, we prove that the pre-image of the skeleton 'stabilizes after a finite, separable ground field extension', and that if~$\phi_1,\ldots, \phi_m$ are finitely many morphisms from~$X$ to~$\gman$, the union~$\bigcup \phi_j(S_n)$ also inherits a canonical piecewise-linear structure.
}

\section*{Introduction} 
Cet article est consacré, d'une manière générale, à l'étude des liens étroits qui existent entre la géométrie analytique ultramétrique (ici, au sens de Berkovich) et la géométrie linéaire par morceaux ; ces liens, dont le polygone de Newton constitue une première manifestation, se rencontrent essentiellement à travers deux types de phénomènes. 

\medskip
\begin{itemize}

\item[A)] L'image d'un espace analytique compact~$X$ par une application de «tropicalisation» est un polytope de dimension majorée par celle de~$X$ (\cite{loc2}, cor. 6.2.2 pour le cas localement algébrique, et \cite{imrsk}, \S 3.30 pour l'extension au cas général). La preuve de~\cite{imrsk} consiste essentiellement à algébriser la situation par des techniques standard (arguments de densité, lemme de Krasner...), pour se ramener au corollaire~6.2.2 de~\cite{loc2}, lui-même démontré par réduction, {\em via} les altérations de de Jong, au cas où $X$ possède un modèle polystable. 

\medskip
\item[B)] Certains sous-ensembles d'un espace de Berkovich héritent d'une structure d'espace linéaire par morceaux ; donnons quelques exemples. 

\medskip
\begin{itemize}
\item[B1)] Le cas archétypal est celui du «squelette standard»~$S_n$ de~$\gman$, défini comme l'ensemble~$\{\eta_{\bf r}\}_{{\bf r}\in (\RR^*_+)^n}$, où~$\eta_{\bf r}$ est la semi-norme~$\sum a_I{\bf T}^I\mapsto \max |a_I|\cdot {\bf r}^I$. 

\medskip
\item[B2)] Si~$\got X$ est un schéma formel polystable et quasi-compact sur~${\rm Spf}\,k\zero$, le «polytope d'incidence» de sa fibre spéciale s'identifie naturellement à un fermé~$S(\got X)$ de sa fibre générique~$\got X_\eta$ : c'est l'objet principal de l'article \cite{loc2} de Berkovich, voir notamment le théorème 5.1 de {\em loc. cit.}

\medskip
\item[B3)] Si~$X$ est un espace strictement~$k$-analytique de dimension~$n$ et si~$\phi : X\to \gman$ est un morphisme,~$\phi\inv (S_n)$ possède une unique structure linéaire par morceaux telle que~$\phi\inv(S_n)\to S_n$ soit linéaire par morceaux (et c'est alors une  immersion par morceaux). Ceci a été démontré par l'auteur dans \cite{imrsk}, par une démarche un peu analogue à celle évoquée à propos de l'assertion A) : algébrisation, réduction au cas de réduction polystable {\em via} de Jong, et utilisation de B2). 
\end{itemize}
\end{itemize}

\medskip
\noindent
{\bf Remarques} 

\medskip
$\bullet$ En 1984, Bieri et Groves avaient donné une preuve d'un avatar algébrique de l'assertion~A), fondée sur des calculs explicites à base de polygones de Newton  (\cite{bgrov}, th. 5.2) ; l'assertion A) dans le cas localement algébrique (autrement dit, le corollaire~6.2.2 de~\cite{loc2}), peut s'en déduire de manière à peu près formelle (sans recours aux altérations de de Jong). 

$\bullet$ Le théorème principal de~\cite{imrsk} est d'apparence plus générale que l'énoncé~B3), mais il s'y ramène en réalité : c'est l'objet des paragraphes~3.2--3.11 de~\cite{imrsk}.

\medskip
\noindent
{\bf Le contenu de cet article}
 
Dans ce texte, nous généralisons les énoncés A) et B3) tout en en donnant de nouvelles preuves : pour chacun d'eux nous utilisons, une fois la situation algébrisée, 
des outils de {\em théorie des modèles des corps valués}. Outils élémentaires pour notre généralisation de A), qui ne fait essentiellement appel qu'à l'élimination des quantificateurs ; outils plus avancés pour notre généralisation de B3), qui utilise des propriétés de finitude de certains des espaces de types que Hrushovski et Loeser ont récemment introduits pour étudier le type d'homotopie des espaces de Berkovich algébriques (\cite{hl}).

\medskip
Expliquons maintenant un peu plus avant de quoi il retourne.

\subsection*{Vocabulaire et conventions}

{\em Le point de vue multiplicatif.} Précisons pour commencer que nous suivons les conventions adoptées par Berkovich dans \cite{loc2} en matière de géométrie linéaire par morceaux. La principale originalité de son point de vue est l'usage de la {\em notation multiplicative}, qui est la plus naturelle lorsqu'on s'intéresse à des normes de fonctions holomorphes inversibles : elle dispense du choix totalement arbitraire d'une base de logarithmes, et évite une profusion de symboles «$\log$» et «$\exp$»~dans les formules et démonstrations. Dans ce qui suit, nous considérerons donc~$(\RR^*_+)^n$ comme un {\em espace vectoriel réel}, par le biais de sa structure de groupe abélien (loi interne) et de l'exponentiation coordonnée par coordonnée (loi externe) ; on peut alors parler de {\em polytopes} de~$(\RR^*_+)^n$, et modeler des {\em espaces linéaires par morceaux} abstraits sur ces derniers. Le lecteur trouvera  plus bas (\ref{defnotmult} {\em et sq.}) le rappel de toutes les définitions précises. 

\medskip
{\em Paramètres de définition des espaces analytiques et des polytopes.} On fixe un corps ultramétrique complet~$k$ et un sous-groupe~$\Gamma$ de~$(\RR^*_+)$ tel que~$\Gamma\cdot |k^*|$ soit non trivial. La classe des espaces~$k$-analytiques {\em~$\Gamma$-stricts} est celle des espaces analytiques modelés sur les lieux des zéros de fonctions holomorphes sur des polydisques dont le rayon appartient à~$\Gamma$ (pour une définition précise, {\em cf}. \ref{gammastr} {\em et sq.}).

On pose~$c=(\QQ,\sqrt{|k^*|\cdot \Gamma})$. Les  {\em~$c$-polytopes} de~$(\RR^*_+)^n$ sont les polytopes que l'on peut définir au moyen de formes affines «dont la partie linéaire est à coefficients rationnels et dont le terme constant appartient à~$\sqrt{|k^*|\cdot \Gamma}$» ; pour une définition précise, {\em cf}. \ref{defagtop}. Ils donnent naissance par recollement à la catégorie des {\em~$c$-espaces linéaires par morceaux} (\ref{defpl} {\em et sq.}). 

\medskip
{\em Algèbre commutative graduée.} Nous faisons un usage intensif dans ce texte du formalisme de l'algèbre commutative {\em graduée}, introduit et utilisé par Temkin dans \cite{tmk2} pour développer en toute généralité la théorie de la réduction des germes d'espaces~$k$-analytiques. Les rappels de base sur ces sujets sont faits aux  \ref{rappgrad} {\em et sq.}, ainsi qu'au \ref{germe}.  

\medskip
\subsection*{Tropicalisations d'un espace analytique}

Nous démontrons deux théorèmes sur le sujet. 

\medskip
{\bf Théorème \ref{tropicglob}}. {\em Soit~$X$ un espace~$k$-analytique compact et~$\Gamma$-strict de dimension~$d$, et soit~$(f_1,\ldots, f_n)$ une famille de fonctions holomorphes inversibles sur~$X$. Notons ~$|{\bf f}|$ l'application~$(|f_1|, \ldots, |f_n|) : X\to (\RR^*_+)^n$. Le compact~$|{\bf f}|(X)$ est alors un~$c$-polytope de~$(\RR^*_+)^n$ de dimension au plus~$d$ ; le compact~$|{\bf f}|(\partial X)$ est contenu dans un~$c$-polytope de~$(\RR^*_+)^n$ de dimension au plus~$d-1$, et est un~$c$-polytope de dimension au plus~$d-1$ si~$X$ est affinoïde.}

\medskip
{\bf Théorème \ref{theotroploc}}. {\em Soit~$X$ un espace~$k$-analytique compact et~$\Gamma$-strict et soient~$(f_1,\ldots, f_n)$ une famille de fonctions holomorphes inversibles sur~$X$. Soit~$x\in X$ et soit~$d$ le degré de transcendance sur le corps gradué~$\red k$ du sous-corps gradué~$\red k(\widetilde{f_1(x)}, \ldots, \widetilde{f_n(x)})$ de~$\widetilde{\hres(x)}$. Posons~$\xi=|{\bf f}|(x)$.

\medskip
\begin{itemize}

\item[1)] Il existe un voisinage~$k$-analytique compact et~$\Gamma$-strict~$U$ de~$x$ dans~$X$ possédant la propriété suivante : pour tout voisinage analytique compact~$V$ de~$x$ dans~$U$, les germes de polytopes~$(|{\bf f}|(U),\xi)$ et~$(|{\bf f}|(V),\xi)$ coïncident.

\medskip
De plus,~$|{\bf f}|(U)$ est de dimension~$\leq d$ en~$\xi$. Si~$(X,x)$ est sans bord,~$|{\bf f}|(U)$ est purement de dimension~$d$ en~$\xi$, et si de surcroît~$n=d$ alors~$|{\bf f}|(U)$ est un voisinage de~$\xi$ dans~$\rst n$.

\medskip
\item[2)] Supposons que~$X$ est affinoïde et que~$x\in \partial X$. Il existe un voisinage~$k$-analytique compact et~$\Gamma$-strict~$Y$ de~$x$ dans~$X$ possédant la propriété suivante : pour tout voisinage analytique compact~$Z$ de~$x$ dans~$Y$, les germes de polytopes~$(|{\bf f}|(Y\cap \partial X),\xi)$ et~$(|{\bf f}|(Z\cap \partial X),\xi)$ coïncident. 
\end{itemize}}

\medskip
Faisons quelques commentaires. 

\medskip
$\bullet$ Chacun de ces deux théorèmes comporte un énoncé relatif à~$X$, et un autre relatif à son bord~$\partial X$. Dans les deux cas, le second se déduit très simplement du premier, grâce au lemme \ref{lemmebord} qui assure que si~$X$ est affinoïde, son bord~$\partial X$ s'écrit comme une réunion finie d'espaces affinoïdes~$\Gamma$-stricts définis sur de «gros» corps et de dimension majorée par~$d-1$. 

$\bullet$ Hormis son assertion relative au bord (qui est nouvelle), le théorème \ref{tropicglob} est simplement l'énoncé A)  mentionné plus haut, et déjà connu. 

$\bullet$ Le théorème \ref{theotroploc} est entièrement nouveau. C'est la variante locale du théorème \ref{tropicglob} : il assure en gros que l'image d'un germe d'espace analytique~$\Gamma$-strict par une tropicalisation est un germe de~$c$-polytope. 

$\bullet$ Notons que le degré de transcendance de~$\red k(\widetilde{f_1(x)}, \ldots, \widetilde{f_n(x)})$ est exactement égal à~$n$ si et seulement si l'image de~$x$ sur~$\gman$ par le morphisme qu'induisent les~$f_i$ appartient au squelette standard~$S_n$ (\ref{transeteta}).

\medskip
Disons maintenant quelques mots des preuves de ces théorèmes. Elles reposent toutes deux sur l'élimination des quantificateurs dans la théorie des corps non trivialement valués algébriquement clos. Celle-ci permet tout d'abord d'établir un avatar du théorème \ref{tropicglob} portant sur une famille finie de fonctions inversibles sur une variété {\em algébrique} au-dessus d'un corps valué quelconque (non nécessairement de hauteur~$\leq 1$) : c'est le théorème \ref{propeq}. Il est peu ou prou équivalent au théorème de Bieri et Groves, démontré par des calculs explicites à base de polygones de Newton, que nous avons évoqué plus haut (\cite{bgrov}, th. 5.2)~-- le recours à l'élimination des quantificateurs rend la preuve nettement plus concise, mais moins effective.

\medskip
Pour prouver le théorème \ref{tropicglob}, on se ramène au théorème \ref{propeq} par les méthodes usuelles d'algébrisation déjà évoquées. Quant au théorème \ref{theotroploc}, il se démontre également à l'aide du théorème \ref{propeq}, mais de manière plus indirecte : on déduit tout d'abord de celui-ci un énoncé qui lui est apparenté et porte sur les espaces de Riemann-Zariski {\em gradués} (th. \ref{theoconetem}), que l'on combine ensuite avec la théorie de la réduction des germes d'espaces~$k$-analytiques pour aboutir à nos fins. 

\subsection*{Images réciproques du squelette standard de~$\gman$}

Avant d'énoncer notre théorème principal sur la question, nous commençons par introduire la notion de{\em~$c$-squelette} d'un espace~$k$-analytique~$\Gamma$-strict et topologiquement séparé~$X$ (\ref{defskel}). Sans entrer dans les détails, indiquons que si une partie localement fermée de~$X$ est un~$c$-squelette, elle possède une structure naturelle d'espace~$c$-linéaire par morceaux qui peut être décrite {\em purement en termes de l'espace analytique~$X$} (la définition précise de cette structure fait intervenir les domaines analytiques~$\Gamma$-stricts de~$X$ et les normes de fonctions inversibles sur ces derniers). 

\medskip
Le prototype du~$c$-squelette est le squelette standard~$S_n$ de~$\gman$ (\ref{exsquel}) ; nous aurons également besoin dans ce qui suit de considérer le squelette standard de~${\mathbb G}^{n,\rm an}_{m,F}$ lorsque~$F$ est une extension complète de~$k$, et nous le noterons~$S_{n,F}$. 

\medskip
Nous démontrons alors le théorème suivant. 

\medskip
{\bf Théorème \ref{imrecsquels}.} {\em Soit~$n$ un entier et soit~$X$ un espace~$k$-analytique topologiquement séparé,~$\Gamma$-strict et de dimension~$\leq n$. Soit~$(\phi_1,\ldots, \phi_m)$ une famille finie de morphismes de~$X$ vers~$\gman$.

\medskip
1) La réunion des~$\phi_j^{-1}(S_n)$ est un~$c$-squelette de~$X$, vide si~$\dim {} X<n$, et pour tout~$j$, l'application~$\phi_j^{-1}(S_n)\to S_n$ est une immersion par morceaux.

\medskip
2) Si~$X$ est compact il existe une extension finie séparable~$F_0$ de~$k$ telle que pour toute extension complète~$F$ de~$F_0$ , la flèche naturelle~$$\bigcup \phi_{j,F}\inv(S_{n,F})\to \bigcup \phi_{j,F_0}\inv(S_{n,F_0})$$ soit un homéomorphisme.}

\medskip
Faisons quelques commentaires. 

\medskip
$\bullet$ Supposons que~$j=1$ (on écrit alors~$\phi$ au lieu de~$\phi_1$) et que le groupe~$\Gamma$ est trivial. L'assertion 1) redonne l'énoncé B3) mentionné au début de l'introduction, mais est plus précise : elle assure par surcroît, en affirmant que~$\Sigma:=\phi\inv(S_n)$ est un~$c$-squelette, que la structure~$c$-linéaire par morceaux de~$\Sigma$ ne dépend pas du morphisme choisi pour le décrire comme image réciproque du squelette. Nous répondons ainsi à une question que nous avait posée Temkin. 

\medskip
$\bullet$ L'assertion 1) lorsque~$\Gamma\neq\{1\}$ ou lorsque~$m>1$ est nouvelle. 

\medskip
$\bullet$ L'assertion 2) est nouvelle.

\medskip
Donnons maintenant les grandes lignes la preuve du théorème \ref{imrecsquels}. Nous le prouvons tout d'abord dans le cas où~$m=1$, dans lequel nous nous plaçons donc pour le moment ; on écrit~$\phi$ au lieu de~$\phi_1$. La propriété à établir est locale, ce qui autorise à raisonner au voisinage d'un point~$x\in \phi\inv(S_n)$. Par les techniques standard d'algébrisation, on se ramène au cas où~$X$ est un voisinage de~$x$ dans l'analytification~$\sch X\an$ d'une~$k$-variété algébrique irréductible et lisse~$\sch X$, et où~$\phi$ provient d'une famille~$(f_1,\ldots, f_n)$ de fonctions inversibles (algébriques) sur~$\sch X$. 

\medskip
L'étape suivante consiste à montrer qu'il existe un voisinage affinoïde~$V$ de~$x$ dans~$X$ et d'une famille finie~$(g_1,\ldots, g_r)$ de fonctions analytiques inversibles sur~$V$ telle que la restriction de~$(|f_1|,\ldots, |f_n|, |g_1|, \ldots, |g_m|)$ à~$\phi_{|V}\inv(S_n)$ soit injective. On le déduit du théorème suivant. 

\medskip
{\bf Théorème \ref{propvalg}.} {\em Soit~$K$ un corps valué et soit ~$L$ une extension finie de~$K(T_1,\ldots, T_n)$. Il existe un sous-ensemble fini~$E$ de~$L^*$ tel que pour tout groupe abélien ordonné~$G$ contenant~$|K^*|$ et tout~${\bf g}\in G^n$, l'ensemble~$E$ sépare les prolongements de~$\eta_{\bf g}$ à~$L$, où~$\eta_{\bf g}$ est la valuation de Gauß~$\sum a_I{\bf T}^I\mapsto \max |a_I\cdot |{\bf g}^I$.}

\medskip
Pour prouver ce théorème, on raisonne par récurrence sur~$n$. Le passage de~$n$ à~$n+1$ repose de manière cruciale sur deux résultats de définissabilité en théorie des modèles des corps valués (dans le langage multisorte~$\mathscr L_{\rm val}$
que Haskell-Hrushovski-Macpherson ont introduit dans \cite{hhmelim}).

\medskip
$\bullet$ Le premier est  dû à Hrushovski et Loeser et porte sur les courbes relatives ; nous allons en dire quelques mots. 
Dans \cite{hl}, ils associent à tout morphisme~$\sch Z\to \sch Y$ entre variétés algébriques sur~$k$ un foncteur~$\widehat{\sch Z/ \sch Y}$ de la catégorie des extensions non trivialement valuées et algébriquement closes de~$k$ vers celle des ensembles ; le foncteur ~$\widehat{\sch Z/ \sch Y}$ est muni d'une application vers~$\sch Y$. On peut y penser comme un avatar modèle-théorique d'une fibration dont la base serait la variété {\em algébrique}~$\sch Y$ et les fibres les {\em analytifiés} à la Berkovich des fibres de~$\sch Z\to \sch Y$. 
Ainsi,~$\widehat{\Aff^1_{\sch Y}/ \sch Y}$ envoie~$F$ sur l'ensemble des couples~$(y,\phi)$ où~$y\in \sch Y(F)$ et où~$\phi$ est une \og semi-norme\fg~sur~$F[T]$ de la forme~$$\sum a_i(T-\lambda)^i \mapsto \max |a_i|\cdot r^i$$ pour un certain~$(\lambda, r)\in F\times |F|$ (on voit apparaître ici les valuations de Gauß, ce qui explique le rôle de ces espaces chapeautés dans la preuve du théorème~\ref{propvalg}). Hrushovski et Loeser montrent que~$\widehat{\sch Z/ \sch Y}$ est {\em prodéfinissable} en général, et {\em définissable} lorsque~$\sch Z\to \sch Y$ est de dimension relative majorée par 1 ; c'est cette dernière assertion que nous utilisons dans la preuve
du théorème~\ref{propvalg}. Elle découle elle-même du théorème de Riemann-Roch pour les courbes, et plus précisément de l'une de ses conséquences : si~$Z$ est une courbe algébrique projective, irréductible et lisse de genre~$g$ sur un corps algébriquement clos~$F$, le groupe~$F(Z)^*$ est engendré par les fonctions ayant au plus~$g+1$ pôles avec multiplicités. 

\medskip
$\bullet$ Le second est dû à Haskell, Hrushovski et Macpherson. Il assure que si~$F$ est un corps valué algébriquement clos, et si~$r$ est un élément d'un groupe ordonné~$G$ contenant~$|F^*|$, la structure engendrée par~$F$ et~$r$ est algébriquement close au sens de la théorie des modèles. Cela signifie la chose suivante. Donnons-nous une extension non trivialement valuée et algébriquement close~$F'$ de~$F$ telle que~$G\subset |(F')^*|$, ainsi qu'un produit fini~$\mathscr S$ de sortes. Soit~$\Phi$ une formule de~$\mathscr L_{\rm val}$ à paramètres dans~$(k,r)$, 
portant sur les éléments de~$\sch S$, 
et telle que 
l'ensemble~$E:=\{x\in \mathscr S(F'),\Phi(x)\}$ soit fini. Pour tout~$x_0\in E$ il existe alors une formule~$\Psi$ 
de~$\mathscr L_{\rm val}$ à paramètres dans~$(k,r)$ tel que~$\{x\in \mathscr S(F'),\Psi(x)\}=\{x_0\}.$

\medskip
Revenons à la preuve du théorème \ref{imrecsquels}. Une fois exhibés~$V$ et~$g_1,\ldots,g_r$ comme ci-dessus, on restreint un peu~$V$ de sorte qu'il satisfasse les conditions suivantes : 

\medskip
$\bullet$~$\phi$ induit un morphisme fini et plat de~$V$ sur un domaine affinoïde de~$\gman$ ;

$\bullet$ l'intersection~$\phi(V)\cap S_n$ est un pavé~$n$-dimensionnel. 

\medskip
Si~$y$ est un point de~$V$ qui n'est pas situé sur~$\phi\inv(S_n)$, on déduit du théorème \ref{theotroploc} qu'il existe un voisnage analytique compact de~$y$ dans~$V$ dont l'image par~$(|f_1|,\ldots, |f_n|)$ est de dimension~$\leq n-1$. Ce fait, joint au caractère ouvert des morphismes fini et plats et aux conditions imposées à~$V$, permet de donner une description de l'image de~$\phi_{|V}\inv(S_n)$ par ~$(|f_1|,\ldots, |f_n|, |g_1|, \ldots, |g_m|)$, que nous allons maintenant expliciter. 

L'image de~$V$ par~$(|f_1|,\ldots, |f_n|, |g_1|, \ldots, |g_m|)$ est un~$c$-polytope~$P$ de dimension~$\leq n$. Triangulons-le (on a simplement besoin de cellules convexes, il n'est pas nécessaire que ce soient des simplexes), et appelons~$Q$ la réunion des cellules fermée de dimension~$n$ de~$P$ en restriction auxquelles la projection sur les~$n$ premières coordonnées est injective. Le~$c$-polytope~$Q$ est alors précisément l'image de~$\phi_{|V}\inv(S_n)$ par ~$(|f_1|,\ldots, |f_n|, |g_1|, \ldots, |g_m|)$. 

On dispose ainsi d'un homéomorphisme entre~$\phi_{|V}\inv(S_n)$ et le~$c$-polytope~$Q$, homéomorphisme qui est construit de façon suffisamment «analytique» pour faire de~$\phi_{|V}\inv(S_n)$ un~$c$-squelette ({\em cf.} lemme \ref{testpolytanal}).  

\medskip
{\em Preuve de l'assertion 2) du théorème \ref{imrecsquels}, toujours dans le cas où~$m=1$.} Elle repose sur ce qui précède, et sur le fait suivant : le théorème \ref{propvalg} mentionné ci-dessus assure également que l'on peut, après passage à une extension finie séparable convenable de~$K$, exhiber un ensemble~$E$ qui sépare {\em universellement} ({\em i.e.} après n'importe quelle extension des scalaires) les prolongements des valuations de Gauß à~$L$. 

\medskip
{\em Preuve du  théorème \ref{imrecsquels} dans le cas où~$m$ est quelconque.} Le principe est le suivant : on fabrique un morphisme~$\psi : \Aff^N_X\to {\mathbb G}_{m,k}^{N+n,\rm an}$ pour un entier~$N$ convenable, et une section continue~$\sigma$ (de type «section de Shilov à rayon variable») du morphisme~$\Aff^N_X\to X$ qui identifie chacun des~$\phi_j\inv(S_n)$ à un sous-espace~$c$-linéaire par morceaux de~$\psi\inv(S_{n+N})$, lequel est un~$c$-squelette en vertu du cas~$m=1$ déjà traité. La réunion des~$\phi_j\inv(S_n)$ hérite ainsi d'une structure~$c$-linéaire par morceaux ; et là encore, la construction est suffisamment analytique pour que l'on puisse en déduire que~$\bigcup \phi_j\inv(S_n)$ est un~$c$-squelette. 

\subsection*{Remerciements} Je tiens à faire part de ma gratitude à Antoine Chambert-Loir. C'est en effet lors de discussions avec lui, au sujet d'un travail commun consistant à développer une théorie des formes différentielles et courants réels sur les espaces de Berkovich (\cite{fdceb}), que le projet du présent article a pris forme -- certains de ses résultats sont d'ailleurs utilisés dans \cite{fdceb}. 

\setcounter{section}{-1}

\section{Rappels, notations, préliminaires}

Afin d'obtenir un article raisonnablement auto-suffisant, 
nous avons choisi de faire figurer dans cette 
(longue) section un certain nombre de rappels, pour l'essentiel sans preuves, 
sur différentes théories que nous allons utiliser. On pourra la sauter en première lecture,  et ne s'y référer qu'en cas de besoin.  
Nous y présentons plus précisément : 

- la théorie de la réduction des germes de Temkin, et le formalisme de l'algèbre commutative graduée sur lequel elle se fonde (concernant ce dernier, nous avons démontré quelques assertions bien connues, faute de références dans la littérature) ; 

- nos conventions en matière d'espaces de Berkovich, quelques résultats de base sur le sujet dont nous nous servirons souvent, et la notion d'espace~$\Gamma$-strict, utile pour garder une trace des paramètres réels en jeu ; 

- la théorie des polytopes et espaces linéaires par morceaux «multiplicatifs», due à Berkovich. 

\medskip

\subsection*{Algèbre commutative graduée} 

Dans \cite{tmk2}, Temkin a introduit un certain nombre d'outils extrêmement efficaces pour l'étude locale des espaces analytiques. Ils reposent sur le formalisme de l'algèbre commutative {\em graduée}, dont nous allons rappeler les bases, sans démonstrations (celles-ci consistent essentiellement à retranscrire les preuves classiques, en rajoutant les adjectifs «homogène» ou «gradué» un peu partout). Le lecteur intéressé pourra consulter \cite{tmk2},  \cite{form} ou \cite{angel}. 

\medskip
On fixe pour ce paragraphe un groupe abélien divisible~$D$ noté multiplicativement ; en pratique, ce qui suit sera surtout utilisé lorsque ~$D=\RR^*_+$. Un {\em polyrayon} est une famille finie d'éléments de~$D$. 

\deux{rappgrad} Un {\em anneau~$D$-gradué} est un anneau (commutatif, unitaire)~$A$ muni d'une décomposition~$A=\bigoplus\limits_{r\in D}A^r$, telle que~$A^r\cdot A^s\subset A^{rs}$ pour tout~$(r,s)$ (attention : la graduation est {\em multiplicative}) ; on dit que~$A^r$ est l'ensemble des éléments homogènes de degré~$r$. Dans ce qui suit, nous dirons simplement «gradué» au lieu de «$D$-gradué». 

\medskip
Tout anneau~$A$ peut être vu comme un anneau gradué : il suffit de le munir de la {\em graduation triviale} pour laquelle~$A^1=A$ et~$A^r=\{0\}$ si~$r\neq 1$. 

\medskip
Un idéal d'un anneau gradué~$A$ est dit homogène s'il est engendré par des éléments homogènes ou, ce qui revient au même, s'il est somme directe de son intersection avec les~$A^r$. 

\deux{polgrad} Soit~$A$ un anneau gradué et soit~${\bf d}=(d_1,\ldots,d_n)$ un polyrayon. On note~$A[{\bf d}\inv{\bf T}]$ l'anneau gradué défini comme suit : l'anneau sous-jacent est l'anneau de polynômes~$A[{\bf T}]=A[T_1,\ldots, T_n]$, et l'ensemble des éléments homogènes de degré~$s$ de~$A[{\bf d}\inv{\bf T}]$ est l'ensemble des polynômes de la forme~$\sum a_I{\bf T}^I$ avec~$a_I\in A^{s{\bf d}^{-I}}$ pour tout~$I$. Lorsqu'on voudra évoquer le degré usuel d'un élément de~$A[{\bf d}\inv{\bf T}]$, on parlera de son degré {\em monomial}. 

\medskip
Si~$P$ est un élément homogène de degré~$s$ de~$A[{\bf d}\inv{\bf T}]$ et si~${\bf b}=(b_1,\ldots, b_n)$ est une famille d'éléments homogènes d'une~$A$-algèbre graduée~$B$, chaque~$b_i$ étant de degré~$d_i$, alors~$P({\bf b})$ est un élément homogène de degré~$s$ de~$B$.

\deux{corspgrad} Un {\em corps gradué} est un anneau gradué non nul dans lequel tout élément {\em homogène} non nul est inversible (en tant qu'anneau abstrait, un corps gradué  n'est pas nécessairement un corps). Si~$K$ est un corps gradué, un polyrayon~$\bf r$ sera dit {\em~$K$}-libre s'il constitue une famille libre de~$\QQ\otimes_{\ZZ}(D/D_0)$, où~$D_0$ désigne le groupe des degrés des éléments homogènes non nuls de~$K$. 

\medskip
Si~$\bf d$ est un polyrayon~$K$-libre alors~$K_{\bf d}:=K[{\bf d}\inv{\bf T}, {\bf d}{\bf T}\inv]$ est un corps gradué. 

\deux{polgrad} Soit~$K\hookrightarrow F$ une extension de corps gradués, soit~${\bf d}=(d_1,\ldots, d_n)$ un polyrayon, et soit~$(x_1,\ldots, x_n)$ une famille d'éléments homogènes de~$F$, chaque~$x_i$ étant de degré~$r_i$. On dit que les~$x_i$ sont algébriquement indépendants sur~$K$ si~$P(x_1,\ldots, x_n)\neq 0$ pour tout élément homogène non nul~$P\in K[{\bf T}/{\bf d}]$. Si~$x$ est un élément homogène de~$F$ on dira qu'il est transcendant sur~$K$ si la famille singleton~$\{x\}$ est algébriquement indépendante ; dans le cas contraire, on dira que~$x$ est algébrique. 

\medskip
On dispose dans ce cadre d'une théorie du degré de transcendance, analogue à celle bien connue dans le cas non gradué. 

\subsection*{Valuations graduées}

\deux{valgrad} Soit~$K$ un corps gradué. Une {\em valuation graduée sur~$K$} est une application~$\val .$ définie sur l'ensemble des éléments homogènes de~$AK$ et à valeurs dans un groupe abélien ordonné~$G$ (la notation est multiplicative) auquel on adjoint un plus petit élément absorbant~$0$, telle que :

\medskip

$\bullet$~$\val 1=1$,~$\val 0=0$ et~$\val {ab}=\val a \cdot \val b$ pour tout~$(a,b)$  (ce qui implique que~$\val a \neq 0$ dès que~$a\neq 0$) ; 
 
~$\bullet$  pour tout~$r>0$ et tout couple~$(a,b)$ d'éléments de~$K^r$ on a l'inégalité ultramétrique~$\val {a+b}\leq \max(\val a, \val b)$. 
 
 \medskip
Deux valuations graduées~$\val .: \bigcup K^r\to G \cup\{0\}$ et~$\val .' : \bigcup K^r\to G'\cup\{0\}$ sont dites équivalentes s'il existe une valuation graduée~$\val .'': \bigcup K^r\to G'' \cup\{0\}$ et deux morphismes strictement croissants~$i: G''\hookrightarrow G$ et~$i': G''\hookrightarrow G'$ tels que~$\val .=i\circ \val . ''$ et~$\val .'=i'\circ \val .''$.

\medskip
Soit~$\val .~$ une valuation graduée sur~$K$. L'ensemble des éléments homogènes~$x$ de~$K$ tels que~$\val x\leq 1$ est l'ensemble des éléments homogènes d'un sous-anneau gradué~$\sch O_{\val .}$ de~$K$, qui est appelé {\em l'anneau gradué de~$\val.$.} L'anneau ~$\sch O_{\val .}$ est un anneau gradué local : il a un et un seul idéal homogène maximal, qui est appelé l'idéal maximal  de~$\val .$. Il possède la propriété suivante : si~$x$ est un élément homogène non nul de~$K$ alors ou bien~$x\in \sch O_{\val .}$ ou bien~$x\inv\in \sch O_{\val .}$ ; cela équivaut à dire que~$\sch O_{\val .}$ est maximal pour la relation de domination entre sous-anneaux gradués locaux de~$K$ (si~$A$ et~$B$ sont deux sous-anneaux gradués locaux de~$K$, on dit que~$B$ domine~$A$ si~$A\subset B$ et si l'idéal homogène maximal de~$A$ est contenu dans celui de~$B$). 

\medskip
Inversement, toute anneau gradué de~$K$ qui est maximal pour la relation de domination est l'anneau d'une valuation graduée de~$K$, unique à équivalence près. On dispose ainsi d'une bijection entre l'ensemble des classes d'équivalences de valuations graduées sur~$K$ et celui des sous-anneaux gradués locaux de~$K$ maximaux pour la relation de domination. 

\medskip
La valuation gradué {\em triviale} sur~$K$ est celle dont l'anneau gradué est~$K$ tout entier ; elle envoie tout élément homogène non nul sur~$1$. 

\medskip
Si~$K\hookrightarrow F$ est une extension de corps gradués, le lemme de Zorn assure que toute valuation graduée~$\val .$ sur~$K$ s'étend en une valuation graduée~$\val .'$ sur~$F$ ; si~$\val .$ prend ses valeurs dans un groupe~$G$, on peut choisi ~$\val .'$ à valeurs dans le groupe ordonné~$G\otimes_{\ZZ}\QQ$, que l'on notera plus volontiers~$\sqrt G$.

\deux{valgausgrad} {\em Valuations de Gauß.} Soit~$K$ un corps gradué, muni d'une valuation graduée à valeurs dans un groupe ordonné~$G$, et soit~${\bf d}=(d_1,\ldots, d_n)$ un polyrayon. Soit~$K({\bf d}\inv{\bf T})$ le corps des fractions gradué de~$K[{\bf d}^{-1}{\bf T}, {\bf d}{\bf T}^{-1}]$, et soit~${\bf g}=(g_1,\ldots, g_n)$ un~$n$-uplet d'éléments de~$G$. On notera~$\eta_{K,{\bf d},{\bf g}}$ l'unique valuation graduée de~$K({\bf d}\inv{\bf T})$ qui envoie tout élément homogène~$\sum a_I{\bf T}^I$ de~$K[{\bf d}^{-1}{\bf T}]$ sur~$\max |a_I|{\bf g}^I$. Si~${\bf d}=(1,\ldots,1)$ on écrira simplement~$\eta_{K,\bf g}$, voire~$\eta_{\bf g}$ si le corps~$K$ est clairement indiqué par le contexte.

\deux{valtrvsalg} Soit~$K\hookrightarrow F$  une extension de corps gradués soit~$\val .~$ une valuation graduée sur~$K$. Soit~$t$ un élément homogène de~$F$. 

\medskip
{\em Supposons~$t$ algébrique sur~$K$}.  Pour toute valuation graduée~$\val .$ sur~$F$, il existe~$n>0$ tel que~$\val t^n\in \val {K}$. Notons un cas particulier important : si la restriction de~$\val .~$ à~$K$ est triviale, on a~$\val t=1$ si~$t\neq 0$. 

\medskip
{\em Supposons~$t$ transcendant sur~$K$.} Il existe alors une valuation graduée~$\val .$ sur~$F$, triviale sur~$K$, et telle que~$\val t>1$. En effet, soit~$d$ le degré de~$t$. Le sous-corps de~$F$ engendré par~$K$ et~$t$ s'identifie à~$K(d\inv T)$, et l'on peut alors modulo cette identification considérer un prolongement à~$F$ de la valuation gradué~$\eta_{K,d,r}$, où~$K$ est considéré comme trivialement valué et où~$r$ est un élément~$>1$ d'un groupe abélien ordonné quelconque, par exemple~$\RR^*_+$.

\deux{zr} Soit~$K$ un corps gradué. On note~$\PP_K$ l'ensemble l'espace de Riemann-Zariski gradué de~$K$, c'est-à-dire l'ensemble des classes d'équivalence de valuations graduées de~$F$. Pour toute partie~$X$ de~$\PP_K$ et tout ensemble~$E$ d'éléments homogènes de~$F$, on note~$X\{E\}$ (resp.~$X\{\{E\}\}$) le sous-ensemble de~$X$ formé des valuations~$\val .$ dont l'anneau gradué (resp. l'idéal homogène homogène maximal) contient~$E$ ou, si l'on préfère, telles que~$\val e\leq 1$ (resp.~$\val e <1$) pour tout~$e\in E$ On munit~$\PP_K$ de la topologie engendrée par les parties de la forme~$\PP_K\{E\}$ où~$E$ est {\em fini}.

Soient~$F$ et~$F'$ deux ensembles quelconques d'éléments homogènes de~$K$ ; posons~$X=\PP_K\{F\}\{\{F'\}\}$, et munissons-le de la topologie induite. Il est quasi-compact : cela résulte du fait qu'il est fermé dans~$\PP_K\{F\}$, et que ce dernier est quasi-compact d'après le lemme 2.1 de \cite{tmk2}. On qualifiera d'{\em affine} tout ouvert de~$X$ de la forme~$X\{E\}$ où~$E$ est un ensemble fini d'éléments homogènes de~$K$ ; remarquons que~$X$ est lui-même affine (prendre~$E=\emptyset$).

\medskip
Soit~$K_0$ un sous-corps gradué de~$K$ et soit~$\val .$ une valuation graduée de~$K_0$. Soit~$F$ (resp.~$F'$) l'ensemble des éléments homogènes de l'anneau gradué (resp. de l'idéal homogène) de~$\val.$. Le sous-ensemble~$\PP_K\{F\}\{\{F'\}\}$ de~$\PP_K$ est alors l'ensemble des valuations graduées qui prolongent la valuation~$\val .$ de~$K_0$. On le notera le plus souvent~$\PP_{K/(K_0,\val .)}$,  ou simplement~$\PP_{K/K_0}$ s'il n'y a pas d'ambiguïté sur la valuation de~$K_0$ ; lorsque celle-ci n'est pas précisée,~$K_0$ sera implicitement considéré comme muni de la valuation triviale. 

\subsection*{Valuations classiques et corps gradués résiduels} 

\deux{caspartvalclass} Si~$K$ est un corps classique, vu comme trivialement gradué,  une valuation graduée sur~$K$ n'est autre qu'une valuation (de Krull) usuelle. Les espaces~$\PP_K$ et~$\PP_{K/(K_0,\val .)}$ (où~$K_0$ est un sous-corps de~$K$ muni d'une valuation de Krull~$\val .$) sont les espaces de Riemann-Zariski usuels ; si~$K$ est une extension de type fini de~$K_0$, l'espace topologique~$\PP_{K/(K_0,\val .)}$ s'identifie à la limite projective des fibres spéciales de tous les modèles intègres, propres et plats du corps~$K$ sur~$\spec \sch O_{\val .}$.

\deux{redgrad} Soit~$K$ un corps classique, muni d'une valuation~$\val .$, prenant ses valeurs dans un groupe abélien ordonné divisible~$G$. Temkin associe à cette donnée un corps~$G$-gradué~$\red K$, appelé le {\em corps gradué résiduel} de~$\val .$ (ou, par un abus fréquent, de~$K$), et défini comme suit :~$$\red K=\bigoplus_{g\in G}\{z\in k, \val z\leq g\}/\{z\in k, \val z <g\}.$$ Remarquons que~$\red K$ ne change pas si l'on agrandit~$G$, et que~$\red K^1$ n'est autre que le corps résiduel traditionnel de~$K$. Toute extension~$K\hookrightarrow L$ de corps valués induit une extension~$\red K\hookrightarrow \red L$ de corps gradués. 

\medskip
Si~$g\in G$ et si~$\lambda$ est un élément de~$K$ tel que~$|\lambda|\leq g$ on notera~$\red \lambda^g$ l'image de~$\lambda$ dans~$\red K^g$ ; si de plus~$|\lambda|=g$ on écrira simplement~$\red \lambda$ ; si~$\lambda=0$ on pose~$\red \lambda=0$. 

\deux{interpdegtrres} Soit~$L$ une extension de~$K$, munie d'un prolongement de~$\val .$ noté encore~$\val .~$ ; quitte à agrandir~$G$, on suppose que~$|L^*|\subset G$. Soit~$d$ le degré de transcendance de~$\red L^1$ sur~$\red k^1$, et soit~$r$ la dimension du~$\QQ$-espace vectoriel~$\QQ\otimes_{\ZZ}(|L^*|/|K^*|)$. La somme~$d+r$ est alors égale au degré de transcendance de l'extension graduée~$\red k \hookrightarrow \red L$. 

\medskip
Pour le voir, fixons une base de transcendance~${\bf a}$ de~$\red L^1$ sur~$\red k^1$, et une partie~${\bf b}$ de~$L^*$ tels que~$|{\bf b}|$ soit une base de ~$\QQ\otimes_{\ZZ}(|L^*|/|K^*|)$. Nous allons montrer que~${\bf a}\cup\red{\bf b}$ est une base de transcendance de~$\red L~$ sur~$\red k$. 

\medskip
{\em La famille~${\bf a}\cup \red{\bf b}$ est algébriquement indépendante.} Soit~$h\in G$ et soit~$P$ un élément homogène de degré~$h$ de~$\red K[{\bf T}, |{\bf b}|\inv{\bf S}]$ tel que~$P({\bf a}, \red{\bf b})=0$. Nous allons montrer que~$P$ est nul. 

Les monômes qui constituent~$P$ sont tous homogènes de même degré~$h$ ; comme~$|{\bf b}|$ est une~$\QQ$-base de~$|L^*|$ modulo~$|K^*|$, ceci entraîne l'existence d'un monôme unitaire~$Q\in \red K[{\bf g}\inv{\bf S}]$  et d'un élément homogène~$R\in \red K[{\bf T}]$ tels que~$P=QR$. En utilisant encore la~$\QQ$-liberté de~$|{\bf b}|$ on voit que~$Q(\red{\bf b})\neq 0$ ; par conséquent,~$R({\bf a})=0$. On peut toujours écrire~$R=\alpha R^\sharp$, où~$\alpha$ est un élément homogène non nul de~$\red K$ et où~$R^\sharp$ est de degré 1, c'est-à-dire somme de monômes dont tous les coefficients appartiennent à~$\red K^1$. On a alors~$R^\sharp({\bf a})=0$ ; comme~$\bf a$ est une base de transcendance de~$\red L^1$ sur~$\red K^1$, il vient~$R^\sharp=0$ et finalement~$P=0$, ce qu'il fallait démontrer. 

\medskip
{\em Le corps gradué~$\red L$ est algébrique sur~$\red K({\bf a}\cup \red{\bf b})$.} Soit~$\lambda\in L^*$. Comme~$\bf g$ est une base de~$\QQ\otimes_{\ZZ}(|L^*|/|K^*|)$, il existe un entier~$n>0$, un multi-indice~$J$, un élément~$\alpha$ de~$K^*$ et un élément~$\beta\in L^*$ de valeur absolue égale à 1 tels que~$\lambda^n=\alpha\beta\cdot {\bf b}^J$. En réduisant, il vient~$$\red{\lambda}^n=\red \alpha\red\beta\cdot\red {\bf b}^J.$$ Comme~$\red \beta \in \red L^1$, il est algébrique sur~$\red K^1({\bf a})$ ; par conséquent,~$\red \lambda$ est algébrique sur~$\red K({\bf a}\cup \red{\bf b})$, ce qui achève la démonstration. 

\deux{transeteta} Soit~$K\hookrightarrow L$ une extension de corps valués, soit~$n$ un entier et soit~$L$ une extension de~$K$. Soient~$t_1,\ldots, t_n$ des éléments de~$L^*$. Les assertions suivantes sont équivalentes : 

\medskip
i) les~$t_i$ sont algébriquement indépendants sur~$K$, et la valuation de~$K({\bf t})$ est égale à~$\eta_{|{\bf t}|}$ ; 

ii) les éléments~$\red {t_i}$ de~$\red L$ sont algébriquement indépendants sur~$\red K$. 

\medskip
\medskip
En effet  i) est {\em fausse} si et seulement si il existe un polynôme~$P=\sum a_I{\bf T}^I$ tel que~$|P({\bf t})|<\max |a_I|\cdot |{\bf t}|^I$, c'est-à-dire encore s'il existe~$g\in |L^*|$ et une famille finie~$(a_I)~$ d'éléments de~$K$ vérifiant les conditions suivantes : ~$|a_I|\leq |{\bf t}|^{-I}g$ pour tout~$I$ ; il existe~$I$ pour lequel on a égalité ; et~$|\sum a_I{\bf t}^I|<g$. 

On peut reformuler cette condition en demandant qu'il existe~$g\in |L^*|$ et une famille finie~$(\alpha_I)$ d'éléments homogènes de~$\red K$ satisfaisant les conditions suivantes :~$\alpha_I$ est de degré~$|{\bf t}|^{-I}g$ ; les~$\alpha_I$ sont non tous nuls ; et l'élément~$\sum \alpha_I\red{\bf t}^I$ de~$\red K^g$ est nul (prendre~$\alpha_I=(\red {a_I})^{|{\bf t}| ^{-I}g}$). Mais cela signifie exactement qu'il existe~$g\in |L^*|$ et un élément~$Q$ non nul homogène de degré~$g$ dans~$\red K[|{\bf t}|\inv{\bf S}]$ tel que~$Q(\red{\bf t})=0$. 
 
 \medskip
 Autrement dit, i) est fausse si est seulement si ii) est fausse, ce qu'il fallait démontrer. 
 
 \deux{remresidgauss} {\em Remarque.} Supposons que les conditions équivalentes ci-dessus soient satisfaites, et soit~$\lambda=\sum a_I{\bf t}^I$ un élément non nul de~$K[t_1,\ldots, t_n]$. Soit~$\sch J$ l'ensemble des indices~$J$ tels que~$|a_J|\cdot {|\bf t}|^J=\max |a_I|\cdot |{\bf t}^J|=|\lambda|.$ On a~$$\red \lambda=\red{\sum_{J\in \sch J}a_J {\bf t}^J}=\sum_{J\in \sch J}\red {a_J}\red {\bf t}^J.$$ Il s'ensuit que~$\red{K({\bf t})}=\red K(\red{\bf t})$. Notons que ce dernier est par ailleurs, en vertu de nos hypothèses, naturellement isomorphe à~$\red K(|{\bf t}|\inv{\bf T})$. 
 
 \deux{corextgauss} Soit~$K\hookrightarrow L$ une extension finie de corps valués, et soit~$G$ un groupe abélien ordonné contenant~$|L^*|$. Soit~$n$ un entier, et soit~${\bf g}\in G^n$. Le seul prolongement de~$\eta_{K,{\bf g}}$ à~$L(T_1,\ldots, T_n)$ est~$\eta_{L, {\bf g}}$. En effet, fixons un tel prolongement. En vertu du \ref{transeteta} ci-dessus, les~$\red{T_i}$ sont algébriquement indépendants sur~$\red K$ ; l'extension~$K\hookrightarrow L$ étant finie,~$\red L$ est algébrique sur~$\red K$, et les~$\red{T_i}$ sont donc encore algébriquement indépendants sur~$\red K$. En utilisant à nouveau le \ref{transeteta}, on voit que la valuation de~$L(T_1,\ldots, T_n)$ est égale à~$\eta_{L,{\bf g}}$.

\subsection*{Géométrie analytique}

\deux{introgeanal} On fixe pour toute la suite du texte un corps ultramétrique complet~$k$ (sa valeur absolue peut être triviale) et un sous-groupe~$\Gamma$ de~$\RR^*_+$ tel que~$\Gamma\cdot |k^*|\neq\{1\}$ (autrement dit,~$\Gamma$ est non trivial si~$|k^*|=\{1\}$). Dans ce texte, la notion d'espace~$k$-analytique sera à prendre \emph{au sens de Berkovich} (\cite{brk1}, \cite{brk2}). Conformément aux conventions introduites plus haut,~$\red k$ désignera le corps {\em gradué} résiduel de~$k$. Si~$\sch X$ est une~$k$-variété algébrique, on notera~$\sch X\an$ son analytifiée. 

\deux{polir} Dans ce contexte, un {\em polyrayon} sera simplement une famille finie~$(r_1,\ldots, r_n)$ de réels strictement positifs ; on  dira {\em~$k$-libre} si les~$r_i$ sont linéairement indépendants dans~$\QQ\otimes_{\ZZ}(\RR^*_+/|k^*|)$. Si~${\bf r}=(r_1,\ldots, r_n)$ est un polyrayon~$k$-libre, l'algèbre normée~$k_{\bf r}:=k\{{\bf r}\inv {\bf T},{\bf r}{\bf T}\inv\}$ est un corps valué ; le foncteur d'extension des scalaires à~$k_{\bf r}$ sera simplement noté par un~$\bf r$ en indice.

On dira qu'un polyrayon~$k$-libre~$\bf r$ {\em déploie} une algèbre~$k$-affinoïde~$\sch A$ si  la valeur absolue de~$k_{\bf r}$ n'est pas triviale et si~$\sch A_{\bf r}$ est strictement~$k_{\bf r}$-affinoïde.

\deux{defdkx} Soit~$X$ un espace~$k$-analytique et soit~$x\in X$. On notera~$\hres x$ le corps résiduel complété de~$x$ ; lorsque~$X$ est bon, on notera~$\sch O_{X,x}$ l'anneau local de~$X$ en~$x$ (il aurait un sens dans le cas général, mais ne semble pertinent à considérer que dans le cas bon). 

On désignera par~$d_k(x)$ le degré de transcendance de l'extension de corps gradués~$\red k \hookrightarrow \red{\hres(x)}$ ; on peut également le définir comme la somme du degré de transcendance de l'extension de corps résiduels classiques~$\red k^1\hookrightarrow \red{\hres(x)}^1$ et de la dimension de~$\QQ\otimes_{\ZZ}(|\hres(x)^*|/|k^*|)$ (\ref{interpdegtrres}).

\medskip
On a l'égalité~$\dim{}X=\sup\limits_{x\in X}d_k(\hres(x))$. 

\deux{etaretsquel} Soit~${\bf r}=(r_1,\ldots, r_n)$ un polyrayon. La valuation~$\eta_{\bf r}$ du corps~$k(T_1,\ldots,T_n)$ induit un point de~$\gman$, encore noté~$\eta_{\bf r}$ ; l'application~${\bf r}\mapsto \eta_{\bf r}$ induit un homéomorphisme entre~$(\RR^*_+)^n$ et un fermé de~$\gman$ que l'on notera~$S_n$. 

\medskip
Soit~$x\in S_n$. On déduit de la remarque \ref{remresidgauss}, et de la densité de~$k(T_1,\ldots, T_n)$ dans~$\hres (x)$, que~$d_k(x)=n$. Il s'ensuit que si~$V$ est un domaine analytique de~$\gman$ et si~$Y$ est un fermé de Zariski de~$V$ contenant~$x$, la dimension de~$Y$ est égale à~$n$ (autrement dit,~$Y$ contient une composante connexe de l'espace normal~$V$, qui est purement de dimension~$n$). On déduit alors du corollaire 1.12 de \cite{flatn} que~$\sch O_{X,x}$ est artinien ; comme~$\gman$ est lisse en en particulier réduit, ~$\sch O_{X,x}$ est un corps.

\deux{gammastr} Soit~$\sch A$ une algèbre~$k$-affinoïde. Nous suivrons les conventions de \cite{flatn} : l'algèbre~$\sch A$ sera dite {\em~$\Gamma$-stricte} si elle peut s'écrire comme un quotient admissible de~$k\{r_1\inv T_1,\ldots, r_n\inv T_n\}$ pour une certaine famille~$(r_i)$ d'éléments de~$\Gamma$ ; un espace~$k$-affinoïde sera qualifié de~$\Gamma$-strict si son algèbre des fonctions analytiques est~$\Gamma$-stricte.

Si~$\sch A$ est une algèbre~$\Gamma$-stricte, il existe un polyrayon~$k$-libre~${\bf r}$ qui déploie~$\sch A$ et est constitué d'éléments de~$\Gamma$ ; et pour toute extension complète~$L$ de~$k$, l'algèbre~$L$-affinoïde~$\sch A_L$ est~$\Gamma$-stricte. 

Notons que lorsque~$|k^*|\neq\{1\}$ une algèbre~$k$-affinoïde est~$\{1\}$-stricte si et seulement si elle est strictement~$k$-affinoïde. 

\trois{racgamma} Pour qu'une algèbre~$k$-affinoïde~$\sch A$ soit~$\Gamma$-stricte, il suffit qu'elle s'écrive comme un quotient admissible de~$\sch A\simeq k\{r_1\inv T_1,\ldots, r_n\inv T_n\}$ pour une certaine famille~$(r_i)$ d'éléments de~$\sqrt{|k^*|\cdot \Gamma}$ (\cite{flatn}, 0.24.1). 

\trois{gammastrns} Un espace~$k$-affinoïde~$X$ est~$\Gamma$-strict si et seulement si~$$\sup\limits_{x\in X}|f(x)|\in \sqrt{|k^*|\cdot \Gamma}\cup\{0\}$$ pour toute fonction analytique~$f$ sur~$X$ (\cite{flatn}, 0.24.2).

\trois{gammastrgen} On dira qu'un espace~$k$-analytique~$X$ est~$\Gamma$-strict s'il satisfait les deux propriétés suivantes : 

\medskip
1)~$X$ possède un G-recouvrement par des espaces~$k$-affinoïdes~$\Gamma$-stricts ; 

2) l'intersection de deux domaines affinoïdes~$\Gamma$-stricts de~$X$ est G-recouverte par des domaines affinoïdes~$\Gamma$-stricts. 

\medskip
{\em Remarque.} La condition 2) est automatiquement satisfaite dès que~$X$ est séparé, auquel cas l'intersection de deux domaines affinoïdes~$\Gamma$-stricts de~$X$ est un domaine affinoïde~$\Gamma$-strict. 

\medskip
Cette définition est compatible avec la précédente lorsque~$X$ est affinoïde : cela résulte de la remarque précédente et de \ref{gammastrns}. 

\trois{classedense} En vertu de \ref{racgamma} et du fait que~$|k^*|.\Gamma\neq\{1\}$, tout point d'un espace affinoïde~$\Gamma$-strict a une base de voisinages affinoïdes~$\Gamma$-stricts ; par conséquent, tout ouvert d'un espace~$k$-analytique~$\Gamma$-strict est~$\Gamma$-strict. 

\trois{pleinfid} Soit~$f: Y\to X$ un morphisme entre deux espaces~$k$-analytiques~$\Gamma$-stricts. Si~$X'$ est un domaine analytique~$\Gamma$-strict de~$X$ alors~$f\inv(X')$ est~$\Gamma$-strict d'après un résultat de Temkin (\cite{tmk2}, cor. 4.10 ; \cf. aussi \cite{flatn}, lemma 0.29 iii) ). 

\trois{bonsdom} Si~$X$ est un bon espace~$k$-analytique~$\Gamma$-strict, tout point de~$X$ possède un voisinage affinoïde~$\Gamma$-strict dans~$X$ (\cite{flatn}, 0.30.2). 

\trois{gerrgrauertgam} Soit~$X$ un espace~$k$-affinoïde~$\Gamma$-strict et soit~$V$ un domaine affinoïde de~$X$. Nous dirons que~$V$ est {\em~$\Gamma$-rationnel} s'il peut être défini par une condition de la forme ~$$|f_1|\leq \lambda_1|g|\;{\rm et}\;|f_2|\leq \lambda_2|g|\;{\rm et}\; \ldots\;{\rm et}\; |f_n|\leq \lambda_n|g|,$$ où~$g$ et les~$f_i$ sont des fonctions analytiques sans zéro commun sur~$X$ et où les~$\lambda_i$ sont des réels appartenant à~$\Gamma$. 

Si~$V$ peut être défini par une telle condition avec les~$\lambda_i$ appartenant à~$\sqrt{|k^*|.\Gamma}$, il est~$\Gamma$-rationnel : on choisit un entier~$N>0$ tel que~$\lambda_i ^N=|a_i|.\gamma_i$ pour tout~$i$ avec~$a_i\in k^*$ et~$\gamma_i\in \Gamma$, et l'on décrit~$V$ par la condition~$$|f_1^N/a_1|\leq \gamma_1|g^N|\;{\rm et}\;|f_2^N/a_2|\leq \gamma_2|g^N|\;{\rm et}\; \ldots\;{\rm et}\; |f_n^N/a_n|\leq \gamma_n|g^N|.$$

\medskip
Un domaine affinoïde rationnel~$V$ de~$X$ est~$\Gamma$-strict si et seulement si il est~$\Gamma$-rationnel. En effet, si~$V$ est~$\Gamma$-rationnel il est~$\Gamma$-strict par définition. Réciproquement, supposons~$V$ rationnel. Il peut alors être défini par une condition de la forme ~$$|f_1|\leq \lambda_1|g|\;{\rm et}\;|f_2|\leq \lambda_2|g|\;{\rm et}\; \ldots\;{\rm et}\; |f_n|\leq \lambda_n|g|,$$ où~$g$ et les~$f_i$ sont des fonctions analytiques sans zéro commun sur~$X$ et où les~$\lambda_i$ sont des réels. Comme les~$f_i$ et~$g$ sont sans zéro commun sur~$X$, la fonction~$g$ ne s'annule pas sur~$V$. On peut dès lors, dans la définition de~$V$, remplacer chacun des~$\lambda_i$ par la borne supérieure de~$|f_i/g|$ sur~$V$, laquelle appartient à~$\sqrt{|k^*|\cdot \Gamma}$ parce que~$V$ est~$\Gamma$-strict. 

Le théorème de Gerritzen-Grauert classique (\cite{bgr}, \S 7.3.5, cor. 3 du th. 1) s'étend à ce contexte : tout domaine affinoïde~$\Gamma$-strict de~$X$ est réunion finie de domaines~$\Gamma$-rationnels. Pour le voir, on peut ou bien reprendre la preuve du lemme 2.4 de \cite{semialg} en choisissant un polyrayon~$\bf r$ constitué d'éléments de~$\Gamma$, ou bien reprendre la preuve du théorème 3.1 de \cite{tmk3} en utilisant la réduction~$\Gamma$-graduée des germes~$\Gamma$-stricts (au lieu de la réduction~$\RR^*_+$-graduée des germes généraux) introduite et étudiée dans \cite{flatn}, 0.29 {\em et sq.}

\medskip
Terminons cette section par un lemme bien connu dont nous rappelons la preuve pour la commodité du lecteur. 

\deux{algquasil} {\em Soit~$X$ un bon espace~$k$-analytique quasi-lisse~$\Gamma$-strict et soit~$x\in X$. Il existe un voisinage~$k$-affinoïde~$\Gamma$-strict~$V$ de~$x$ dans~$X$ et une~$k$-variété algébrique affine lisse~$\sch X$, de dimension égale à~$\dim x X$, tels que~$V$ s'identifie à un domaine affnoïde de~$\sch X\an$.}

\medskip
{\em Démonstration.} Posons~$d=\dim x X$. Comme~$X$ est quasi-lisse en~$x$, le point~$x$ possède un voisinage affinoïde~$\Gamma$-strict~$Y$ dans~$X$ qui s'identifie à un domaine affinoïde d'un espace lisse~$X'$ purement de dimension~$d$. Par définition de la lissité, il existe un morphisme étale d'un voisinage de~$x$ dans~$X'$ vers~$\Aff^{d,\rm an}_k$ ; soit~$\xi$ l'image de~$x$ sur~$\Aff^{d,\rm an}_k$ et soit~$\bf x$ l'image de~$\xi$ sur~$\Aff^d_k$. En vertu du lemme de Krasner et de la densité du corps~$\kappa({\bf x})$ dans~$\hres(\xi)$, il existe un morphisme étale~$\sch X\to \Aff^d_k$ dont l'image contient~$\bf x$, et un point~$\bf y$ appartenant à~$\sch X$ tel que~$\hres(x)\simeq \hres(\xi)\otimes_{\kappa({\bf x})}\kappa({\bf y})$ ; on peut toujours supposer~$\sch X$ affine, quitte à la restreindre autour de~$\bf y$. Il résulte dès lors du théorème 3.4.1 de \cite{brk2} que~$x$ possède dans~$X'$ un voisinage ouvert~$U$ qui se plonge dans~$\sch X\an$ (au-dessus de~$\Aff^{d,\rm an}_k$). On peut alors prendre pour~$V$ l'intersection de~$Y$ et d'un voisinage affinoïde~$\Gamma$-strict de~$x$ dans~$U$.~$\Box$ 

\subsection*{Germes d'espaces analytiques : la réduction de Temkin}

\deux{germe} À tout germe~$(X,x)$ d'espace~$k$-analytique, Temkin associe de manière fonctorielle dans \cite{tmk2} un espace topologique~$\red{(X,x)}$ connexe, non vide, et quasi-compact, muni d'une application continue~$\widetilde{(X,x)}\to \PP_{\red{\hres(x)}/\red k}$, qui est un homéomorphisme local. Le germe~$(X,x)$ est séparé (resp. bon, resp. sans bord) si et seulement si ~$\widetilde{(X,x)}$ s'identifie à un ouvert de~$\PP_{\red{\hres(x)}/\red k}$ (resp. à un ouvert affine de ~$\PP_{\red{\hres(x)}/\red k}$, resp. à ~$\PP_{\red{\hres(x)}/\red k}$ tout entier). La flèche~$(Y,x)\mapsto \red{(Y,x)}$ met en bijection l'ensemble des domaines analytiques de~$(X,x)$ et l'ensemble des ouverts quasi-compacts et non vides de~$\widetilde{(X,x)}$. 

\medskip
Si~$(X,x)$ est séparé (resp. bon), il est~$\Gamma$-strict si et seulement si~$\widetilde{(X,x)}$ s'écrit comme une réunion finie d'ouverts de la forme (resp. est de la forme)~$\PP_{\red{\hres(x)}/\red k}\{f_1,\ldots, f_n\}$, où chaque~$f_i$ est homogène de degré appartenant à~$\Gamma$ ({\em cf.} \cite{flatn}, lemma 0.28). 

\subsection*{Géométrie linéaire par morceaux}

\deux{defnotmult} Comme rappelé en introduction, nous adoptons de point de vue exposé par Berkovich dans \cite{loc2} en matière de géométrie linéaire par morceaux : on considère, un entier~$n$ étant donné,~$(\RR^*_+)^n$ comme un {\em espace vectoriel réel}, par le biais de sa structure de groupe abélien (loi interne) et de l'exponentiation coordonnée par coordonnée (loi externe). Si~$(a_1,\ldots, a_n)$ sont des nombres réels et si~$g\in \RR^*_+$, l'application de~$(\RR^*_+)^n$ vers~$\RR^*_+$ qui envoie~$(t_1,\ldots, t_n)$ sur~$g\prod t_i^{a_i}$ est alors affine, et linéaire si~$g=1$. 

\medskip
On fixe un couple~$c=(A,G)$, où~$A$ est un sous-anneau de~$\RR$ et où~$G$ est un sous-$A$-module  de~$\RR^*_+$. On note~$\mathsf {Aff}_c(\rst n)$ l'ensemble des applications affines de~$\rst n$ vers~$\RR^*_+$ de la forme~$$(t_1,\ldots, t_n)\mapsto g\prod t_i^{a_i},$$ où les~$a_i$ appartiennent à~$A$ et où~$g\in G$.

\trois{defagtop} Un {\em~$c$-polytope} de~$\rst n$ est une partie compacte de~$\rst n$ définie par une condition de la forme~$$\bigvee_i \bigwedge_j \phi_{i,j}\leq 1,$$ où les~$\phi_{i,j}$ appartiennent à~$\mathsf {Aff}_c(\rst n)$. 

\trois{dimpol} Soit~$P$ un~$c$-polytope convexe et non vide de~$\rst n$. La {\em dimension} de~$P$ est par définition celle du sous-espace affine qu'il engendre. 

\medskip
Si~$P$ est un~$c$-polytope quelconque de~$\rst n$, sa dimension est le maximum des dimensions des~$c$-polytopes convexes qu'il contient.

\trois{defpl} Si~$P$ est un~$c$-polytope de~$\rst n$, on notera~$\lgk Pc$ l'ensemble des applications continues~$\phi$ de~$P$ dans~$\RR^*_+$ possédant la propriété suivantes : il existe un recouvrement fini~$(P_i)$ de~$P$ par des~$c$-polytopes, et pour tout~$i$ une application~$\psi_i\in \mathsf {Aff}_{(A,G)}(\rst n)$, tels que~$\phi_{|P_i}=\psi_{i,|P_i}$ quel que soit~$i$.

\medskip
Si~$Q\subset \rst m$ est un~$c$-polytope on notera~$\pl cPQ$ l'ensemble des applications continues de~$P$ dans~$Q$ qui sont de la forme~$(\phi_1,\ldots, \phi_m)$ où~$\phi_j\in \lgk P  \Gamma$ pour tout~$j$. 

\medskip
La classe~$\pl c..$ d'applications continues entre~$c$-polytopes est stable par composition. 

\deux{structpolyt} Soit~$X$ un espace topologique compact. 

\trois{defstrpol} Une {\em structure~$c$-polytopale} sur~$X$ est la donnée d'un ensemble~$\lgk X c$ de fonctions continues de~$X$ vers~$\RR^*_+$ telles qu'il existe un entier~$n,$ un~$c$-polytope~$P\subset \rst n$, et un homéomorphisme~$X\simeq P$ modulo lequel~$\lgk X c$ s'identifie à~$\lgk P c$. 

On dira qu'un tel homéomorphisme~$X\simeq P\subset \rst n$ est une {\em présentation} de~$(X,\lgk X c)$. 

\trois{souspol} Supposons~$X$ muni d'une structure polytopale~$c$-stricte~$\lgk X c$. On dira qu'un sous-ensemble compact~$Y$ de~$X$ est un {\em~$c$-polytope} de~$(X,\lgk X c)$ s'il existe une présentation~$X\simeq P\subset \rst n$ de~$(X,\lgk X c)$ qui identifie~$Y$ à un~$c$-polytope de~$\rst n$ (évidemment contenu dans~$P$) ; c'est alors le cas pour {\em toute} présentation de~$(X,\lgk X c)$. 

\medskip
Soit~$Y$ un~$c$-polytope de~$X$ et soit~$\iota : X\simeq P\subset \rst n$ une présentation de~$(X,\lgk X c)$. L'ensemble des fonctions continues sur~$Y$ de la forme~$\phi\circ \iota_{|Y}$, où~$\phi\in \lgk {\iota(Y)} c$, ne dépend pas de la présentation choisie et sera noté~$\lgk Y c$ ; il définit une structure~$c$-polytopale sur~$Y$. 

Pour toute présentation ~$\iota : X\simeq P\subset \rst n$ de~$(X,\lgk X c)$, la restriction~$\iota_{|Y}: Y\simeq \iota(Y)\subset \rst n$ est une présentation de~$(Y,\lgk Y c)$. 

\medskip
Si~$(X_i)$ est un recouvrement fini de~$X$ par des~$c$-polytopes, une fonction continue~$\phi : X\to \RR^*_+$ appartient à~$\lgk X c$ si et seulement si sa restriction à~$X_i$ appartient à~$\lgk {X_i} c$ pour tout~$i$. 

\trois{morpol} Soit~$(Z,\lgk Z \Gamma)$ un espace topologique compact muni d'une structure~$c$-polytopale. Soient~$\iota : X\simeq P\subset \rst n$ et~$\upsilon : Z\simeq q \subset \rst m$ des présentations respectives de~$(X,\lgk X c)$ et~$(Z,\lgk Z c)$ et soit~$\phi$ une application continue de~$X$ vers~$Z$. 

\medskip
Les assertions suivantes sont équivalentes : 

\medskip
i) l'application~$P\to Q$ déduite de~$\phi$ {\em via} les homéomorphismes~$\iota$ et~$\upsilon$ appartient à~$\pl \Gamma c$ ; 

ii) pour toute~$\psi\in \lgk Z c~$, la composée~$\psi\circ \phi$ appartient à~$\lgk X c$. 

\medskip
En particulier, la validité de i) ne dépend pas du choix des présentations. 

\medskip
On notera~$\pl c XZ$ l'ensemble des applications continues de~$X$ vers~$Z$ qui satisfont ces propriétés équivalentes. 

\medskip
Si~$Y$ est un~$c$-polytope de~$X$, l'inclusion~$Y\subset X$ appartient à~$\pl c YX$. 

\medskip
La classe~$\pl c ..$ d'applications continues entre espaces topologiques munis d'une structure~$c$-polytopaleest stable par composition. 

\deux{cartepol} Soit~$X$ un espace topologique séparé et localement compact. 

\trois{defgrec} Soit~$Y$ une partie de~$X$ et soit~$(Y_i)$ une famille de sous-ensembles de~$Y$. On dit que la famille~$(Y_i)$ est  un {\em G-recouvrement} de~$Y$ si  tout point~$y$ de~$Y$ possède un voisinage dans~$Y$ de la forme~$\bigcup\limits_{i\in I}Y_i$ où~$I$ est un ensemble fini d'indices et où~$y\in \bigcap\limits_{i\in I}Y_i$. 

\trois{defcartepol} Une {\em carte~$c$-polytopale} sur~$X$ est la donnée d'une partie compacte~$Y$ de~$X$ et d'une structure~$c$-polytopale~$\lgk Y  c$ sur~$Y$. Deux cartes~$c$-polytopales~$(Y,\lgk Y c)$ et~$(Z, \lgk Z c)$ sont dite {\em compatibles} si~$Y\cap Z$ est un~$c$-polytope de~$Y$ aussi bien que de~$Z$, et si les structures~$c$-polytopales induites sur~$Y\cap Z$ par celle de~$Y$ et celle de~$Z$ coïncident.

\trois{atlpolyt} Un {\em atlas~$c$-polytopal}~$\got A$ sur~$X$ est une famille~$(Y_i, \lgk {Y_i} c)$ de cartes~$c$-polytopales deux à deux compatibles sur~$X$, telles que~$(Y_i)$ constitue un G-recouvrement de~$X$.

\medskip
Une carte~$c$-polytopale~$(Z,\lgk Z c)$ de~$X$ sera dite {\em compatible} avec~$\got A$ si elle l'est avec chacune des cartes de~$\got A$. Si~$(Z,\lgk Z c)$ est une carte polytopale~$c$-stricte de~$X$ compatible avec~$\got A$ alors~$\lgk Z c$ est nécessairement l'ensemble des applications continues~$\phi : Z\to \RR^*_+$ telles que~$\phi_{|Z\cap Y}\in \lgk {Z\cap Y} c$ pour toute carte~$(Y,\lgk Y c)$ de~$\got A$. Il y a donc au plus une structure~$c$-polytopale compatible avec~$\got A$ sur un compact~$Z$ de~$X$. 

Il s'ensuit que deux cartes~$c$-polytopales compatibles avec~$\got A$ sont compatibles entre elles. 

\medskip
Si~$\got A$ et~$\got B$ sont deux atlas~$c$-polytopaux  sur~$X$, ils seront dits {\em équivalents} si toute carte de~$\got A$ est compatible avec toute carte de~$\got B$. Si c'est le cas, une carte~$c$-polytopale de~$X$ est compatible avec~$\got A$ si et seulement si elle est compatible avec~$\got B$.

\deux{defpl} On définit comme suit la catégorie des {\em espaces~$c$-linéaires par morceaux}.

\trois{obpl} {\em Les objets.}  Un espace~$c$-linéaire par morceaux est un espace topologique séparé et localement compact~$X$ muni d'une classe d'équivalence~$\sch E$ d'atlas~$c$-polytopaux. 

\trois{convpl} {\em Quelques conventions.} Soit~$(X,\sch E)$ un espace~$c$-linéaire par morceaux (en pratique,~$\sch E$ sera bien entendu le plus souvent omis des notations, si cela ne prête pas à ambiguïté). 

\medskip
$\bullet$ Un atlas~$c$-polytopal sur~$X$ sera toujours, sauf mention expresse du contraire, supposé appartenir à~$\sch E$.

$\bullet$ Si~$Z$ est un compact de~$X$ qui possède une structure~$c$-polytopale compatible avec l'un des atlas de~$\sch E$, cette structure est unique, et est compatible avec tous les atlas de~$\sch E$. Nous dirons qu'un tel~$Z$ est un~$c$-polytope de~$X$. L'ensemble des~$c$-polytopes de~$X$ est le plus grand atlas~$c$-polytopal sur~$X$. 

\trois{moral} {\em Les flèches.} Soient~$Y$ et~$X$ deux espaces~$c$-linéaires par morceaux. Un morphisme de~$Y$ vers~$X$ est une application continue~$\phi : Y\to X$ telle qu'il existe un atlas~$c$-polytopal$\got A$ sur~$X$ et un atlas polytopal~$c$-polytopal~$\got B$ sur~$Y$ satisfaisant la condition suivante : pour tout~$Z\in \got B$, il existe~$T\in \got A$ contenant~$\phi(Z)$ et tel que~$\phi_{|Z}: Z\to T$ appartienne à~$\pl c ZT$. 

\deux{exagpl} {\bf Exemples.} Si~$X$ est un espace topologique compact muni d'une structure~$c$-polytopale, l'atlas~$\{X\}$ fait de~$X$ un espace~$c$-linéaire par morceaux. Ses~$c$-polytopes au sens de \ref{convpl} coïncident avec ceux définis au \ref{souspol}. 

L'espace~$\rst n$, muni de l'atlas de {\em tous} ses~$c$-polytopes (au sens de \ref{defagtop}), est un espace~$c$-linéaire par morceaux. Ses~$c$-polytopes au sens de \ref{convpl} sont les mêmes que ceux définis au \ref{defagtop}. 
 
\deux{notat} Si~$Y$ et~$X$ sont deux espaces~$c$-linéaires par morceaux, on notera~$\pl c YX$ l'ensemble des morphismes de~$Y$ vers~$X$ ; on qualifiera parfois ces morphsimes d'applications~$c$-linéaires par morceaux. Si~$X=\RR^*_+$, on écrira~$\lgk Y c$ au lieu de~$\pl c YX$ ; on peut alternativement définir~$\lgk Y c$ comme l'ensemble des fonctions continues~$\phi : X\to \RR^*_+$ telles que~$\phi_{|Z}\in \lgk Z c$ pour tout~$c$-polytope~$Z$ de~$Y$. 

Ces notations sont bien entendu compatibles avec elles précédemment introduites dans des cas particuliers. 

\deux{dompl} Soit~$X$ un espace~$c$-linéaire par morceaux et soit~$Y$ une partie de~$X$. On dit que~$Y$ est un {\em sous-espace~$c$-linéaire par morceaux} de~$X$ s'il est G-recouvert par les~$c$-polytopes de~$X$ qu'il contient. Tout sous-espace~$c$ linéaire par morceaux de~$X$ est localement fermé. 

\trois{exdompl} Les ouverts de~$X$ et les~$c$-polytopes de~$X$ sont des sous-espaces~$c$-linéaires par morceaux de~$X$. 

\trois{ensembstruct} Tout sous-espace~$c$-linéaire par morceaux~$Y$ de~$X$ hérite d'une structure naturelle d'espace~$c$-linéaire par morceaux, et l'inclusion de~$Y$ dans~$X$ est une application~$c$-linéaire par morceaux.  

\medskip
Soit~$X'$ un espace~$c$-linéaire par morceaux, soit~$\phi\in \pl c {X'}X$ et soit~$Y$ un sous-espace~$c$-linéaire par morceaux de~$X$. L'image réciproque~$Y':=\phi\inv(Y)$ est alors un sous-espace~$c$-linéaire par morceaux  de~$X'$, et~$Y'\to Y$ appartient à~$\pl c {Y'} Y$.

\trois{gtoppl} On définit~$X_{\grot}$ comme le site dont les objets sont les sous-ensembles~$c$-linéaires par morceaux, les morphismes les inclusions, et les recouvrements les G-recouvrements. Cette topologie de Grothendieck (que l'on appelle parfois la G-topologie) est plus fine que la topologie usuelle : cela résulte de \ref{exdompl}, et du fait que si~$(U_i)$ est une famille d'ouverts de~$X$, elle constitue un G-recouvrement de~$\bigcup U_i$. 

\trois{immerspl} Si~$X'$ est un espace~$c$-linéaire par morceaux, on dira qu'une application~$\phi \in \pl c {X'}X$ est une {\em immersion} si elle identifie~$X'$ à un sous-espace~$c$-linéaire par morceaux de~$X$ ; on dira que~$\phi$ est G-localement une immersion s'il existe un G-recouvrement~$(X_i)$ de~$X'$ par des sous-espaces~$c$-linéaires par morceaux tels que~$\phi_{|X_i}$ soit une immersion pour tout~$i$. 

\deux{changeant} Soit~$A'$ un sous-anneau de~$\RR$ contenant~$A$, et soit~$G'$ un sous-$A'$-module de~$\RR^*_+$ contenant~$G$. Posons~$c'=(A',G')$. Tout~$c$-polytope~$P$ de~$\rst n$ est {\em a fortiori} un~$c'$-polytope,  et~$\lgk P c\subset \lgk P {c'}$. Ceci permet de définir un foncteur «d'extension du domaine des paramètres autorisés»~de la catégorie des espaces~$c$-linéaires par morceaux vers celle des espaces~$c'$-linéaires par morceaux ; il ne modifie pas les espaces topologiques en jeu.

\subsection*{Théorie des modèles}

Nous utiliserons dans ce texte le langage~$\sch L_{\rm val}$ des corps valués, et plus précisément sa variante multisorte introduite par Haskell, Hrushovski et Macpherson dans \cite{hhmelim} afin d'éliminer les imaginaires. Ce n'est toutefois qu'à la section \ref{Gau}, consacrée aux prolongements des valuations de Gauß, que nous nous servirons vraiment de celle-ci -- pas directement, mais par le biais des espaces de types de Hrushovski et Loeser (\cite{hl}). En ce qui concerne la section \ref{elimacvf}, elle ne met en jeu que le langage plus classique à trois sortes (corps valué, corps résiduel, groupe des valeurs).

\deux{notcategor} Introduisons quelques notations pour un certain nombre de catégories dont nous aurons abondamment besoin ; rappelons que les les groupes abéliens ordonnés et les valuations 
sont notés multiplicativement. 

\trois{catcetd} On désigne par~$\mathsf C$ la catégorie des corps algébriquement clos non trivialement valués, et par~$\mathsf D$ celle des groupes abéliens divisibles ordonnés non triviaux. 

\trois{dgetck} Soit~$G$ un groupe abélien ordonné quelconque. On note~$\mathsf D_G$ la catégorie des groupes abéliens divisibles ordonnés non triviaux~$H$ munis d'un plongement (croissant)~$G\hookrightarrow H$. 

Soit~$K$ un corps valué. Supposons-le muni d'un plongement croissant~$|K^*|\hookrightarrow G$ ; on note alors~$\mathsf C_{K,G}$ la catégorie des corps algébriquement clos non trivialement valués~$L$ munis d'un plongement isométrique~$K\hookrightarrow L$ et d'un~$|K^*|$-plongement croissant~$G\hookrightarrow |L^*|$. Lorsque~$G=|K^*|$, on écrira~$\mathsf C_K$ au lieu de~$\mathsf C_{K,G}$ ; la catégorie~$\mathsf C_K$ est simplement celle des extensions non trivialement valuées algébriquement closes de~$K$.

\trois{deffonctdef} Soit~$\Phi$ une formule du langage~$\sch L_{\rm val}$ sans paramètres (resp. à paramètres dans~$(K,G)$), dont les variables libres
vivent dans un produit donné~$\sch S$ de sortes. 
La modèle-complétude de la théorie des corps  algébriquement clos non trivialement valués assure que~$$F\mapsto \{{\bf x}\in \sch S(F),\Phi({\bf x})\}$$ est de manière naturelle un 
foncteur de~$\mathsf C$ (resp.~$\mathsf C_{K,G}$) vers~$\mathsf {Ens}$. On dira qu'un tel foncteur est 
un {\em sous-foncteur définissable} de~$\sch S$  (resp.~$\sch S_{|\mathsf C_{K,G}}$) vers~$\mathsf {Ens}$. 
Soient~$\sch S$ et~$\sch S'$ deux produits finis de sortes, soit~$\mathsf X$ un sous-foncteur 
définissable de~$\sch S$ (resp.~$\sch S_{|\mathsf C_{K,G}}$), et soit 
$\mathsf X'$ un sous-foncteur 
définissable de~$\sch S'$ (resp.~$\sch S'_{|\mathsf C_{K,G}}$). On dit qu'une transformation
naturelle de~$\mathsf X$ vers~$\mathsf X'$ est définissable si son graphe est un sous-foncteur
définissable de~$\sch S\times \sch S'$ (resp.~$(\sch S\times \sch S')_{|\mathsf C_{(K,G)}}$) ; l'image
d'une transformation naturelle définissable est un sous-foncteur définissable du but. 

\trois{fontdefabstr} Il existe une version «abstraite», c'est-à-dire non plongée, de la notion de sous-foncteur définissable. 
Un foncteur~$\mathsf X$
de~$\mathsf C$ (resp.~$\mathsf C_{K,G}$) vers~$\mathsf{Ens}$ est dit {\em définissable} s'il
existe un produit fini de sortes~$\sch S$ et un sous-foncteur définissable~$\mathsf X_0$ 
de~$\sch S$  (resp.~$\sch S_{|\mathsf C_{K,G}}$) tel que
$\mathsf X\simeq \mathsf X_0$. 

\medskip
{\em A priori},~$\mathsf X_0$ n'a pas de raison 
d'être unique à unique isomorphisme définissable près. En pratique, il le sera : lorsqu'on dira
qu'un foncteur~$\mathsf X$ de~$\mathsf C$ (resp.~$\mathsf C_{K,G}$)
vers~$\mathsf{Ens}$ est définissable, il
sera toujours sous-entendu qu'on sait décrire non seulement
$\mathsf X(F)$ pour tout~$F$ appartenant à~$\mathsf C$ ou
$\mathsf C_{K,G}$, mais également les 
«familles définissables plongées d'objets paramétrés par~$\mathsf X$»~, ce qui assurera
la canonicité de~$\mathsf X_0$. Pour donner un sens rigoureux à cette 
dernière phrase nous renvoyons le lecteur
à la section 2 de \cite{moder}. 

\medskip
Par exemple, si~$\sch X$ est un~$K$-schéma de type fini, 
il définit de manière naturelle un foncteur définissable de~$\mathsf C_K$ vers~$\mathsf{Ens}$. 
Et un sous-foncteur de~$\sch X|_{\mathsf C_{K,G}}$ 
est définissable si et seulement si il est 
défini Zariski-localement 
sur~$\sch X$ par une combinaison booléenne d'inégalités de la forme
$|f|\Join \gamma |g|$, où~$f$ et~$g$ sont des fonctions régulières, où~$\gamma\in G$ 
et où~$\Join\in \{<,>,\leq,\geq\}$ : c'est une
conséquence de l'élimination des quantificateurs
pour la théorie des corps non trivialement valués algébriquement clos, 
dans le langage classique des corps valués à trois sortes (le corps valué, le corps 
résiduel et le groupe des valeurs). 

\trois{defdoag} On dispose de définitions analogues sur les catégories~$\mathsf D$ et~$\mathsf D_G$, en se restreignant au langage des groupes abéliens ordonnés. 

\deux{gencorpsval} Soit~$\GG$ le foncteur d'oubli de~$\mathsf D$ vers~$\mathsf{Ens}$ ; pour tout groupe abélien ordonné~$G$, on note~$\GG_G$ la restriction de~$\GG$ à~$\mathsf D_G$. 

\deux{fonctdef} Fixons un groupe abélien ordonné~$G$. Il résulte de l'élimination des quantificateurs dans la théorie des groupes abéliens ordonnés divisibles non triviaux qu'un sous-foncteur de~$\GG_G^n$ est définissable si et seulement si il peut s'écrire~$$H\mapsto\left \{{\bf h}\in H^n, \bigvee _{i\in I} \bigwedge_{j\in I_i} \phi_{i,j}({\bf h})  \Join_{ij} 1 \right\}$$ où les ensembles~$I$ et~$J_i$ sont tous finis, où chacune des~$\phi_{i,j}$ est de la forme~$at_1^{e_1}\ldots,t_n^{e_n}$ avec~$a\in G$ et~$e_\ell\in \ZZ$ pour tout~$\ell$, et où~$\Join_{i,j}\in \{<,\leq \}$ pour tout~$(i,j)$. 

\medskip
Si~$\sch D$ est un sous-foncteur définissable de~$\GG_G^n$ : 

\medskip
i) il est entièrement déterminé par sa valeur sur n'importe quel~$H\in \mathsf D_ G$ ;

\medskip
ii) l'expression~$\sch D(H)$ garde un sens pour {\em tout} groupe abélien ordonné~$H$ muni d'un plongement~$G\hookrightarrow H$.

\trois{fonctdefferm} On dira qu'un sous-foncteur définissable de~$\GG_G^n$ est {\em fermé} s'il admet une description comme ci-dessus avec~$\Join_{i,j}$ égal à~$\leq$ pour tout~$(i,j)$.

\trois{pavesval} Si~$r$ et~$R$ sont des éléments de~$G$ avec~$r<R$ on note~$]r;R[$ le sous-foncteur définissable de~$\GG_G$ qui envoie~$H$ sur~$\{h\in H, r<h<R\}$. Un {\em pavé de dimension~$d$ de~$\GG_G^n$} est un sous-foncteur définissable de~$\GG_G^n$ de la forme~$$M\left(\prod_{i=1}^d ]r_i;R_i[\times \{(g_{d+1},\ldots, g_n)\}\right)$$ où les~$r_i$, les~$R_i$ et les~$g_j$ sont des éléments fixés de~$G$, avec~$r_i<R_i$ pour tout~$i$, et où~$M\in {\rm GL}_n(\QQ)$.

\trois{defdimfonct} Soit~$\sch D$ un sous-foncteur définissable et non vide de~$\GG_G^n$. Soit~$H\in \mathsf D_G$ et soit ~$\sch E$ l'ensemble des entiers~$d$ tels que~$\sch D_{|\mathsf D_H}$ contienne un pavé de dimension~$d$ de~$\GG_H^n$. L'ensemble~$\sch E$ est non vide (il contient~$0$) et est inclus dans~$\{0,\ldots, n\}$. Il ne dépend pas de~$H$, et peut alternativement être décrit comme l'ensemble des entiers~$d$ possédant la propriété suivante :  il existe un pavé~$\sch P$ de dimension~$d$ dans~$\GG_H^n$ tel que~$\sch D(H)$ contienne~$\sch P(H)$. Le plus grand élément de~$\sch E$ est appelé la {\em dimension} de~$\sch D$. Par convention, on pose~$\dim{} \emptyset=-\infty$. Il n'y a pas de conflit de terminologie : un pavé de dimension~$d$ au sens de \ref{pavesval} est de dimension~$d$ au sens que nous venons de définir. 

\medskip
Soit~$H\in \mathsf D_G$ et soit~$x\in \sch D(H)$. 

La dimension de~$\sch D$ {\em en~$x$} est la dimension minimale de~$\sch P\cap \sch D_{|\mathsf D_H}\subset \GG_H^n$, où~$\sch P$ est un pavé de dimension~$n$ de~$\GG_H^n$ tel que~$x\in \sch P(H)$. Si~$\sch D$ est de dimension~$d$ en chacun des points de~$\sch D(H)$, il est de dimension~$d$ en chacun des points de~$\sch D(H')$ pour tout~$H'\in \mathsf D_G$ ; on dit alors que~$\sch D$ est {\em purement de dimension~$d$.}

On dit que~$\sch D$ est purement de dimension~$d$ {\em en~$x$} s'il existe un pavé~$\sch P$ de dimension~$n$ dans~$\GG_H^n$ tel que~$x\in \sch P(H)$ et tel que le sous-foncteur~$\sch P\cap \sch D_{|\mathsf D_H}$ de~$\GG_H^n$ soit purement de dimension~$d$.

\trois{fonctdefconn} Un sous-foncteur définissable~$\sch D$ de~$\GG_G^n$ est dit {\em connexe} si pour tout couple~$(\sch D_1,\sch D_2)$ de sous-foncteurs définissables et fermés de~$\GG_G^n$ tels que~$$\sch D=(\sch D_1\cap \sch D)\coprod (\sch D_2\cap \sch D),$$ on  a~$\sch D\subset \sch D_1$ ou~$\sch D\subset \sch D_2$.

\section{Applications de l'élimination des quantificateurs sur les corps valués}\label{elimacvf} 

\deux{lemmouv}{\bf Lemme.} {\em Soit~$K\in \mathsf C$, soit~$\sch X$ un~$K$-schéma intègre, soit~$U$ un ouvert non vide de~$\sch X(K)$ et soit~$\sch V$ un ouvert de Zariski non vide de~$\sch X$. L'intersection~$U\cap \sch V(K)$ est non vide.}

\medskip
{\em Démonstration.} On raisonne par récurrence sur la dimension~$d$ du schéma~$\sch X$. Si~$d=0$ l'assertion est triviale. Supposons que~$d>0$ et que le lemme est vrai en dimension~$<d$. 

\medskip
Quitte à remplacer~$\sch X$ par son normalisé~$\sch X'$ (et~$U$ par son image réciproque sur~$\sch X'$) on peut supposer que~$\sch X$ est normal. Comme~$U$ est non vide, il existe un fermé de Zariski~$\sch Y$ de~$\sch X$, irréductible et de codimension~$1$, tel que~$U$ rencontre~$\sch Y(K)$. Le schéma~$\sch X$ étant normal sur un corps algébriquement clos, il est régulier en codimension~$1$ ; par conséquent, l'intersection~$\sch W$ de~$\sch Y$ et du lieu régulier de~$\sch X$ est non vide. L'hypothèse de récurrence assure alors que~$U\cap \sch W(K)$ est non vide. Il s'ensuit que~$U$ contient un point lisse~$x$ de~$\sch X$.

\medskip
Soit~$(f_1,\ldots, f_d)$ un système de paramètres au voisinage de~$x$. Le théorème des fonctions implicites sur les corps valués algébriquement clos assure l'existence d'un élément~$r\in |K^*|$ et d'un voisinage ouvert définissable~$U'$ de~$x$ dans~$\sch X(K)$, que l'on peut supposer contenu dans~$U$, tel que~$(f_1,\ldots, f_d)$ induise pour toute extension valuée algébriquement close~$L$ de~$K$ une bijection~$U'(L)\simeq B(L)$, où~$B$ désigne le polydisque ouvert de polyrayon~$(r,\ldots, r)$. 

\medskip
Soit~$H$ le groupe abélien ordonné~$|K^*|\oplus \varpi_1^\QQ\oplus\ldots\oplus \varpi_d^\QQ$, où chacun des~$\varpi_j$ est infiniment proche inférieurement de~$1$ relativement au groupe engendré par~$|K^*|$ et les~$\varpi_i$ pour~$i<j$. Soit~$L$ une extension valuée algébriquement close de~$K$ de groupe des valeurs~$H$ (on peut par exemple prendre pour~$L$ une clôture algébrique de~$K(H)$, munie d'une extension quelconque de la valuation de Gauß~$\sum a_h h\mapsto \max |a_h|\cdot h$) ; pour tout~$j$, soit~$\lambda_j$ un élément de~$L$ de valeur absolue~$r\varpi_j$. 

\medskip
Soit~$y$ le point de~$U'(L)$ correspondant au~$d$-uplet~$(\lambda_1,\ldots, \lambda_d)$ de~$B(L)$. Le point~$y$ induit un point~$\bf y$ du schéma~$\sch X$ et une valuation sur~$K({\bf y})$ dont le groupe des valeurs contient les~$\varpi_j$ ; il est donc de rang rationnel au moins~$d$ sur~$|K^*|$, ce qui entraîne que le degré de transcendance de~$K({\bf y})$ sur~$K$ est au moins~$d$. Par conséquent,~$\bf y$ est le point générique de~$\sch X$, et est en particulier situé sur~$\sch V$ ; il s'ensuit que~$y\in \sch V(L)$. 

\medskip
Ainsi,~$U(L)$ rencontre~$\sch V(L)$. Par modèle-complétude,~$U(K)$ rencontre~$\sch V(K)$, ce qu'il fallait démontrer.~$\Box$

\deux{propeq} {\bf Théorème.} {\em Soit~$K$ un corps valué muni d'un plongement de~$|K^*|$ dans~$G$. Soit~$\sch X$ une~$K$-variété algébrique, soit~$(g_{\ell,j})$ une famille finie de fonctions sur~$\sch X$, et soit~$(\gamma_{\ell,j})$ une famille d'éléments de~$G$. Soit~$X$ le foncteur~$$L\mapsto \left\{P\in \sch X(L),\bigvee_{\ell}\bigwedge_j |g_{\ell,j}|(P)\leq \gamma_{\ell,j} \right\}$$ de~$\mathsf C_{K,G}$ dans la catégorie des ensembles. 

\medskip
Soit~${\bf f}=(f_1,\ldots, f_n)$ une famille de fonctions inversibles sur le schéma~$\sch X$, soit~$\pi:\sch X\to \GG_{m,K}^n$ le morphisme induit, et soit~$d$ la dimension de la partie constructible~$\pi(\sch X)$ de~$\GG_{m,K}^n$ (si~$\sch X$ est intègre,~$d$ est le degré de transcendance de~$K(f_1,\ldots, f_n)$ sur~$K$). 

\medskip
\begin{itemize}
\item[A)] Il existe un sous-foncteur définissable et fermé~$\sch D$ de~$\GG_G^n$ tel que~$$|{\bf f}|(X(L))=\sch D(|L^*|)\subset |L^*|^n$$ pour tout~$L\in \mathsf C_{K,G}$, et la dimension de~$\sch D$ est majorée par~$d$. 

\medskip
\item[B)] Supposons~$\sch X$ intègre. Pour tout~$H\in \mathsf D_G$ et tout~${\bf h}=(h_1,\ldots, h_n)\in H^n$, les assertions suivantes sont équivalentes : 

\begin{itemize} 

\medskip
\item[i)]~${\bf h}\in \sch D(H)$ ; 

\item[ii)] il existe une valuation~$|.|$ sur~$K(\sch X)$ prolongeant la valuation de~$K$, à valeurs dans un groupe ordonné~$H'\in \mathsf D_H$, telle que~$$\bigvee_\ell\bigwedge_j |g_{\ell,j}|\leq \gamma_{\ell,j}$$ et telle que~$|f_i|=h_i$ pour tout~$i$. 

\medskip
En particulier,~$\sch D$ ne change pas si l'on remplace~$\sch X$ par un de ses ouverts de Zariski non vides. 

\end{itemize}
\medskip
\item[C)] Supposons que~$X=\sch X_{|\mathsf C_{K,G}}$, et soit~$K^h$ le hensélisé de~$(K,|.|)$. 

\medskip
\begin{itemize}
\item[C1)] Si~$\sch X_{K^h}$ est connexe alors~$\sch D$ est connexe.
\item[C2)] Supposons que~$\pi(\sch X)$ est purement de dimension~$d$. Le sous-foncteur~$\sch D$ est alors purement de dimension~$d$, et égal à~$\GG_G^n$ si~$\sch X\neq \emptyset$ et~$n=d$.
\end{itemize}
\end{itemize}
}

\medskip
{\em Démonstration.} Commençons par prouver A). L'image de~$|{\bf f}|$ est un sous-foncteur définissable du foncteur~$L\mapsto |L^*|$ de~$\mathsf C_{(K,G)}$ dans~$\mathsf{Ens}$, et un tel sous-foncteur est nécessairement de la forme~$L\mapsto \sch D(|L^*|)$ pour un certain sous-foncteur définissable~$\sch D$ de~$\GG_G^n$ (qui est uniquement déterminé). Il reste à s'assurer que~$\sch D$ est fermé et de dimension~$\leq d$. 

\trois{schdferme} {\em Le foncteur~$\sch D$ est fermé.} On peut décrire~$\sch D$ par une combinaison booléenne~$$\bigvee_i \bigwedge_j \phi_{i,j}\Join_{i,j} 1$$ où les~$\phi_{i,j}$ et les~$\Join_{i,j}$ sont comme en \ref{fonctdef}. On peut supposer que pour tout~$i$, le sous-foncteur de~$\sch D$ décrit par la conjonction~$\bigwedge_j \phi_{i,j}\Join_{i,j} 1$ est non vide (sinon, on le retire de la description de~$\sch D$). Nous allons montrer que~$\sch D$ peut être décrit par la condition~$$\bigvee_i \bigwedge_j \phi_{i,j}\leq 1.$$ Il suffit pour cela de montrer que pour tout~$i$ le sous-foncteur de~$\GG_G^n$ défini par la conjonction~$\bigwedge_j \phi_{i,j}\leq 1$ est contenu dans~$\sch D$. 

\medskip
Fixons~$i$, et écrivons~$\phi_j$ et~$\Join_j$ au lieu de~$\phi_{i,j}$ et~$\Join_{i,j}$. Soit~$L\in \mathsf C_{K,G}$ et soit~$(r_1,\ldots, r_n)$ un point de~$|L^*|^n$ en lequel~$\phi_j\leq 1$ pour tout~$j$. Nous allons prouver qu'il appartient à~$\sch D(|L^*|)$, ce qui permettra de conclure. Soit~$\sch D'$ le sous-foncteur de~$\sch D$ défini par la conjonction~$\bigwedge_j \phi_j\Join_j1$ ; il est non vide par hypothèse. 

La formule~$(h_1,\ldots, h_n)\mapsto \max (h_i\inv r_i, h_i r_i\inv)_i$ définit un morphisme de foncteurs de~$\sch D'_{|\mathsf D_{|L^*|}}$ dans~$\GG_{|L^*|}$.  Son image est un sous-foncteur définissable de~$\GG_{|L^*|}$ ; par o-minimalité, celui-ci admet une borne inférieure~$a$, appartenant à~$|L^*|$ et supérieure ou égale à~$1$. 

\medskip
{\em Montrons que~$a=1$. } On raisonne par l'absurde. Supposons~$a>1$.  Par définition
de~$a$, il existe~$(s_1,\ldots,s_n)\in |L^*|^n$ tel que les propriétés suivantes soient vérifiées : 

\medskip
$\alpha)$~$a^{-3/2}\leq s_i\leq a^{3/2}$  pour tout~$i$ ; 

$\beta)$ il existe~$i$ tel que~$s_i\leq a^{-1}$ ou~$s_i\geq a$ ; 

$\gamma)$~$(r_1s_1,\ldots, r_ns_n)\in \sch D'(|L^*|)$. 

\medskip
Par convexité, le point~$(r_1s_1^{1/2},\ldots, r_ns_n^{1/2})$ satisfait l'inégalité~$\phi_j\leq 1$ pour tout~$j$. Par minimalité de~$a$ et comme~$a^{3/4}<a$ puisque~$a>1$, il ne peut appartenir à~$\sch D'(|L^*|)$. Il existe donc~$j$ tel que~$\Join_j$ soit égal à~$<$ et tel que~$\phi_j$ vaille~$1$ sur~$(r_1s_1^{1/2},\ldots, r_ns_n^{1/2})$.  L'application~$\phi_j$ prend alors une valeur inférieure ou égale à~$1$ en~$(r_1,\ldots, r_n)$, la valeur 1 en~$(r_1s_1^{1/2},\ldots, r_ns_n^{1/2})$, et une valeur strictement inférieure à~$1$ sur~$(r_1s_1,\ldots, r_ns_n)$ ; mais c'est contradictoire avec le caractère affine de~$\phi$. 

\medskip
Soit~$H$ le groupe ordonné divisible~$|L^*|\oplus \varpi^\QQ$, où~$\varpi$ est infiniment proche supérieurement de~$1$. Soit~$F$ une extension valuée algébriquement close de~$L$ de groupe des valeurs~$H$. 

\medskip
Comme~$a=1$, il existe~$(s_1,\ldots, s_n)\in \varpi ^\QQ$ tels que le~$n$-uplet~$(r_1s_1,\ldots, r_ns_n)$ appartienne à~$\sch D'(H)=\sch D'(|F^*|)$. Comme~$\sch D'$ est un sous-foncteur de~$\sch D$, il existe~$P\in X(F)$ tel que~$|{\bf f}|(P)=(r_1s_1,\ldots, r_ns_n)$. 

\medskip
Le groupe~$\varpi ^{\QQ}$ est un sous-groupe convexe de~$H$ ; le quotient~$ H/\varpi^\QQ$ s'identifie à~$|L^*|$. La valuation de~$F$ peut donc être rendue un peu plus grossière en la composant avec la flèche quotient~$H\to |L^*|$ ; le corps valué~$F_0$ ainsi obtenu a même corps sous-jacent que~$F$, et est toujours une extension valuée de~$L$. Comme le foncteur~$X$ est défini par des inégalités larges,~$P\in X(F_0)$ ; son image par~$|{\bf f}|$ dans~$|F_0^*|^n$ appartient donc à~$\sch D(F_0^*)$. Mais~$|F_0^*|=|L^*|$ ; et l'image de~$P$ dans~$|F_0^*|^n$ est égale à l'image de~$(r_1s_1,\ldots, r_ns_n)$ par la flèche quotient~$H^n\to |L^*|^n$, c'est-à-dire à~$(r_1,\ldots, r_n)$. Ainsi, ce dernier appartient à~$\sch D(|L^*|)$, ce qu'il fallait démontrer.

\trois{dimschdmaj} {\em Majoration de la dimension de~$\sch D$.} Soit~$L\in \mathsf C_{K,G}$, et soit~$\sch P$ un pavé de~$\GG_{|L^*|}^n$ contenu dans~$\sch D_{|\mathsf D_{|L^*|}}$. On note~$\delta$ la dimension de~$\sch P$ ; nous allons montrer que~$\delta\leq d$. Quitte à composer~$\pi$ avec une isogénie du tore~$ {\mathbb G}_{m,K}^n$, on peut supposer que~$\sch P$ est décrit  (universellement) par des conditions de la forme~$$r_i<h_i<R_i, 1\leq i\leq \delta\;{\rm et}\;h_i=\gamma_i, i\geq \delta+1,$$ où les~$r_i, R_i$ et~$g_i$ appartiennent à~$|L^*|$. 

\medskip
Soit~$H$ le groupe ordonné divisible~$|L^*|\oplus \varpi_1^\QQ\oplus\ldots\oplus \varpi_{\delta}^\QQ$, où chacun des~$\varpi_j$ est infiniment proche supérieurement de~$1$ relativement au groupe engendré par~$|L^*|$ et les~$\varpi_i$ pour~$i<j$. Soit~$F$ une extension valuée algébriquement close de~$L$ de groupe des valeurs~$H$. 

\medskip
L'élément~$(r_1\varpi_1,\ldots, r_\delta\varpi_{\delta},g_{\delta+1},\ldots, g_n)$ de~$ H^n$ appartient à~$\sch P(H)$, et donc à~$\sch D(H)$. Il est dès lors égal à~$|{\bf f}|(P)$ pour un certain~$P\in X(F)$. Le point~$P$ induit un point~$x$ du schéma~$\sch X_L$, et une valuation~$\val.$ sur~$L(x)$. 

\medskip
Si l'on pose~$L_0=L(f_1(x),\ldots, f_n(x))\subset L(x)$ alors~$|L_0^*|$ contient les~$r_i\varpi_i$, qui sont~$\QQ$-linéairement indépendants modulo~$|L^*|$ ; par conséquent, le degré de transcendance de~$L_0$ sur~$L$ est au moins égal à~$\delta$. 

\medskip
Le corps~$L_0$ s'identifie par ailleurs au corps~$L(\pi(x))$ ; il s'ensuit que le degré de transcendance de~$L_0$ sur~$L$ est majoré par~$d$, puis que~$\delta\leq d$, ce qu'il fallait démontrer.

\trois{dinvarbir} {\em Preuve de l'assertion B)}. Soit~$H\in \mathsf D_G$ et soit~${\bf h}=(h_1,\ldots, h_n)\in H^n$. 

\medskip
{\em Supposons que i) est vraie.} Soit~$L\in \mathsf C$ un corps tel que~$|L^*|=H$. Soit~$X_{\bf h}$ le foncteur qui associe à un corps~$F\in \mathsf C_L$ l'ensemble des points de~$X(F)$ en lesquels~$|{\bf f}|$ est égal à~$\bf h$. Comme~$\bf h$ appartient à~$\sch D(H)$, le foncteur~$X_{\bf h}$ est non vide ; pour tout~$F\in \mathsf C_L$, l'ensemble~$X_{\bf h}(F)$ est un ouvert de~$\sch X(F)$. 

\medskip
Soit~$\sch Y$ un fermé de Zariski strict de~$\sch X$. En vertu du lemme \ref{lemmouv}, le sous-foncteur~$X_{\bf h}$ de~$\sch X_L$ n'est pas contenu dans le sous-foncteur~$\sch Y_L$ de~$\sch X_L$. Par le théorème de compacité en théorie des modèles,  il existe~$F\in \mathsf C_L$ et un point~$x\in X_{\bf h}(F)$ qui n'appartient à~$\sch Y(F)$ pour aucun fermé de Zariski strict~$\sch Y$ de~$X$. 

\medskip
Cela signifie que~$x$ induit le point générique de~$\sch X$. Par conséquent,~$K(\sch X)$ se plonge dans~$F$, et si l'on note encore~$|.|$ la restriction à~$K(\sch X)$ de la valeur absolue de~$F$ on a~$$\bigvee_\ell\bigwedge_j |g_{\ell,j}|\leq \gamma_{\ell,j}$$ et~$|f_i|=|h_i|$ pour tout~$i$ (puisque~$x\in \sch X_{\bf h}(F)$). Ainsi, ii) est vraie en prenant~$H'=|F^*|$. 

\medskip
{\em Supposons que ii) est vraie.} Soit~$F\in \mathsf C_{K,G}$ une extension valuée de~$(K(\sch X),\val.)$ telle que~$|F^*|=H'$ 
et soit~$P$ le point de~$\sch X(F)$ défini par la flèche composée~$\spec F\to \spec K(\sch X)\to \sch X$. On a par construction~$$\bigvee_\ell\bigwedge_j |g_{\ell,j}|(P)\leq \gamma_{\ell,j}.$$ 
Par conséquent,~$P\in X(F)$. Dès  lors,~${\bf h}=|{\bf f}|(P)\in \sch D(|F^*|)\cap H^n=\sch D(H)$. Ainsi, i) est vraie. 

\trois{connpurdim} {\em Preuve de C).} Comme~$X=\sch X_{|\mathsf C_{K,G}}$, on peut définir~$X$ sans faire appel à des éléments de~$G\setminus |K^*|$, ce qui permet de supposer que~$G=|K^*|^{\QQ}$, et que~$\mathsf C_{K,G}=\mathsf C_K$. Remarquons qu'on aurait pu tout aussi bien prendre~$G$ égal à~$|K^*|$ ; mais en optant pour~$|K^*|^{\QQ}$, on ne modifie ni la catégorie~$\mathsf D_G$, ni la notion de foncteur définissable sur~$\mathsf D_G$, et l'on se donne le droit de remplacer~$K$ par n'importe laquelle de ses extensions algébriques, dont le groupe se plonge naturellement dans~$\sqrt{|K^*|}$. 

\medskip
Soit~$\sch E$ un sous-foncteur définissable et fermé de~$\GG_G^n$. Comme~$G=|K^*|^{\QQ}$, la condition~$|{\bf f}|\in \sch E(|K(\sch X)^*|)$ définit un ouvert~$\PP_{K(\sch X)/K}\{\sch E\}$ de l'espace de Zariski-Riemann~$\PP_{K(\sch X)/K}$ (les corps~$K$ et~$K(\sch X)$ sont considérés comme trivialement gradués ; il s'agit donc d'un espace de Zariski-Riemann classique). 

\medskip
Pour démontrer C1) on peut, quitte à étendre les scalaires, se ramener au cas où~$K=K^h$. Supposons que le foncteur~$\sch D$ s'écrive comme la réunion disjointe de deux sous-foncteurs définissables fermés~$\sch D_1$ et~$\sch D_2$ de~$\GG_G^n$. En vertu de l'assertion B) déjà établie, et plus précisément de ii)$\Rightarrow$i), on a alors~$$\PP_{K(\sch X)/K}=\PP_{K(\sch X)/K}\{\sch D_1\}\coprod \PP_{K(\sch X)/K}\{\sch D_2\}.$$

\medskip
L'espace topologique~$\PP_{K(\sch X)/K}$ s'identifie à la limite projective des fibres spéciales des modèles propres et intègres de~$K(\sch X)$ sur l'anneau des entiers~$K\zero$ de~$K$. Or si~$\sch Y$ est un tel modèle, sa fibre spéciale est connexe. En effet, soit~$$\sch Y\to \spec B\to \spec K\zero$$ la factorisation de Stein de~$\sch Y\to \spec K\zero$. La~$K\zero$-algèbre finie~$B$ est intègre ; c'est d'autre part un produit d'anneaux locaux puisque~$K\zero$ est hensélien, et elle est de ce fait locale ; la fibre spéciale de~$\spec B\to \spec K\zero$ est donc un singleton, ce qui équivaut à la connexité de la fibre spéciale de~$\sch Y$. 

\medskip
Il en découle que~$\PP_{K(\sch X)/K}$ est connexe. En conséquence, il existe~$i$ tel que~$\PP_{K(\sch X)/K}\{\sch D_i\}=\emptyset$. En utilisant à nouveau l'assertion B), mais cette fois-ci à travers l'implication i)$\Rightarrow$ii),  il vient~$\sch D_i=\emptyset$ ; il s'ensuit que~$\sch D$ est connexe.

\medskip
Démontrons maintenant C2). On peut, quitte à étendre les scalaires, supposer~$K$ algébriquement clos. Commençons par une remarque, que nous utiliserons à plusieurs reprises : si~$T_1,\ldots, T_n$ désignent les fonctions coordonnées sur~$\GG_{m,K}^n$ alors~$|{\bf  f}|(\sch X(L))=|{\bf T}|(\pi(\sch X)(L))$ pour tout~$L\in \mathsf C_K$. 

\medskip
On procède par récurrence sur~$d$. Si~$d=0$ alors~$\pi(\sch X)$ est une union finie de~$K$-points, d'où l'assertion. Supposons que~$d>0$ et que l'assertion a été démontrée en dimension~$<d$. On peut écrire~$\pi(\sch X)$ comme une union finie de parties localement fermées irréductibles de dimension~$d$ ; cela permet de se ramener au cas où~$\sch X$ est lui-même un sous-schéma irréductible de~$\GG_{m,K}^n$ de dimension~$d$ et où~$f_i=T_i$ pour tout~$i$. On sait d'après C1) que~$\sch D$ est alors connexe. 

\medskip
Fixons~$L\in \mathsf C_K$.  Soit~$\xi\in \sch D(|L^*|)$ et soit~$\delta$ la dimension de~$\sch D$ en~$\xi$. L'assertion A) assure que~$\delta\leq d$. 

\medskip
{\em La dimension~$\delta$ ne peut être nulle.} En effet si c'était le cas~$\sch D$ serait par connexité un singleton. Mais comme~$\sch X$ est de dimension strictement positive, l'une au moins des projections de~$\sch X$ sur~$\GG_{m,K}$ est dominante et la fonction~$T_i$ correspondante n'est dès lors pas bornée sur~$\sch X(L)$, d'où une contradiction. 

\medskip
Comme~$\delta>0$, il existe une application~$\phi : |L^*|^n\to |L^*|$ de la forme~$(h_1,\ldots, h_n)\mapsto a\prod h_i^{e_i}$ où~$a\in |L^*|$, où les~$e_i$ appartiennent à~$\ZZ$, et qui possède la propriété suivante :~$\phi(\xi)=1$, et si~$\sch H$ désigne le sous-foncteur de~$\GG^n_{|L^*|}$ défini par la condition~$\phi=1$, alors~$\sch H\cap  \sch D$ est de dimension~$\delta-1$ en~$\xi$. 

\medskip
Soit~$x\in \sch X(L)$ tel que~$|{\bf T}|(x)=\xi$, et soit~$\alpha\in L^*$ tel que~$|\alpha|=a$. Posons~$\beta=\alpha\prod T_i(x)^{e_i}$ ; comme~$\xi\in \sch H(|L^*|)$, on a~$|\beta|=1$. Soit~$\sch Y$ l'intersection de~$\sch X$ et du sous-schéma fermé de~$\GG_{m,K}^n$ d'équation~$\alpha\prod T_i^{e_i}-\beta=0$. 

\medskip
Comme~$\sch H\cap \sch D$ est de dimension~$\delta-1$ en~$x$, on ne peut avoir~$\sch D\subset \sch H$ ; par conséquent,~$\sch X$ n'est pas contenu dans~$\sch Y$. Le {\em Hauptidealsatz} assure alors que~$\sch X\cap \sch Y$ est purement de dimension~$d-1$. 

\medskip
Soit~$\sch D'$ le sous-foncteur définissable de~$\GG^n_{|L^*|}$ qui décrit l'image de~$\sch Y$ par~$|{\bf T}|$. Comme~$x\in \sch Y$, on a~$\xi\in \sch D'(|L^*|)$ ; la définition même de~$\sch Y$ assure par ailleurs que~$\sch D'\subset \sch H\cap \sch D$. Il s'ensuit que la dimension de~$\sch D'$ en~$\xi$ est majorée par~$\delta-1$. 

\medskip
Mais en vertu de notre hypothèse de récurrence, cette dimension est égale à~$d-1$. On a donc~$d-1\leq \delta-1\leq d-1$, et partant~$\delta=d$, ce qu'il fallait démontrer. 

\medskip
Supposons pour terminer que~$\sch X\neq \emptyset$ et que~$d=n$, et montrons que~$\sch D=\GG_G^n$. 
Sous nos hypothèses,~$\pi(\sch X)$ contient un ouvert dense~$\sch U$ de~$\GG^n_{m,K}$ ; quitte à remplacer~$\sch X$ par~$\sch U$, on peut supposer que~$\sch X$ est un ouvert dense de~$\GG^n_{m,K}$, et que~${\bf f}={\bf T}$. Mais l'assertion B) permet de remplacer~$\sch X$ par~$\GG^n_{m,K}$ tout entier, et l'égalité~$\sch D=\GG_G^n$ est dès lors évidente.~$\Box$ 

%
%

\subsection*{Application aux espaces de Zariski-Riemann gradués}

\deux{theoconetem} {\bf Théorème.} {\em Soit~$K$ un corps {\em gradué}, soit~$L$ une extension graduée de~$K$, et soit~$\sch U$ un ouvert quasi-compact de~$\PP_{L/K}$. Soit~$(f_1,\ldots, f_n)$ une famille d'éléments homogènes non nuls de~$K$ et soit~$d$ le degré de transcendance gradué de~$K(f_1,\ldots, f_n)$ sur~$K$. 

\medskip
Soit~$\sch D_{K,L,\sch U,(f_i)}$ le sous-foncteur de~$\GG^n$ qui associe à tout~$H\in \mathsf D$ l'ensemble des~$n$-uplets~$(h_1,\ldots, h_n)\in  H^n$ tels qu'il existe une valuation graduée~$\val .\in \sch U$ à valeurs dans 
un groupe~$H'\in \mathsf D_H$, et satisfaisant les conditions~$\val {f_i}=h_i$ pour tout~$i$.

\medskip
1) Le sous-foncteur~$\sch D_{K,L,\sch U,(f_i)}$ de~$\GG^n$ est définissable, fermé, et de dimension majorée par~$d$. 

\medskip
2) Si~$\sch U= \PP_{L/K}$ alors~$\sch D_{K,L,\sch U,(f_i)}$ est purement de dimension~$d$, et si de plus~$n=d$ alors~$\sch D_{K,L,\sch U,(f_i)}=\GG^n$. 
}

\medskip
{\em Démonstration.} L'ouvert~$\sch U$ est de la forme~$\bigcup\limits_\ell\PP_{L/K}\{g_{\ell,j}\}_j$ pour une certaine famille finie~$(g_{\ell,j})$ d'éléments homogènes de~$L$. 

\trois{redstrictgrad} {\em Réduction au cas «strict».} Soit~$\bf r$ un polyrayon~$K$-libre tel que les degrés des~$f_i$ et des~$g_{\ell,j}$ soient de torsion modulo le groupe des degrés de~$K_{\bf r}$ ; soit~$\sch U_{\bf r}$ l'ouvert~$\bigcup\limits_\ell\PP_{L({\bf r}\inv {\bf T})/K_{\bf r}}\{g_{\ell,j}\}_j$ ; l'ouvert~$\sch U_{\bf r}$ est l'image réciproque de~$\sch U$ sur~$\PP_{L({\bf r}\inv {\bf T})/K_{\bf r}}$. 

\medskip
Soit~$H\in \mathsf D$ et soit~$\val .$ une~$K$-valuation graduée sur~$L$, à valeurs dans~$H$. Il existe une~$K_{\bf r}$-valuation graduée sur~$L_{\bf r}$ qui est à valeurs dans~$H$ et étend~$\val .$ : il suffit par exemple de considérer la valuation de Gauß~$\eta_{(L,\val .), {\bf r}, {\bf 1}}$. 

\medskip
Il s'ensuit que~$\sch D_{K,L,\sch U,(f_i)}$ est égal à~$\sch D_{K_{\bf r},L_{\bf r},\sch U_{\bf r},(f_i)}$. Par ailleurs toute famille d'éléments homogènes de~$L$ algébriquement indépendante sur~$K$ reste algébriquement indépendante sur le corps gradué~$K_{\bf r}$ ; le degré de transcendance de~$K_{\bf r}(f_1,\ldots, f_n)$ sur~$K_{\bf r}$ est donc égal à~$d$. Enfin, comme~$\sch U_{\bf r}$ est l'image réciproque de~$\sch U$ sur~$\PP_{L({\bf r}\inv {\bf T})/K_{\bf r}}$, on a~$\sch U_{\bf r}=\PP_{L({\bf r}\inv {\bf T})/K_{\bf r}}$ dès que~$\sch U=\PP_{L/K}$. 

\medskip  On peut donc, quitte à remplacer~$K$ par~$K_{\bf r}$,~$L$ par~$L({\bf r}\inv {\bf T})$, et~$\sch U$ par~$\sch U_{\bf r}$, supposer que les degrés des~$f_i$ et des~$g_{\ell,j}$ sont de torsion modulo le groupe des degrés de~$K$.

\trois{rednongrad} {\em Réduction au cas non gradué.} Comme~$$\PP_{L/K}\{g_1,\ldots, g_m\}=\PP_{L/K}\{a_1g_1^N, \ldots, a_ng_m^N\}$$ pour tout entier~$N>0$, toute famille~$(a_i)$ d'éléments homogènes de~$K$ et toute famille~$(g_i)$ d'éléments homogènes de~$L$,  on peut supposer que les~$g_{\ell,j}$ sont tous de degré 1. 

\medskip
D'autre part, si~$N>0$ et si~$(a_i)$ est une famille d'éléments homogènes non nuls de~$K$, le degré de transcendance de~$K(a_1f_1^N,\ldots, a_nf_n^N)$ sur~$K$ est égal à~$d$, et~$\sch D_{K,L,\sch U, (a_if_i^N)}(H)$ est pour tout~$H\in \mathsf D$ l'image de~$\sch D_{K,L,\sch U, (f_i)}(H)$ par l'homothétie~${\bf h}\mapsto {\bf h}^N$. Par conséquent, on peut, quitte à remplacer~$(f_i)$ par une famille~$(a_if_i^N)$ convenable, supposer que les~$f_i$ sont tous de degré 1. 

\medskip
Posons~$\sch U^1=\bigcup\limits_\ell \PP_{L^1/K^1}\{g_{\ell,j}\}_j$. L'ouvert~$\sch U$ est l'image réciproque de~$\sch U^1$ sur~$\PP_{L/K}$. 

\medskip
Soit~$H\in \mathsf D$ et soit~$\val .$ une~$K^1$-valuation sur~$L^1$, à valeurs dans~$H$. Il existe une~$K$-valuation graduée sur~$L$ qui est à valeurs dans~$H$ et dont la restriction à~$L^1$ est égale à~$\val.$ : par exemple, on peut prolonger~$\val .$ à~$L^1\cdot K$, en envoyant~$\lambda \mu$ sur~$\val \lambda$ pour tout~$\lambda\in L^1$ et tout élément homogène non nul~$\mu\in K$, puis l'étendre de façon quelconque à~$L$. 

\medskip
Il s'ensuit que~$\sch D_{K^1,L^1,\sch U_1,(f_i)}$ est égal à~$\sch D_{K,L,\sch U_,(f_i)}$. Par ailleurs toute famille d'éléments homogènes de~$L^1$ algébriquement indépendante sur~$K^1$ reste algébriquement indépendante sur le corps gradué~$K$ ; le degré de transcendance de~$K^1(f_1,\ldots, f_n)$ sur~$K^1$ est donc égal à~$d$. Enfin, comme~$\PP_{L/K}\to \PP_{L^1/K^1}$ est surjective par ce qui précède, on a~$\sch U_1=\PP_{L^1/K^1}$ dès que~$\sch U=\PP_{L/K}$.

\medskip
On peut ainsi remplacer~$K$ par~$K^1$,~$L$ par~$L^1$, et~$\sch U$ par~$\sch U^1$, et se ramener par ce biais à un problème de théorie des valuations traditionnelle (non graduée). 

\trois{conetemknongrad} {\em Réduction au cas où~$L$ est de type fini sur~$K$.} Soit~$L'$ le sous-corps de~$L$ engendré par~$K$, les~$f_i$ et les~$g_{\ell,j}$, et soit~$\sch U'$ l'ouvert quasi-compact~$\bigcup\limits_\ell\PP_{L'/K}\{g_{\ell,j}\}_j$. L'ouvert~$\sch U$ est égal à l'image réciproque de~$\sch U'$ sur~$\PP_{L/K}$.

\medskip
Soit~$H\in \mathsf D$. 
Toute~$K$-valuation sur~$L'$ à valeurs dans~$H$ s'étend en une~$K$-valuation sur~$L$ à valeurs dans~$H$. On en déduit que~$\sch D_{K,L',\sch U',(f_i)}$ est égal à~$\sch D_{K,L,\sch U_,(f_i)}$, et que~$\sch U'=\PP_{L'/K}$ dès que~$\sch U=\PP_{L/K}$. On peut donc remplacer~$L$ par~$L'$ et~$\sch U$ par~$\sch U'$, c'est-à-dire supposer~$L$ de type fini sur~$K$.

\trois{ltfsurk} {\em Preuve dans le cas non gradué lorsque~$L$ est de type fini sur~$K$.} Il existe alors un~$K$-schéma intègre et de type fini~$\sch X$ dont le corps des fonctions est isomorphe à~$L$, sur lequel les~$g_{\ell,j}$ et les~$f_i$ sont définies, et sur lequel les~$f_i$ sont inversibles. Considérons le corps~$K$ comme trivialement valué. 

Soit~$X$ le foncteur qui associe à~$F\in \mathsf C_K$ le sous-ensemble de~$\sch X(F)$ formé des points~$P$ tels que~$$\bigvee_\ell\bigwedge_j |g_{\ell,j}|(P)\leq 1.$$

\medskip
En vertu de l'assertion A) du théorème \ref{propeq}, il existe un sous-foncteur définissable fermé~$\sch D$ de~$\GG^n$, de dimension~$\leq d$, tel que pour toute~$F\in \mathsf C_K$, l'on ait~$|{\bf f}|(X(F))=\sch D(|F^*|)$.

\medskip
L'assertion B) de {\em loc. cit.} assure que~$\sch D_{K,L,\sch U, (f_i)}$ coïncide avec~$\sch D$, ce qui prouve 1). 

\medskip
{\em Remarque.} Le schéma~$\sch X$ et le foncteur~$X$ dépendent, entre autres choix, de l'écriture~$\sch U=\bigcup\limits_\ell\PP_{L/K}\{g_{\ell,j}\}_j$ fixée au départ. 

\medskip
Il reste à démontrer 2). Supposons que~$\sch U=\PP_{L/K}$. On peut alors opter pour l'écriture~$\sch U=\PP_{L/K}\{\emptyset\}$, auquel cas~$X=\sch X$. Les propriétés requises découlent dès lors de l'assertion C2) de {\em loc. cit.}~$\Box$ 

\section{Prolongements des valuations de Gauß}\label{Gau}

\subsection*{Séparation des prolongements d'une valuation}

\deux{valanneau} Suivant Huber, on définit une valuation sur un {\em anneau} (commutatif, unitaire)~$A$ de deux façons équivalentes : 

i) c'est la donnée d'un point~$x$ de~$\spec A$ et d'une valuation~$\val .$ sur le corps résiduel~$\kappa(x)$ ; 

ii) c'est une application~$\val . :A\to G\cup\{0\}$, où~$G$ est un groupe abélien ordonné, telle que~$\val 0=0, \val 1=1, \val {ab}=\val a \cdot \val b$ et~$\val{a+b}\leq \max \val a, \val b$
pour tout couple~$(a,b)$.

On passe de i) à ii) en considérant la composée de~$A\to \kappa(x)$ et de~$\val.$, et de~ii) à~i)
en prenant pour~$x$ le point de~$\spec A$ correspondant à
 l'idéal premier~$\{a\in A,\val a=0\}$, et en munissant~$\kappa(x)$
de l'application induite par~$\val .~$ par passage au quotient. 

On dispose d'une notion d'équivalence entre deux valuations sur~$A$ : si on utilise la définition i), on demande qu'elles
soient données par le même point de~$\spec A$ et par deux valuations équivalentes sur son corps résiduel ; si on utilise ii), 
on impose une condition analogue à celle de~\ref{valgrad}. Le plus souvent, les valuations seront implicitement considérées
à équivalence près. 

\deux{valanneauextgal} Soit~$A$ un anneau, et soit~$B$
une~$A$-algèbre galoisienne de groupe~$\Sigma$.  
L'ensemble des valuations de~$A$ s'identifie alors au quotient
de l'ensemble des valuations de~$B$ par l'action de~$\Sigma$. 

En effet, soit~$\val .$
une valuation sur~$A$.
Elle est donnée par un couple~$(x,\val ._0)$, où~$x$ est un point
de~$\spec A$ et où~$\val ._0$ est une valuation sur~$\kappa(x)$.
En choisissant un antécédent~$y$ de~$x$
sur~$\spec B$ et un prolongement de la
valuation~$\val ._0$
à~$\kappa(y)$, on obtient un prolongement de la
valuation~$\val .$ de~$A$ à~$B$. 

Donnons-nous maintenant deux 
prolongements~$\val .'$ et~$\val .''$ de
la valuation~$\val .$ 
de~$A$ à 
à~$B$, le premier correspondant
à un couple~$(y,\val .'_0)$ et le second 
à un couple~$(z,\val .''_0)$. Nous allons montrer
qu'ils sont dans la même orbite sous l'action de~$\Sigma$. 
Comme~$y$ et~$z$ sont
deux antécédents de~$x$, il existe un élément
de~$\Sigma$,
qui envoie~$y$ sur~$z$,
ce qui permet de supposer que~$y=z$. Les deux valuations
~$\val'_0$ et~$\val ''_0$ prolongeant
$\val ._0$, elles appartiennent à la même orbite
sous~${\rm Gal}\;(\kappa(y)/\kappa(x))$,
dont tout élément se relève en un 
élément de~$\Sigma$ ; ceci achève
de prouver l'assertion requise.

\deux{valanneauextsuite} Soit~$K$ un corps, soit~$A$ une~$K$-algèbre et soit~$L$
une extension galoisienne de~$K$. Soit~$\val .$ une valuation 
sur~$K$ ; fixons un prolongement~$\val .'$ de~$\val .$ à~$L$. L'ensemble des prolongements
de~$\val .$ à~$A$ s'identifie alors au quotient de
l'ensemble des prolongements de~$\val .'$ à $L\otimes_KA$ par l'action de ${\rm Gal}\;(L/K)$. 

En effet, en vertu de~\ref{valanneauextgal}, il suffit de vérifier que tout prolongement~$\val .''$
de~$\val .$ à~$A$ s'étend en un prolongement de~$\val .'$ à~$L\otimes_KA$. Soit donc~$\val .''$ un prolongement
de~$\val .$ à~$A$. D'après~\ref{valanneauextgal}, il admet un prolongement~$\val .'''$ à~$L\otimes_KA$. Il existe
 $g\in {\rm Gal}\;(L/K)$ tel que~$\val .'''_{| L}\circ g=\val .'$ ; la valuation~$\val .'''\circ g$ sur~$L\otimes_KA$ répond alors
 aux conditions souhaitées.

\deux{intropresc} Soit~$K$ un corps et soit~$A$ une~$K$-algèbre {\em finie}. Soit~$E$ un sous-ensemble fini de~$A$, soit~$\Sigma$ un sous-groupe de~$\got S_E$, et soit~$\val .$ une valuation sur~$K$. Tout prolongement de~$\val .$ à~$A$ induit une fonction de~$E$ dans~$|K^*|^\QQ\cup\{0\}$. Le nombre de fonctions de~$E$ dans~$|K^*|^\QQ\cup\{0\}$ 
distinctes modulo l'action de~$\Sigma$ que l'on peut obtenir par ce biais sera  noté~$\got n(\val ., A, E,\Sigma)$. 

\deux{lemmsepval} {\bf Lemme.} {\em Soit~$K$ un corps, soit~$A$ une~$K$-algèbre finie, soit~$E$ un sous-ensemble fini de~$A$ et soit~$\Sigma$ un sous-groupe de~$\got S_E$ ; soit~$N$ un entier. 

1) L'ensemble des valuations~$\val .$ de~$A$ telles que~$\got n(\val., A,E,\Sigma)=N$ est une partie constructible~$\sch U_{K,A,E,\Sigma,N}$ de~$\PP_K$. 

2) Pour toute extension~$F$ de
$K$, le constructible~$\sch U_{F,F\otimes_KA,E,\Sigma,N}$ est égal à l'image réciproque de~$\sch U_{K,A,E,\Sigma,N}$ sur~$\PP_F$. }

\medskip
{\em Démonstration.}  Quitte à remplacer~$A$ par sa sous-algèbre engendrée
par~$E$, on peut supposer que~$A=K[E]$. Choisissons une famille génératrice ~$(P_1,\ldots, P_r)$ du noyau du morphisme 
surjectif~$K[X_e]_{e\in E}\to A$ d'évaluation en les éléments de~$E$.

\medskip
Soit~$F$ une extension de~$K$ et soit~$\val .$ une valuation sur~$F$. Choisissons une extension
non trivialement valuée et algébriquement close~$\FF$ de~$F$. 
Toute valuation de~$F\otimes_KA$ prolongeant~$\val .~$ est donnée par un~$F$-morphisme de~$F\otimes_KA$ dans~$\FF$, 
c'est-à-dire par une famille ~$(x_e)_{e\in E}$ d'éléments de~$\FF$ annulant les~$P_j$. Deux prolongements de~$\val .$,  correspondant à deux familles~$(x_e)$ et~$(y_e)$,
définissent 
deux fonctions distinctes sur~$E$ si et seulement si il existe~$e$ tel que~$|x_e|\neq |y_e|$. 

\medskip
Soit~$(b_i)_{i\in I}$ la famille des coefficients des polynômes~$P_j$. En vertu de ce qui précède, l'élimination des quantificateurs dans la théorie des corps non trivialement valués algébriquement clos assure qu'il existe une formule du premier ordre~$\Phi_N$ dans le langage des corps valués, sans paramètres ni quantificateurs, possédant les propriétés suivantes : 

$\bullet$ l'ensemble des variables libres de~$\Phi_N$ est indexé par~$I$, et chacune d'elles
 vit dans le corps valué ; 

$\bullet$ pour toute extension~$F$ de~$K$ et toute valuation
$\val .~$ sur~$F$, l'énoncé~$\Phi_N(b_i)_{i\in I}$ est vrai dans~$(F,\val.)$ si et seulement si l'on a~$\got n(\val .,F\otimes_K A, E,\Sigma)=N$.

\medskip
Le lemme est alors établi, en prenant pour~$\sch U_{K,A,E,\Sigma,N}$ la partie constructible
de~$\PP_K$ définie par la condition~$\Phi_N(b_i)_i$.~$\Box$

\deux{defesepval} Soit~$K$ un corps, soit~$A$ une extension~$K$-algèbre finie et soit~$E$ un sous-ensemble fini de~$A^*$. Si~$\val .$ est une valuation sur~$K$, on dira que~$E$ {\em sépare les prolongements de~$\val .~$ à~$A$} si~$\got n(\val ., A, E,\{{\rm Id}_E\})$ est égal au nombre de valuations de~$A$ prolongeant~$\val .$. 

\subsection*{Prolongement des valuations de Gauß}

\deux{equivprogauss} {\bf Proposition.} {\em Soit~$K$ un corps valué et soit~$A$ une~$K(T_1,\ldots, T_n)$-algèbre finie. Soit~$\sch C$ une classe d'extensions valuées de~$K$ telle que~$K\in \sch C$. Pour tout entier~$N$, tout~$\Lambda\in \sch C$ et tout groupe~$G\in \mathsf D_{|\Lambda^*|}$, notons~$\sch D_{\Lambda, N}(G)$ le sous-ensemble de~$G^n$ formé des~$n$-uplets~$\bf g$ tels que~$\eta_{\Lambda, \bf g}$ admette exactement~$N$ prolongements à~$\Lambda({\bf T})\otimes_{K({\bf T})}A$. Les assertions suivantes sont équivalentes :

\medskip
i) pour tout~$N$, le sous-foncteur~$\sch D_{K,N}$ de~$\GG_{|K^*|}^n$ est définissable, et ~$\sch D_{\Lambda,N}$ est la restriction de ~$\sch D_{K,N}$ à~$\GG_{|\Lambda^*|}^n$ pour tout~$\Lambda\in \sch C$ ;

ii) il existe un sous-ensemble fini~$E$ de~$A$ qui sépare les prolongements 
de~$\eta_{\Lambda, \bf g}$ à~$\Lambda({\bf T})\otimes_{K({\bf T})}A$ pour tout~$\Lambda\in \sch C$, tout~$G\in \mathsf D_{|\Lambda^*|}$ et tout~${\bf g}\in G^n$.} 

\medskip
{\em Démonstration.} Pour alléger un peu les notations, nous écrirons~$\sch D_N$ au lieu de~$\sch D_{K,N}$. On procède en deux temps. 

\trois{schck} {\em Le cas où~$\sch C=\{K\}$}. 

\medskip
Supposons que i) est vraie. Soit~$N$ un entier. Soit~$G\in \mathsf D_{|K^*|}$ et soit~${\bf g}$ appartenant à~$\sch D_N(G)$. La valuation~$\eta_{\bf g}$ admet~$N$ prolongements à~$A$. Il existe un ensemble fini~$E_{\bf g}\subset A$ qui sépare les prolongements de~$\eta_{\bf g}$ à~$A$, c'est-à-dire tel que~$\got n(\eta_{\bf g}, A, E_{\bf g}, \{{\rm Id}_{E_{\bf g}}\})=N$. 

En vertu du lemme \ref{lemmsepval}, il existe une famille finie
~$(b_{{\bf g},N,i})_{i\in I_{{\bf g},N}}$ d'éléments 
de~$K(T_1,\ldots, T_n)$ et une formule du premier ordre~$\Phi_{{\bf g},N}$ dans le langage des corps valués, sans paramètres ni quantificateurs, possédant les propriétés suivantes : 

$\bullet$ l'ensemble des variables libres de~$\Phi_{{\bf g},N}$ est indexé par~$I_{{\bf g},N}$, et chacune d'elles vit dans le corps valué ; 

$\bullet$  pour tout~$H\in \mathsf D_{|K^*|}$ et tout~${\bf h}\in H^n$,
l'énoncé~$\Phi_{\bf g,N}(b_{{\bf g},N,i})_{i\in I_{{\bf g},N}}$ est vrai dans~$(k(T_1,\ldots,T_n),\eta_{\bf h})$ si et seulement si l'on a~$\got n(\eta_{\bf h}, A, E_{\bf g}, \{{\rm Id}_{E_{\bf g}}\})=N$. En
particulier,~$\Phi_{\bf g,N}(b_{{\bf g},N,i})_{i\in I_{{\bf g},N}}$ est vrai dans~$(k(T_1,\ldots,T_n),\eta_{\bf g})$. 

Compte-tenu de la définition même de~$\eta_{\bf h}$, il existe une formule~$\Psi_{{\bf g}, N}$ en~$n$ variables, qui est du premier ordre dans le langage des groupes abéliens ordonnés et à paramètres dans~$|K^*|$, telle que pour tout~$H\in \mathsf D_{|K^*|}$ et tout~${\bf h}\in H^n$, l'énoncé~$\Phi_{\bf g,N}(b_{{\bf g},N,i})_{i\in I_{{\bf g},N}}$ soit vrai
dans~$(K(T_1,\ldots, T_n),\eta_{\bf h})$ si et seulement si~$\Psi_{{\bf g},N}({\bf h})$ est vraie. En particulier, 
$\Psi_{{\bf g},N}({\bf g})$ est vraie.  

Par le théorème de compacité, il existe un ensemble {\em fini}~$\{{\bf g}_{N,j}\}_{j\in J_N}$ où chaque~${\bf g}_{N,j}$ est un~$n$-uplet d'éléments d'un 
groupe appartenant à~$\mathsf D_{|K^*|}$, tel que~$\bigvee \Psi_{{\bf g}_{N,j},N}({\bf g})$ vaille pour tout~$G\in \mathsf D_{|K^*|}$ et tout~${\bf g}\in \sch D_N(G)$. Cela entraîne que la réunion des~$E_{{\bf g}_{N,j}}$ pour~$j\in J_N$ sépare les prolongements de~$\eta_{\bf g}$ à~$A$ pour tout~$G\in \mathsf D_{|K^*|}$ et tout~${\bf g}\in \sch D_N(G)$. 

\medskip
Il s'ensuit que la réunion des~$E_{{\bf g}_{N,j}}$, où~$N$ parcourt l'ensemble des entiers compris entre~$1$ et~$\dim{K({\bf T})}A$,  et où~$j$ parcourt~$J_N$ pour tout~$N$, 
sépare les prolongements de~$\eta_{\bf g}$ à~$L$ pour tout~$G\in \mathsf D_{|K^*|}$ et tout~${\bf g}\in G^n$. 
Ainsi, ii) est vraie. 

\medskip
Supposons maintenant que ii) est vraie. Soit~$G\in \mathsf D_{|K^*|}$ et soit~${\bf g}\in G^n$. Le nombre de prolongements de~$\eta_{\bf g}$ à~$A$ est égal à~$\got n(\eta_{\bf g}, A, E, \{{\rm Id}_E\})$ ; l'assertion i) découle dès lors du lemme \ref{lemmsepval}, ce qui termine la preuve lorsque~$\sch C=\{K\}$.

\medskip
\trois{schcgen} {\em Le cas général.} Supposons que i) est  vraie. Par le cas particulier déjà traité, il existe un sous-ensemble fini~$E$ de~$A$
qui sépare les plongements de~$\eta_{\bf g}$ à~$A$ pour tout~$G\in \mathsf D_{|K^*|}$ et tout~${\bf g}\in \sch G^n$. On a donc pour tout entier~$N$ et tout~$G\in \mathsf D_{|K^*|}$ l'égalité~$$\sch D_N(G)=\{{\bf g}\in G^n, \got n(\eta_{\bf g}, A, E,\{{\rm Id}_E\})=N\}.$$ 

Soit~$\Lambda\in \sch C$, soit~$G\in \mathsf D_{|\Lambda^*|}$, soit~${\bf g}\in G^n$ et soit
$N$ un entier. D'après i),~$\sch D_{\Lambda, N}(G)$ est égal à~$\sch D_N(G)$, et donc à~$$\{{\bf g}\in G^n, \got n(\eta_{\bf g}, A, E,\{{\rm Id}_E\})=N\}.$$ Le lemme \ref{lemmsepval} assure que ce dernier ensemble coïncide avec~$$\{{\bf g}\in G^n, \got n(\eta_{\Lambda, \bf g}, \Lambda({\bf T})\otimes_{K({\bf T})}A, E,\{{\rm Id}_E\})=N\}.$$ Ceci valant pour tout~$N$, l'ensemble~$E$ sépare les prolongements de 
$\eta_{\Lambda, \bf g}$ à~$\Lambda({\bf T})\otimes_{K({\bf T})}A$. 

\medskip
Supposons maintenant que ii) est  vraie. Par le cas particulier déjà traité,~$\sch D_N$ est définissable. Soit~$\Lambda\in \sch C$, soit~$G\in \mathsf D_{|\Lambda^*|}$ et soit~${\bf g}\in G^n$. D'après ii),~$\sch D_{\Lambda, N}(G)=\{{\bf g}\in G^n, \got n(\eta_{\Lambda, \bf g}, \Lambda({\bf T}) \otimes_{K({\bf T})} A, E,\{{\rm Id}_E\})=N\}.$ Le lemme \ref{lemmsepval} assure que ce dernier ensemble coïncide avec~$\{{\bf g}\in G^n, \got n(\eta_{\bf g}, A, E,\{{\rm Id}_E\})=N\}$, qui n'est autre que~$\sch D_N(G)$ par définition de~$E$. Ceci achève la démonstration.~$\Box$

\deux{propvalg} {\bf Théorème.} {\em Soit~$K$ un corps valué et soit~$A$ une~$K(T_1,\ldots, T_n)$-algèbre finie.

\medskip
1) Il existe un sous-ensemble fini de~$A$ qui sépare les prolongements de~$\eta_{\bf g}$ à~$A$ pour tout~$G\in \mathsf D_{|K^*|}$ et tout~${\bf g}\in G^n$. 

\medskip
2) Il existe une extension valuée finie séparable~$K_0$ de~$K$ et un sous-ensemble fini de~$K_0\otimes_K A$ qui sépare les prolongements de~$\eta_{\Lambda, \bf g}$ à~$\Lambda({\bf T})\otimes_{K({\bf T})}A$ pour toute extension valuée~$\Lambda$ de~$K_0$,  pour tout~$G\in \mathsf D_{|\Lambda^*|}$ et tout~${\bf g}\in G^n$. }

\medskip
{\em Démonstration.} On procède en plusieurs étapes. {\em On suppose tout d'abord que $K$ est algébriquement clos}. Il suffit alors de montrer~2), avec $K_0$ évidemment
égal à $K$. On procède par récurrence sur $n$.

\trois{casplalgclosnzero} {\em Le cas où ~$n=0$.} Pour toute extension valuée $\Lambda$ de $K$, le corps résiduel de tout point
de~$\spec (\Lambda\otimes_KA)$ est égal à $\Lambda$, et l'ensemble
des valuations de~$\Lambda\otimes_KA$ prolongeant celle de~$\Lambda$
s'identifie donc naturellement à {\em l'ensemble} fini~$\spec \Lambda\otimes_KA$,
lui-même en bijection naturelle avec~$\spec A$.
L'assertion 2) est alors évidente : il suffit de prendre un ensemble
d'idempotents de $A$ séparant les points de $\spec A$. 

\medskip
\trois{casplalgclosnqque1} {\em Utilisation d'un théorème de finitude de Hrushovski et Loeser pour exhiber une 
première famille
finie séparante.} On suppose~$n>0$, et le résultat vrai au rang~$n-1$. Choisissons 
un~$K$-schéma affine et de type fini~$\sch X$ purement
de dimension~$n$, 
muni d'un~$K$-morphisme génériquement fini~$\sch X\to \Aff^n_K$ 
dont la fibre générique s'identifie à la
flèche~$\spec A\to \spec K(T_1,\ldots, T_n)$. On considère~$\Aff^n_K$, et partant~$\sch X$, 
comme un~$\Aff^{n-1}_K$-schéma
{\em via} la projection sur les~$n-1$ premières coordonnées.

\medskip
Soit~$\widehat {\Aff^n_K/\Aff^{n-1}_K}$ (resp.~$\widehat {\sch X/\Aff^{n-1}_K}$ ) le foncteur de~$\mathsf C_K$ vers les ensembles qui associe à~$\Lambda\in \mathsf C_K$ l'ensemble des~$\Lambda$-types stablement dominés situés sur~$\Aff^n_\Lambda$ (resp.~$\sch X_\Lambda$) qui induisent un~$\Lambda$-point (ou un type simple, si l'on préfère) sur~$\Aff^{n-1}_\Lambda$.  Pour la définition de type stablement dominé, nous renvoyons par exemple au paragraphe 2.5 de l'article \cite{hl} de Hrushovski et Loeser. 

\medskip
En vertu d'un résultat de Hrushovski et Loeser (c'est le lemme 7.1.3 de {\em loc. cit.}, qui repose en dernière analyse
sur le théorème de Riemann-Roch pour les courbes),
les foncteurs~$\widehat {\Aff^n_K/\Aff^{n-1}_K}$ et~$\widehat{\sch X/\Aff^{n-1}_K}$ sont~$K$-définissables. Pour tout~$\Lambda\in \mathsf C_K$, tout~$x\in \Lambda^{n-1}$, tout~$G\in \mathsf D_{|\Lambda^*|}$ et tout~$g\in G$, on note~$\eta_{x,g}$ le type sur~$\Aff^n_\Lambda$ situé au-dessus de~$x$ et défini par la valuation~$\eta_{\Lambda, g}$ sur~$\Lambda(T_n)$. 
Le foncteur~$$\sk:=\Lambda\mapsto \{\eta_{x,g}\}_{x\in \Lambda^{n-1}, g\in |\Lambda^*|}$$ est un sous-foncteur de~$\widehat {\Aff^n_K/\Aff^{n-1}_K}$ qui est~$K$-définissable : il est par sa définition même naturellement isomorphe au produit de~$\Aff^{n-1}_K$ par le groupe des valeurs. Posons~$${\bf T}=\widehat{\sch X/\Aff^{n-1}_K}\times_{\widehat{\Aff^n_K/\Aff^{n-1}_K}}\sk.$$

Le morphisme de foncteurs~${\bf T}\to \sk$ est~$K$-définissable et à fibres finies. Soit~$\Lambda\in \mathsf C_K$, soit~$x\in \Lambda^{n-1}$ et soit~$g\in| \Lambda^*|$. Le couple~$(x,g)$ définit un point de~$\sk(\Lambda)$ ; tout antécédent de ce point dans~${\bf T}(\Lambda)$ est algébrique sur~$\overline{K(x)} \cup\{ g\}$, où~$\overline {K(x)}$ est la fermeture algébrique de~$K(x)$ dans~$\Lambda$. Il est donc,  en vertu du lemme 3.4.12 de \cite{hhmelim}, {\em définissable} sur~$\overline{K(x)} \cup\{g\}$. Autrement dit, il s'écrit comme l'image de~$(x,g)$ par une fonction qui est~$\overline{K(x)}$-définissable, et donc définissable sur une extension finie galoisienne de~$K(x)$ ; on peut choisir une telle extension qui convienne pour tous les antécédents de~$(x,g)$. 

\medskip 
Par compacité, on en déduit qu'il existe : 

\medskip
$\bullet$ une famille finie~$(\sch U_i)_i$ de sous-schémas localement fermés, intègres et affines de~$\Aff^{n-1}_K$ ; 

$\bullet$ pour tout~$i$, un revêtement fini galoisien connexe~$\sch V_i\to \sch U_i$ et un sous-foncteur
définissable~$V_i$ de~$\sch V_i$, de sorte que les images des~$V_i$ recouvrent~$\Aff^{n-1}_K$ ; 

$\bullet$ pour tout~$i$, une famille finie~$(\sigma_{ij})_j$ de sections~$K$-définissables de l'application~${\bf T}\times_{\Aff^{n-1}_K}V_i\to \sk \times_{\Aff^{n-1}_K}V_i$ dont les images recouvrent ~${\bf T}\times_{\Aff^{n-1}_K}V_i$.

\medskip
Pour tout couple d'indices~$(i,j)$, notons~${\bf T}_{ij}$ le sous-foncteur~$K$-définissable de~${\bf T}\times_{\Aff^{n-1}_K}V_i$ égal à l'image de~$\sigma_{ij}$. 

\medskip
Pour établir la définissabilité de~$\widehat{\sch X/\Aff^{n-1}_K}$, Hrushovski et Loeser montrent l'existence d'un sous-$\sch O_{\Aff^{n-1}_K}$-module de type fini~$\sch E$ de~$\sch O_{\sch X}$ tel que~$\widehat{\sch X/\Aff^{n-1}_K}$ s'identifie comme suit à un sous-foncteur définissable du foncteur ${\cal S}_{\sch E}$ des «réseaux relatifs de~$\sch E$ sur~$\Aff^{n-1}_K$» : si $\Lambda\in \mathsf C_K$, si $x\in \Lambda^{n-1}$, et si~$y$ est un point de~$\widehat{\sch X/\Aff^{n-1}_K}(\Lambda)$ situé au-dessus de $x$, on l'envoie sur le point 
de~${\cal S}_{\sch E}(\Lambda)$, situé au-dessus de $x$, qui correspond au réseau de $\sch E_x$ constitué des fonctions $f$ telles que $|f(y)|\leq 1$. 

\medskip
Fixons~$i$. Par ce qui précède, il existe pour tout~$j$ un sous-foncteur fini et $K$-définissable 
de~$\sch O_{\sch X\times_{\Aff^{n-1}_K}\sch V_i}(\sch X\times_{\Aff^{n-1}_K}\sch V_i)$ constitué de fonctions 
dont les valuations permettent 
de détecter l'appartenance à~${\bf T}_{ij}$ d'un point de~${\bf T}\times_{\Aff^{n-1}_K}V_i$ ; comme $K$ est algébriquement clos, 
ce sous-foncteur fini est induit par un {\em sous-ensemble} fini
de~$\sch O_{\sch X\times_{\Aff^{n-1}_K}\sch V_i}(\sch X\times_{\Aff^{n-1}_K}\sch V_i)$. 
Il existe donc un sous-ensemble 
fini~$\Theta_i$ de~$\sch O_{\sch X\times_{\Aff^{n-1}_K}\sch V_i}(\sch X\times_{\Aff^{n-1}_K}\sch V_i)$ tel 
que les fonctions~$|f|$, pour~$f$ parcourant~$\Theta_i$, séparent universellement 
les points des fibres de~${\bf T}\times_{\Aff^{n-1}_K}V_i\to \sk \times_{\Aff^{n-1}_K}V_i$. 

\medskip
Soit~$I$ l'ensemble des indices~$i$ tels que~$\sch U_i$ soit un ouvert de~$\Aff^{n-1}_K$. Choisissons un ouvert affine non vide de~$\Aff^{n-1}_K$ et un revêtement fini, connexe et galoisien~$\sch V$ dudit ouvert tel que~$\sch V\to \Aff^{n-1}_K$ se factorise par le schéma~$\sch V_i$ pour tout~$i\in I$. Posons~$\sch Y=\sch X\times_{\Aff^{n-1}_K}\sch V$. 
Pour tout~$i\in I$, les éléments de~$\Theta_i$ peuvent être vus comme 
appartenant à~$\sch O_{\sch Y}(\sch Y)$ ; 
on note~$\Theta$ le sous-ensemble fini~$\bigcup\limits_{i\in I} \Theta_i$ de~$\sch O_{\sch Y}(\sch Y)$. 

\medskip
Avant de poursuivre, fixons quelques notations. 
Pour toute extension~$\Lambda$ du corps algébriquement clos~$K$, on désigne 
par~$F_\Lambda$ le
corps des fonctions du schéma intègre~$\sch V_\Lambda$. Le schéma $\Aff^n_\Lambda\times_{\Aff^{n-1}\Lambda}\sch V_\Lambda$ s'identifie à~$\Aff^1_{\sch V_\Lambda}$ et est
donc intègre, de corps des fonctions~$F_\Lambda(T_n)$. La flèche~$\sch X\to \Aff^n_K$ induit une flèche génériquement finie~$\sch Y\to \Aff^n_K\times_{\Aff^{n-1}_K}\sch V$,
d'où par extension des scalaires une flèche génériquement finie $\sch Y_\Lambda\to \Aff^n_\Lambda\times_{\Aff^{n-1}_\Lambda}\sch V_\Lambda$. Sa fibre générique est de la forme $\spec F'_\Lambda$, 
où~$F'_\Lambda$ est une certaine~$F'_\Lambda(T_n)$-algèbre finie. Si~$\Lambda=K$, on écrira~$F$ et~$F'$ 
au lieu de~$F_\Lambda$ et~$F'_\Lambda$.

\medskip
Soit~$\Lambda$ une extension valuée de~$K$. Fixons un~$K$-prolongement~$\val .$ 
de la valuation de~$\Lambda$ à~$\Lambda(T_1,\ldots, T_{n-1})$, 
et soit~$G\in \mathsf D_{|\Lambda(T_1,\ldots,T_{n-1})^*|}$. Soit~$g\in G$, et soit~$\Lambda^\sharp$ appartenant à~$\mathsf C_{\Lambda(T_1,\ldots, T_{n-1}), G}$. Le plongement~$\Lambda(T_1,\ldots, T_{n-1})\hookrightarrow \Lambda^\sharp$ définit un~$\Lambda^\sharp$-point~$x$ de~$\Aff^{n-1}_K$. Soit~$\xi$ le point générique de~$\Aff^n_K\times_{\Aff^{n-1}_K}x\simeq \Aff^1_{\Lambda^\sharp}$.

\medskip
La valuation~$\eta_{\Lambda^\sharp, g}$ sur~$\kappa(\xi)=\Lambda^\sharp(T_n)$ définit un point de~$\sk(\Lambda^\sharp)$, qui s'identifie à~$(x,g)$ {\em via} l'isomorphisme mentionné plus haut. Par construction, il existe un antécédent~$z$ de~$x$ sur~$\sch V(\Lambda^\sharp)$ tel que les fonctions~$|f|$, où~$ f$ parcourt~$\Theta$, séparent les antécédents de~$(z,g)$ sur~$({\bf T}\times_{\Aff^{n-1}_K }\sch V)(\Lambda^\sharp)$. 

\medskip
Comme~$x$ est situé au-dessus du point générique de~$\Aff^{n-1}_\Lambda$, le point~$z$ induit une valuation sur~$F_\Lambda$, que l'on voit
comme un corps valué par ce biais. Le point~$(z,g)$ de~$\sk(\Lambda^\sharp)$ induit alors la valuation~$\eta_{F_\Lambda,g}$ sur~$F_\Lambda(T_n)$. 

\medskip
Soit~$\Omega$ une clôture algébrique de~$\Lambda^\sharp(T_n)$, munie d'un prolongement de~$\eta_{\Lambda^\sharp,g}$. 
Tout prolongement de~$\eta_{F_\Lambda(T_n),g}$ à~$F_\Lambda'$ est induit par un morphisme de~$F_\Lambda'$ dans~$\Omega$, 
c'est-à-dire par un morphisme
de~$F_\Lambda'\otimes_{F_\Lambda(T_n)}\Lambda^\sharp(T_n)$ dans~$\Omega$,
c'est-à-dire encore par l'un des antécédents de~$(z,g)$ sur~$({\bf T}\times_{\Aff^{n-1}_K}\sch V)(\Lambda^\sharp)$ 
 ; il s'ensuit que les prolongements de~$\eta_{F_\Lambda(T_n),g}$ à~$F'_\Lambda$ sont séparés par les éléments de~$\Theta$.

\trois{casplalgclosnqque2} {\em Utilisation de l'hypothèse de récurrence et conclusion dans le cas algébriquement clos.} L'hypothèse de récurrence
assure qu'il existe un sous-ensemble fini~$\Upsilon$ de~$F$ qui 
possède la propriété suivante : pour toute extension valuée~$\Lambda$ de~$K$,  pour tout~$G\in \mathsf D_{|\Lambda^*|}$ et tout~${\bf g}\in G^{n-1}$, l'ensemble~$\Upsilon$ sépare les prolongements de~$\eta_{\Lambda, \bf g}$ à~$F_\Lambda$. 
Soit~$H$ le groupe de Galois de~$F$ sur~$K(T_1,\ldots, T_{n-1})$. 

\medskip
Soit~$\Lambda$ une extension valuée de~$K$, soit~$G\in \mathsf D_{|\Lambda^*|}$ et soit~${\bf g}\in G^n$. Notons~${\bf g}'$ le~$(n-1)$-uplet constitué des~$n-1$ premières coordonnées de~$\bf g$, et~$g$ sa dernière coordonnée. On a alors~$\eta_{\Lambda, {\bf g}}=\eta_{(\Lambda(T_1,\ldots,T_{n-1}),\eta_{\Lambda, {\bf g}'}),g}$. 

Nous allons travailler avec différents anneaux qui s'organisent selon le diagramme cartésien suivant, dont les flèches horizontales sont galoisiennes de groupe~$H$ : 

$$\diagram \Lambda({\bf T})\otimes_{K({\bf T})}A \rto &F'_\Lambda\\ \Lambda(T_1,\ldots, T_{n-1},T_n)\uto\rto &F_\Lambda(T_n)\uto \enddiagram.$$ 

\medskip
L'ensemble des prolongements de~$\eta_{\Lambda, {\bf g}}$ à ~$\Lambda({\bf T})\otimes_{K({\bf T})}A$ s'identifie au quotient par~$H$ de l'ensemble des prolongements de ~$\eta_{\Lambda, {\bf g}}$ à~$F'_\Lambda$ (\ref{valanneauextgal}). 

\medskip
Par ailleurs, donnons-nous deux prolongements {\em distincts} de~$\val ._1$ et~$\val ._2$ de~$\eta_{\Lambda, {\bf g}}$ à~$F'_\Lambda$. On est alors dans l'un des deux cas suivants.

\medskip
{\em Premier cas.} Les restrictions de~$\val ._1$ et~$\val ._2$ à~$F_\Lambda$ diffèrent ; comme ces deux restrictions sont des prolongements de~$\eta_{\Lambda, {\bf g}'}$, l'ensemble~$\Upsilon$ sépare~$\val ._1$ et~$\val ._2$. 

\medskip
{\em Second cas.} Les restrictions de~$\val ._1$ et~$\val ._2$ à~$F_\Lambda$ coïncident. Si l'on note~$\val .'$ leur restriction commune à~$F_\Lambda$, les restrictions de~$\val ._1$ et~$\val ._2$ à~$F_\Lambda(T_n)$ sont alors toutes deux égales, en vertu de \ref{corextgauss}, à~$\eta_{(F_\Lambda,\val .'),g}$.

\medskip
D'après ce qu'on a vu au~\ref{casplalgclosnqque1}, il existe un prolongement~$\val .''$ 
de~$\eta_{\Lambda, {\bf g}'}$
à~$F_\Lambda$ tel que les prolongements 
de~$\eta_{F_\Lambda,\val .''}$ à~$F'_\Lambda$ soient séparés par les éléments de~$\Theta$.
Il existe~$h\in H$ 
tel que~$\val .'=\val .'' \circ h$. On peut reformuler
ce qui précède en écrivant que tous les prolongements 
de~$\eta_{(F_\Lambda,\val .'),g)}\circ h$ à~$F'_\Lambda$ sont séparés par~$\Theta$. Par conséquent,~$\val ._1$ et~$\val ._2$ sont séparées par~$h.\Theta$.

\medskip
Soit~$\Xi$ la réunion des orbites sous~$H$ des éléments de~$\Theta$ et~$\Upsilon$. Par ce qui précède,~$\Xi$ sépare les prolongements de~$\eta_{\Lambda, {\bf g}}$ 
à~$F'_\Lambda$. Le nombre de prolongements de~$\eta_{\Lambda, {\bf g}}$ à~$ \Lambda({\bf T})\otimes_{K({\bf T})}A$ 
est en conséquence égal à~$$\got n(\eta_{\Lambda, {\bf g}},F'_\Lambda, \Xi,\Sigma),$$ où~$\Sigma$ désigne le groupe de permutations de~$\Xi$ induit par~$H$.

\medskip

%
%
%
En vertu du lemme \ref{lemmsepval}, ceci implique que la première des deux assertions équivalentes de la proposition \ref{equivprogauss} est satisfaite, en prenant pour~$\sch C$ la classe de toutes les extensions valuées de~$K$. Il s'ensuit que la seconde est satisfaite pour cette même classe, ce qui constitue exactement l'énoncé souhaité et achève la preuve lorsque~$K$ est algébriquement clos.

\trois{casnzeropl} {\em Preuve dans le cas général ; on ne suppose plus que $K$ est algébriquement clos. } Soit~$\overline K$ une clôture algébrique de~$K$ ; fixons un prolongement de la valuation de~$K$ à~$\overline K$. En vertu du cas algébriquement clos déjà traité, il existe un sous-ensemble fini~$E$ de~$\overline K\otimes_K A$ qui sépare les prolongements de~$\eta_{\Lambda, \bf g}$ à~$\Lambda({\bf T})\otimes_{K({\bf T})}A$ pour toute extension valuée~$\Lambda$ de~$\overline K$,  pour tout~$G\in \mathsf D_{|\Lambda^*|}$ et tout~${\bf g}\in G^n$. 

\medskip
Quitte à élever les éléments de~$E$ à une puissance convenable de l'exposant caractéristique de~$K$ (ce qui ne modifie pas leur capacité à séparer les valuations), on peut supposer qu'ils sont contenus dans ~$K_{\rm sep}\otimes_K A$, où~$K_{\rm sep}$ est la fermeture séparable de~$K$ dans~$\overline K$. Il existe par conséquent une sous-extension finie séparable~$K_0$ de~$\overline K$ telle que~$E\subset K_0\otimes_K A$. Soit~$\Lambda$ une extension valuée de~$K_0$, soit~$G$ appartenant à $\mathsf D_{|\Lambda^*|}$ et soit~${\bf g}\in G^n$. Soient~$\val .'$ et~$\val .''$ deux prolongements distincts de~$\eta_{\Lambda, \bf g}$ à~$\Lambda({\bf T})\otimes_{K({\bf T})}A$.

\medskip
D'après~\ref{valanneauextsuite} (appliqué à la~$K_0$-algèbre~$\Lambda$ et à la~$K_0$-extension galoisienne~$\overline K$), il existe une valuation sur $\overline K\otimes_{K_0}\Lambda$ prolongeant les valuations fixées sur~$\overline K$ et~$\Lambda$ ; cette valuation est induite par une valuation~$\val .$ sur un quotient~$\Lambda^\sharp$ de~$\overline K\otimes_{K_0}\Lambda$. D'après~\ref{valanneauextsuite}, les valuations~$\val .'$ et~$\val .''$ s'étendent en deux valuations de~$\Lambda^\sharp({\bf T})\otimes_{K({\bf T})}A$, encore notées~$\val .'$ et~$\val .''$, dont les restrictions à~$\Lambda^\sharp$ sont toutes deux égales à~$\val .$. Comme les restrictions de~$\val .'$ et~$\val .''$ à~$\Lambda({\bf T})$ sont toutes deux égales à~$\eta_{\Lambda, \bf g}$, il résulte de~\ref{corextgauss} que les restrictions de~$\val .'$ et~$\val .''$ à~$\Lambda^\sharp({\bf T})$ sont toutes deux égales à~$\eta_{\Lambda^\sharp, \bf g}$. Par ailleurs, les valuations~$\val .'$ et~$\val .''$ de~$\Lambda^\sharp({\bf T})\otimes_{K({\bf T})}A$ sont distinctes, puisque leurs restrictions à~$\Lambda({\bf T})\otimes_{K({\bf T})}A$ le sont. Par choix de~$E$, il existe donc~$x\in E$ tel que~$\val x '\neq \val  x''$. Ainsi, $E$ sépare les extensions de la valuation de~$\Lambda$ à~$\Lambda({\bf T})\otimes_{K({\bf T})}A$, ce qui montre 2).

\medskip
Il reste à montrer 1). Soit~$\Sigma$ l'ensemble des permutations de~$E$ induites par l'action de~${\rm Gal}\;(\overline K/K)$. Pour tout entier~$N$, on note~$\sch D_N$ le foncteur qui envoie un groupe~$G\in \mathsf D_{|K^*|}$ sur le sous-ensemble de~$G^n$ formé des~$n$-uplets~$\bf g$ tels que~$\eta_{\bf g}$ admette exactement~$N$ prolongements à~$A$. 

\medskip
Soit~$G\in \mathsf D_{|K^*|}$ et soit~${\bf g}\in G^n$. Comme~$E$ sépare les prolongements de~$\eta_{\overline K, \bf g}$ à~$\overline K\otimes_KA$, le nombre de prolongements de la valuation~$\eta_{\bf g}$ à~$A$ est égal en vertu de~\ref{valanneauextsuite} à~$\got n(\eta_{\overline K, \bf g}, \overline K\otimes_K A, E, \Sigma)$. Il résulte alors du lemme \ref{lemmsepval} que le sous-foncteur~$\sch D_N$ de~$\GG_{|K^*|}^n$ est définissable pour tout~$N$. En conséquence, la première des deux assertions équivalentes de la proposition \ref{equivprogauss} est satisfaite pour la classe triviale~$\sch C=\{K\}$ ; la seconde l'est donc aussi pour cette même classe, et c'est précisément ce qu'il fallait démontrer.~$\Box$

\section{Tropicalisations globale et locale d'un espace~$k$-analytique}

\deux{lemmebord} {\bf Lemme.} {\em Soit~$X$ un espace~$k$-affinoïde~$\Gamma$-strict et de dimension~$d$. Il existe une famille finie~$(p_1,\ldots, p_n)$ de morphismes de~$X$ vers~$\Aff^{1,\rm an}_k$, et pour tout~$i$ un réel~$s_i$ appartenant à~$\sqrt{\Gamma\cdot |k^*|}$, tels que~$\partial X=\bigcup p_i\inv(\eta_{s_i})$. De plus, pour tout~$i$, l'espace~$\hres(\eta_{s_i})$-affinoïde~$ p_i\inv(\eta_{s_i})$ est~$\Gamma$-strict et de dimension~$\leq d-1$.}

\medskip
{\em Démonstration.}  D'après la définition de la réduction d'un germe à la Temkin, 
il existe une famille finie~$g_1,\ldots, g_r$ d'éléments non
nilpotents de~$\sch A$, de rayons spectraux respectifs notés
$s_1,\ldots, s_r$, tels que pour tout~$x\in X$ la réduction
$\widetilde{(X,x)}$ soit égale à~$\PP_{\red{\hres(x)}/\red
k}\{\widetilde{g_i(x)}^{s_i}\}_i$. Pour tout~$i$, on note~$p_i$ le morphisme~$X\to \Aff^{1,\rm an}_k$ induit par~$g_i$. 

Le point~$x$ appartient
donc à~$\partial X$ si et seulement si il existe un indice~$i$
et une~$\red k$-valuation graduée~$\val .$ sur~$\red{\hres(x)}$ telle que~$\val{\widetilde{g_i(x)}^{s_i}}>1$. Cette dernière condition
équivaut à demander que ~$|g_i(x)|=s_i$ et que~$\widetilde{g_i(x)}$
soit transcendant sur~$\red k$ (\ref{valtrvsalg}) ; en vertu de \ref{transeteta}, elle est
satisfaite si et seulement si~$x\in p_i^{-1}(\eta_{s_i})$. En conséquence,~$\partial X=\bigcup p_i^{-1}(\eta_{s_i})$.

\medskip
Fixons~$i$ ; nous allons prouver que~$p_i^{-1}(\eta_{s_i})$
est~$\Gamma$-strict et de dimension au plus~$d-1$. 

Pour tout point~$z$ de ~$p_i^{-1}(\eta_{s_i})$
on a~$d\geq d_k(z)=d_{\hres(\eta_{s_i})}(z)+1$, et donc
$d_{\hres(\eta_{s_i})}(z)\leq d-1$ ; il s'ensuit que la dimension
$\hres(\eta_{s_i})$-analytique de ~$p_i^{-1}(\eta_{s_i})$ est
inférieure ou égale à~$d-1$.

Comme~$s_i$ est le rayon spectral d'un élément de~$\sch A$,
il appartient à~$\sqrt{\Gamma\cdot|k^*|}$. Par conséquent,~$\sqrt{\Gamma\cdot|\hres(\eta_i)^*|}=\sqrt{\Gamma\cdot|k^*|}$. Il suffit dès lors de vérifier que l'algèbre~$\hres(\eta_i)$-affinoïde des fonctions analytiques sur~$p_i^{-1}(\eta_{s_i})$ est~$\Gamma$-stricte. 
Mais cela découle du fait qu'elle s'identifie au 
quotient~$$\sch A\hotimes_k \hres(\eta_{s_i}) /(g_i-T(\eta_{s_i})).~\Box$$

\deux{tropicglob} {\bf Théorème (le cas global).} {\em Soit~$X$ un espace ~$k$-analytique compact 
et~$\Gamma$-strict de dimension~$d$, et soit~${\bf f}=(f_1,\ldots, f_n)$ une famille de fonctions analytiques inversibles sur~$X$ ; notons~$|{\bf f}|$  l'application~$(|f_1|,\ldots, |f_n|)$ de~$X$ vers~$\rst n$. Posons~$c=(\QQ,\sqrt{|k^*|\cdot \Gamma})$. 

\medskip
1) L'image~$|{\bf f}|(X)$ est un~$c$-polytope de~$\rst n$, qui est de dimension inférieure ou égale à~$d$. 

2) L'image~$|{\bf f}|(\partial X)$ est contenue dans un~$c$-polytope  de~$\rst n$ qui est de dimension inférieure ou égale à~$d-1$ ; lorsque~$X$ est affinoïde ~$|{\bf f}|(\partial X)$ est elle-même un~$c$-polytope de~$\rst n$ de dimension inférieure ou égale à~$d-1$.}

\medskip
{\em Démonstration.} La seconde assertion se déduira de la première.

\trois{preuve1redalg} {\em Preuve de 1) : réduction à une situation algébrique.} On peut supposer que~$X$ est affinoïde et irréductible, et l'on procède alors par récurrence sur~$d$. Si~$d=0$ il n'y a rien à démontrer. Supposons que~$d>0$ et que 1) est vraie en dimension~$<d$.

Le théorème à établir est purement topologique. On peut donc supposer~$k$ parfait (en étendant les scalaires au complété de la clôture radicelle de~$k$), et~$X$ réduit (en le remplaçant par~$X_{\rm red}$). L'espace~$X$ est alors génériquement quasi-lisse. Soit~$Y$ son lieu de non-quasi-lissité. C'est un fermé de Zariski de~$X$ dont toutes les composantes irréductibles sont de dimension~$<d$. Par hypothèse de récurrence,~$|{\bf f}|(Y)$ est un~$c$-polytope de~$\rst n$ de dimension~$\leq d-1$. Son image réciproque~$ |{\bf f}|\inv (|{\bf f}|(Y))$ est un domaine analytique~$\Gamma$-strict de~$X$ contenant~$Y$. Comme~$Y\hookrightarrow X$ est sans bord,~$Y$ est contenu dans l'intérieur topologique de~$|{\bf f}|\inv (|{\bf f}|(Y))$. Ce dernier apparaît ainsi comme un voisinage de~$Y$ dont l'image par~$|{\bf f}|$ est un~$c$-polytope de dimension~$\leq d-1$. 

Soit~$x\in X-Y$. Comme~$X$ est quasi-lisse en~$x$, le lemme \ref{algquasil} assure l'existence d'un voisinage affinoïde~$V$ de~$x$ dans~$X$, et d'une~$k$-variété affine~$\sch X$ de dimension~$d$, telle que~$V$ s'identifie à un domaine affinoïde~$\Gamma$-strict de~$\sch X\an$ (on peut imposer à~$\sch X$ d'être lisse, mais nous ne nous en servirons pas). 

Par compacité de~$X$, on se ramène finalement au cas où~$X$ lui-même est est un domaine affinoïde~$\Gamma$-strict de l'analytification~$\sch X\an$ d'une~$k$-variété algébrique~$\sch X$ affine de dimension~$d$. 

\medskip
Écrivons~$\sch X\simeq \spec A$ où~$A$ est une~$k$-algèbre de type fini, et soit~$(a_1,\ldots, a_r)$ un système de générateurs de~$A$. Il existe~$R\in \Gamma\cdot |k^*|$ tel que~$X$ soit contenu dans le domaine affinoïde~$W$ de~$\sch X\an$ défini par les conditions~$|a_i|\leq R$. En vertu du théorème de Gerritzen-Grauert (\ref{gerrgrauertgam}), le domaine~$X$ est réunion finie de domaines~$\Gamma$-rationnels de~$W$ ; il suffit dès lors de traiter le cas où~$X$ est lui-même un domaine~$\Gamma$-rationnel de~$W$.

Il est alors défini par une conjonction d'inégalités de la forme~$|g_i|\leq \lambda_i |h|$ où les~$\lambda_i$ appartiennent à~$\Gamma$ et où~$(g_1,\ldots, g_r,h)$ engendre l'idéal des fonctions analytiques sur~$W$. Cette dernière condition implique que~$h$ ne s'annule pas sur~$X$ (tout zéro de~$h$ sur~$V$ est un zéro commun à~$h$ et aux~$g_i$), et l'on peut donc, par compacité, rajouter une condition de la forme~$|h|\geq r$ à la définition de~$X$, où~$r$ est un élément convenable de~$\sqrt{\Gamma\cdot |k^*|}$. Il est dès lors immédiat que~$V$ ne change pas si l'on perturbe un peu les~$g_i$ et~$h$ ; on peut par conséquent supposer que toutes ces fonctions appartiennent à~$A$. 

L'anneau des fonctions de la forme~$a/h$, où~$a\in A$, est dense dans l'anneau des fonctions de~$X$. Cela permet de supposer, quitte à remplacer~$\sch X$ par~$D(h)$ et à approcher convenablement chacune des~$f_i$, que~$f_i\in A$ pour tout~$i$.

\trois{preuve1alg} {\em Nature~$c$-polytopale de~$|{\bf f}|(X)$ dans le cas algébrique.} Soit~$(F,|.|)$ une extension valuée de~$k$ (on ne suppose pas que le groupe ordonné~$|F^*|$ est de rang 1). On peut donner un sens naturel à la notation~$X(F)$, coïncidant avec sa signification usuelle lorsque~$F$ est un corps ultramétrique complet : c'est le sous-ensemble de~$\sch X(F)$ formé des points~$P$ tels que~$|g_i(P)|\leq \lambda_i|h(P)|$ pour tout~$i$. 

\medskip
Le théorème \ref{propeq} assure l'existence d'un sous-foncteur fermé et définissable~$\sch D$ de~$\GG_{\Gamma\cdot |k^*|}^n$, de dimension majorée par~$d$ (les notations et définitions sont celles de \ref{gencorpsval} {\em et sq.}) tel que~$$|{\bf f}|(X(F))=\sch D(|F^*|)\subset |F^*|^n$$ pour tout~$F\in \mathsf C_{k,\Gamma\cdot |k^*|}$. 

\medskip
Si~$F$ est une extension ultramétrique complète et algébriquement close de~$k$ munie d'un~$|k^*|$-isomorphisme~$|F^*|\simeq \RR^*_+$, tout point~$P$ de~$X(F)$ définit un point~$x$ de l'espace~$k$-affinoïde~$X$, et~$|{\bf f}|(x)=|{\bf f}|(P)\in \rst n$. 

Réciproquement si~$x\in X$ il existe une extension complète et algébriquement close~$F$ de~$\hres(x)$ et un~$|\hres(x)^*|$-isomorphisme~$|F^*|\simeq \RR^*_+$, et~$x$ est alors induit  par un~$F$-point canonique de~$X$. 

\medskip
En conséquence, le compact~$|{\bf f}|(X)$ est le sous-ensemble~$\sch D(\RR^*_+)$ de~$\rst n$, ce qui achève de prouver 1).

\trois{preuve2} {\em Démonstration de l'assertion 2}. Si~$X$ s'écrit comme une réunion finie~$\bigcup X_i$ où les~$X_i$ sont affinoïdes et~$\Gamma$-stricts, son bord est contenu dans la réunion des bords des~$X_i$. On peut donc supposer que~$X$ est affinoïde, auquel cas le résultat découle de l'assertion 1) déjà établie et du lemme \ref{lemmebord}.~$\Box$

\deux{theotroploc} {\bf Théorème (le cas local)}. {\em Soit~$X$ un espace ~$k$-analytique~$\Gamma$-strict, et soit~${\bf f}=(f_1,\ldots, f_n)$ une famille de fonctions analytiques inversibles sur~$X$ ; notons~$|{\bf f}|$  l'application~$(|f_1|,\ldots, |f_n|)$ de~$X$ vers~$\rst n$. Soit~$x$ un point de~$X$ ; posons~$\xi=|{\bf f}|(x)$, et notons~$d$ le degré de transcendance sur~$\red k$ du sous-corps gradué de~$\red {\hres(x)}$ engendré par~$\red k$ et les~$\red{f_i(x)}$. 

\medskip
\begin{itemize}

\item[1)] Il existe un voisinage~$k$-analytique compact et~$\Gamma$-strict~$U$ de~$x$ dans~$X$ possédant la propriété suivante : pour tout voisinage analytique compact~$V$ de~$x$ dans~$U$, les germes de polytopes~$(|{\bf f}|(U),\xi)$ et~$(|{\bf f}|(V),\xi)$ coïncident.

\medskip
De plus,~$|{\bf f}|(U)$ est de dimension~$\leq d$ en~$\xi$. Si~$(X,x)$ est sans bord,~$|{\bf f}|(U)$ est purement de dimension~$d$ en~$\xi$, et si de surcroît~$d=n$ alors~$|{\bf f}|(U)$ est un voisinage de~$\xi$ dans~$\rst n$. 

\medskip
\item[2)] Supposons que~$X$ est affinoïde et que~$x\in \partial X$. Il existe un voisinage~$k$-analytique compact et~$\Gamma$-strict~$Y$ de~$x$ dans~$X$ possédant la propriété suivante : pour tout voisinage analytique compact~$Z$ de~$x$ dans~$Y$, les germes de polytopes~$(|{\bf f}|(Y\cap \partial X),\xi)$ et~$(|{\bf f}|(Z\cap \partial X),\xi)$ coïncident. 
\end{itemize}}

\medskip
{\em Démonstration.} L'assertion 2) sera, en vertu du lemme \ref{lemmebord}, une conséquence de 1) ; il suffit donc de démontrer cette dernière. 

%
%

Soit~$\sch U$ l'ouvert quasi-compact de~$\PP_{\red{\hres(x)}/\red k}$ égal à l'image de~$\red{(X,x)}$. Le théorème \ref{theoconetem} assure l'existence d'un sous-foncteur définissable~$\sch D$ de~$\GG^n$, fermé et de dimension~$\leq d$, tel que pour tout~$H\in \mathsf D$, et tout~$(h_1,\ldots, h_n)\in H^n$, les assertions suivantes soient équivalentes : 

\medskip
i) il existe une valuation graduée~$\val . \in \sch U$, à valeurs dans un 
groupe~$H'\in \mathsf D_H$, telle que~$\val {\red{f_i(x)}}=h_i$ pour tout~$i$ ; 

ii)~$(h_1,\ldots, h_n)\in \sch D(H)$. 

\medskip
Choisissons une description de~$\sch D$ par une condition de la forme~$$\bigvee_i\bigwedge_j \phi_{i,j}\leq 1,$$ où chaque~$\phi_{i,j}$ est de la forme~$t_1^{e_1}\ldots t_n^{e_n}$ avec les~$e_i\in \ZZ$. 

\medskip
Soit~$P$ un~$(\QQ,\RR_+^*)$-polytope de~$\rst n$ contenant~$\xi$ et dont le germe en~$\xi$ est égal à celui du cône rationnel translaté~$\xi\cdot \sch D(\RR^*_+)$. Au voisinage de~$\xi$, le polytope~$P$ est décrit par la condition~$$\bigvee_i\bigwedge_j \phi_{i,j}\leq \phi_{i,j}(\xi).$$ Sa dimension en~$\xi$ est majorée par~$d$. 

\medskip
{\em L'image réciproque~$X_0:=|{\bf f}|\inv(P)$ est un voisinage de~$x$.} En effet,~$X_0$  est un domaine analytique de~$X$ ; en vertu de la théorie de Temkin,~$\red{(X_0,x)}$ est l'image réciproque sur~$\red{(X,x)}$ de l'ouvert quasi-compact~$$\bigcup_i \PP_{\red{\hres(x)}/\red k}\left\{\phi_{i,j}\left(\red{f_1(x)},\ldots, \red{f_n(x)}\right)\right\}_j.$$ Mais par définition de~$\sch D$, ce dernier ouvert quasi-compact contient~$\sch U$ ; par conséquent,~$\red{(X_0,x)}=\red{(X,x)}$ et~$X_0$ est bien un voisinage de~$x$ dans~$X$. 

\medskip 
Fixons un voisinage~$k$-analytique~$\Gamma$-strict~$U$ de~$x$ dans~$X_0$, et soit~$V$ un voisinage analytique compact de~$x$ dans~$U$. Nous allons démontrer que~$(|{\bf f}|(V),\xi)=(P,\xi)$, ce qui permettra de conclure puisqu'en faisant~$V=U$ il viendra alors~$(|{\bf f}|(U),\xi)=(P,\xi)$. Posons~$Q=|{\bf f}|(V)$ ; il s'agit de démontrer que le~$(\QQ,\RR^*_+)$-polytope~$Q$ est un voisinage de~$\xi$ dans~$P$. 

\medskip
Décrivons~$Q$ au voisinage de~$\xi$ par une condition de la forme~$$\bigvee_i\bigwedge_j \psi_{i,j}\leq \psi_{i,j}(\xi),$$  où chaque~$\psi_{i,j}$ est de la forme~$t_1^{e_1}\ldots t_n^{e_n}$ avec les~$e_i\in \ZZ$. 

\medskip
La réduction à la Temkin du germe de~$|{\bf f}|\inv(Q)$ en~$x$ est égale à l'image réciproque sur~$\red{(X,x)}$ de l'ouvert quasi-compact~$$\sch V:=\bigcup_i \PP_{\red{\hres(x)}/\red k}\left\{\psi_{i,j}\left(\red{f_1(x)},\ldots, \red{f_n(x)}\right)\right\}_j.$$ Comme~$|{\bf f}|\inv(Q)$ contient~$V$, lequel est un voisinage de~$x$, cette image réciproque est égale à~$\red{(X,x)}$ ; il s'ensuit que~$\sch V$ contient~$\sch U$. 

\medskip
Soit~$\eta$ un point de~$P$ différent de~$\xi$ ; écrivons~$\eta\xi\inv=(r_1,\ldots,r_n)$. Si~$\eta$ est suffisamment proche de~$\xi$ il appartient à~$\xi\cdot\sch D(\RR^*_+)$, et~$(r_1,\ldots, r_n)$ appartient donc à~$\sch D(\RR^*_+)$ ; il existe dès lors une valuation graduée~$\val .$ sur~$\red{\hres(x)}$, à valeurs dans un groupe ordonné contenant~$\RR^*_+$, appartenant à~$\sch U$ et qui envoie~$\red{f_i(x)}$ sur~$r_i$ pour tout~$i$. 

\medskip
L'ouvert~$\sch U$ étant contenu dans~$\sch V$, la valuation graduée~$\val .~$ appartient à~$\sch V$ ; cela signifie que~$$\bigvee_i\bigwedge_j \psi_{i,j}(r_1,\ldots, r_n)\leq 1,\;{\rm soit}\;{\rm encore}\;\bigvee_i\bigwedge_j \psi_{i,j}(\eta)\leq \psi_{i,j}(\xi).$$ Mais puisque~$$\bigvee_i\bigwedge_j \psi_{i,j}\leq \phi_{i,j}(\xi)$$ décrit~$Q$ au voisinage de~$\xi$, il en résulte que~$[\xi;\eta[$ est contenu dans~$Q$ au voisinage de~$\xi$. Cela ayant été établi pour tout~$\eta\in P-\{\xi\}$ suffisamment proche de~$\xi$, le polytope~$Q$ est un voisinage de~$\xi$ dans~$P$, ce qu'on souhaitait établir.

\medskip
Supposons que~$(X,x)$ soit sans bord. On a alors~$\red{(X,x)}=\PP_{\red{\hres(x)}/\red k}$. Il découle dès lors de l'assertion 2) du théorème \ref{theoconetem} : que~$\sch D$ est purement de dimension~$d$, ce qui implique que~$P$ est purement de dimension~$d$ en~$\xi$ ; et que si de plus~$d=n$ alors~$\sch D=\GG^n$, ce qui implique que~$P$ est un voisinage de~$\xi$ dans~$\rst n$.~$\Box$

\section{Polytopes analytiques et squelettes} 

On pose~$c=(\QQ,\sqrt{|k^*|\cdot \Gamma})$. 

\deux{introplcan} Soit~$X$ un espace~$k$-analytique~$\Gamma$-strict topologiquement séparé. 

\trois{polytanal} Soit~$P$ un compact de~$X$. Nous dirons qu'une structure~$c$-polytopale~$\Lambda_{c}(P)$ sur~$P$ est {\em analytique} si les  deux conditions suivantes sont satisfaites :

\medskip
a)~$\Lambda_{c}(P)$ admet une présentation~$P\hookrightarrow (\RR^*_+)^n$ de la forme~$(|f_1|, \ldots, |f_n|)_{|P}$, où les~$f_i$ sont des fonctions inversibles sur un domaine analytique~$\Gamma$-strict~$Y$ de~$X$ contenant~$P$ ; 

b) pour tout domaine analytique~$\Gamma$-strict~$Y$ de~$X$, l'intersection~$Y\cap P$ est un sous-espace~$c$-linéaire par morceaux de~$P$, et pour toute fonction analytique inversible~$f$ sur~$Y$, la fonction~$|f|_{|Y\cap P}$ est~$c$-linéaire par morceaux.

\trois{proprpolytanal} Supposons que~$P$ possède une structure polytopale analytique~$\Lambda_{c}(P)$, admettant une présentation~$(|f_1|, \ldots, |f_n|)_{|P}: P\hookrightarrow (\RR^*_+)^n$ comme au a) ci-dessus. Soit~$(g_1,\ldots, g_m)$ une famille quelconque de fonctions analytiques inversibles sur un domaine analytique de~$X$ contenant~$P$. En vertu de b), chacune des~$|g_j|_{|P}$ est~$c$-linéaire par morceaux ; par conséquent,~$(|f_1|, \ldots, |f_n|,|g_1|, \ldots, |g_m|)_{|P} : P\hookrightarrow  (\RR^*_+)^{n+m}$, est une présentation de~$\Lambda _{c}(P)$. 

\medskip
Faisons maintenant l'hypothèse que~$(|g_1|, \ldots, |g_m|)_{|P}$ constitue une présentation d'une structure polytopale analytique~$\Lambda' _{c}(P)$ sur~$P$. En appliquant ce qui précède à la structure~$\Lambda' _{c}(P)$, on voit que~$(|f_1|, \ldots, |f_n|,|g_1|, \ldots, |g_m|)_{|P} : P\hookrightarrow  (\RR^*_+)^{n+m}$ en constitue une présentation ; comme c'est également une présentation de~$\Lambda _{c}(P)$, il vient~$\Lambda' _{c}(P)=\Lambda _{c}(P)$.

\deux{defpolytanal} Soit~$X$ un espace~$k$-analytique~$\Gamma$-strict et topologiquement séparé et soit~$P$ un compact de~$X$. En vertu de \ref{proprpolytanal}, le compact~$P$  possède {\em au plus une} structure~$c$-polytopale analytique. Nous dirons que~$P$ est un {\em~$c$-polytope analytique} si les deux conditions suivantes sont satisfaites :

\medskip
$\bullet$ pour tout~$x\in P$, on a~$d_k(x)=\dim{} X$ ; 

$\bullet$ le compact~$P$ possède une (et partant une seule) une structure ~$c$-polytopale analytique.  

\medskip
Lorsque nous considérerons un ~$c$-polytope analytique, il sera toujours implicitement considéré comme étant muni de son unique structure ~$c$-polytopale analytique. 

\deux{polpolytanal} Soit~$X$ un espace~$k$-analytique~$\Gamma$-strict toplogiquement séparé. Les faits qui suivent résultent immédiatement des définitions. 

\trois{souspol} Soit~$Q$ une partie compacte d'un ~$c$-polytope analytique~$P$ de~$X$. Le compact~$Q$ est un ~$c$-polytope analytique de~$X$ si et seulement si c'est un ~$c$-polytope de~$P$, et si c'est le cas sa structure ~$c$-polytopale analytique coïncide avec sa structure ~$c$-polytopale héritée de~$P$. 

\trois{poldomgrand} Soit~$Y$ un domaine analytique~$\Gamma$-strict de~$X$. Si~$P$ est une partie compacte de~$Y$ alors~$P$ est un ~$c$-polytope analytique de~$Y$ si et seulement si c'est un ~$c$-polytope analytique de~$X$. 

\trois{insensred} Si~$P$ est une partie compacte de~$X$ alors~$P$ est un ~$c$-polytope analytique de~$X$ si et seulement si c'est un ~$c$-polytope analytique de~$X_{\rm red}$. 

\deux{testpolytanal} {\bf Lemme.} {\em Soit~$X$ un espace~$k$-affinoïde~$\Gamma$-strict intègre et soit~$P$ un compact de~$X$ tel que~$d_k(x)=\dim{}X$ pour tout~$x\in P$. 

\medskip
\begin{itemize}
\item[1)] Pour toute fonction non nulle~$f$ sur~$X$, la fonction~$|f|$ ne s'annule pas sur~$P$. 

\item[2)] Pour que~$P$ soit un ~$c$-polytope analytique, il suffit que les conditions suivantes soient satisfaites : 

\medskip
\begin{itemize}
\item[i)] pour toute famille~$(f_1,\ldots,f_n)$ de fonctions non nulles sur~$X$, l'image de~$P$ par~$(|f_1|,\ldots, |f_n|)$ est un ~$c$-polytope de~$(\RR^*_+)^n$ ; 

\item[ii)] il existe une famille~$(f_1,\ldots,f_n)$ de fonctions non nulles sur~$X$ telle que~$(|f_1|,\ldots, |f_n|)_{|P}$ soit injective. 
\end{itemize}
\end{itemize}
}

\medskip
{\em Démonstration.} Prouvons tout d'abord 1). Soit~$f$ une fonction non nulle sur~$X$. Comme~$d_k(x)=\dim{}X$, tout fermé de Zariski de~$X$ contenant~$x$ est de dimension égale à~$\dim {}X$, et coïncide donc avec~$X$. Il s'ensuit que le lieu des zéros de~$f$ ne contient pas~$x$, ce qui montre 1). 

\medskip
Montrons 2). On suppose que les conditions i) et ii) sont satisfaites, et l'on va prouver que~$P$ est un~$c$-polytope analytique. On choisit~$(f_1,\ldots, f_n)$ comme dans ii), et l'on munit~$P$ de la structure polytopale~$\Lambda_{c}(P)$ induite par ~$(|f_1|,\ldots, |f_n|)_{|P}$ ; on note~$Q$ l'image de~$P$ par~$(|f_1|,\ldots, |f_n|)$. Nous allons montrer que~$\Lambda_{c}(P)$ est analytique. Il s'agit de vérifier les conditions a) et b) de \ref{polytanal}.

\trois{condapolyt} {\em La condition a)}. Elle est vérifiée par construction. 

\trois{condbpolytprel} {\em La condition b) : préliminaires}.  Soit~$f$ une fonction non nulle sur~$X$. Nous allons montrer que~$|f|_{|P}\in \Lambda_{c}(P)$. L'assertion i) assure que l'image de~$P$ par~$(|f_1|,\ldots, |f_n|,|f|)$ est un~$c$-polytope~$Q'$ de~$(\RR^*_+)^{n+1}$ ; par choix des~$f_i$, la projection de~$Q$ sur les~$n$ premiers facteurs induit un homéomorphisme~$Q'\to Q$. Comme~$Q'$ et~$Q$ sont des polytopes, il existe une famille finie~$(Q_i)$ de polytopes de~$Q$, recouvrant~$Q$, tels que~$Q'\to Q$ admette une section affine sur chacun des~$Q_i$. 

\medskip
Pour tout entier~$i$, désignons par~$P_i$ l'image réciproque de~$Q_i$ sur~$P$. Les~$P_i$ sont des~$c$-polytopes de~$P$ qui recouvrent ce dernier. Par ailleurs, comme~$Q'\to Q$ admet une section au-dessus de~$Q_i$ pour tout~$i$, la restriction de~$|f|$ à~$P_i$ appartient à ~$\Lambda_{c}(P_i)$ pour tout~$i$. Il s'ensuit que~$|f|_{|P}\in  \Lambda_{c}(P)$. 

\trois{condbpolyt} {\em Vérification de la condition b).} Soit~$Y$ un domaine analytique~$\Gamma$-strict de~$X$ et soit~$f$ une fonction analytique inversible sur~$Y$. Soit~$x\in Y\cap P$. Il résulte de la description locale des domaines analytiques qu'il existe un voisinage affinoïde~$\Gamma$-strict~$V$ de~$x$ dans~$X$ et une famille finie ~$(W_i)$ de domaines~$\Gamma$-rationnels de~$V$, contenant~$x$ et tels que~$Y\cap V=\bigcup W_i$. Quitte à restreindre~$V$, on peut supposer que c'est lui-même un domaine~$\Gamma$-rationnel de~$X$. Posons~$\Pi_i=P\cap W_i$. Comme~$\bigcup \Pi_i$ est un voisinage de~$x$ dans~$Y\cap P$, et comme~$x\in \bigcap \Pi_i$, il suffit de démontrer que chacun des~$\Pi_i$ est un ~$c$-polytope de~$P$, et que~$|f|_{|\Pi_i}$ est~$c$-linéaire par morceaux pour tout~$i$. 

\medskip
Fixons~$i$, et écrivons~$W$ et~$\Pi$ au lieu de~$W_i$ et~$\Pi_i$. Le domaine~$W$ est lui-même un domaine rationnel de~$X$, et est par ailleurs~$\Gamma$-strict. Il est donc~$\Gamma$-rationnel : il peut être défini par une condition de la forme~$$|g_1|\leq \lambda_1|h|\;{\rm et}\;\ldots\;{\rm et}\;|g_m|\leq \lambda_m |h|,$$ où les~$g_j$ et~$h$ sont des fonctions analytiques sans zéro commun sur~$X$ et où les~$\lambda_i$ appartiennent à~$\Gamma$. Tout zéro de~$h$ sur~$W$ serait un zéro commun aux~$g_j$ et à~$h$, ce qui est absurde ; par conséquent~$h$ est inversible sur~$W$. Notons que~$W$ est non vide (il contient~$x$) ; par conséquent,~$h$ n'est pas la fonction nulle. 

\medskip
Comme une inégalité de la forme~$|0|\leq \lambda |h|$ est satisfaite sur tout~$X$, on peut supposer que les~$g_j$ sont tous non nuls. Le compact~$:Pi$ est le sous-ensemble de~$P$ défini par la condition~$$|g_1|\leq \lambda_1|h|\;{\rm et}\;\ldots\;{\rm et}\;|g_m|\leq \lambda_m |h|,$$ et les~$|g_j|_{|P}$ et~$|h|_{|P}$ sont~$c$-linéaires par morceaux d'après \ref{condbpolytprel}. En conséquence,~$\Pi$ est un ~$c$-polytope de~$P$. 

\medskip
Les fonctions de la forme~$g/h$, où~$g$ est une fonction analytique sur~$X$, sont denses dans l'ensemble des fonctions analytiques sur~$W$. La fonction~$f_{|W}$ étant inversible, il existe une fonction analytique~$g$ sur~$X$ telle que~$|f|=|g/h|$ en tout point de~$W$. La fonction~$g$ ne s'annule alors pas sur le domaine affinoïde non vide~$W$, et elle est en particulier non nulle sur~$X$. Les restrictions à~$\Pi$ de~$|f|$ et~$|g|/|h|$ coïncident, et~$|g|/|h|$ appartient à~$\Lambda_{c}(P)$ en vertu de  \ref{condbpolytprel}. Par conséquent,~$|f|_{|\Pi}\in \Lambda_{c}(\Pi)$, ce qu'il fallait démontrer.~$\Box$

\deux{exempolytanal} {\bf Exemple.} Soit~$n\in \NN$ et soit~$R$ un réel strictement  supérieur à~$1$ appartenant à~$\sqrt{|k^*|\cdot \Gamma}$. Soit~$P_R$ le sous-ensemble~$\{\eta_{\bf r}\}_{{\bf r}\in [1/R;R]^N}$ de~$\gman$. C'est un compact, contenu dans le domaine affinoïde~$\Gamma$-strict~$X_R$ de~$\gman$ défini par les conditions~$1/R\leq|T_i|\leq R$ pour tout~$i$ ; et pour tout~$x\in P_R$, on a~$d_k(x)=n$. Il résulte immédiatement de la description explicite des semi-normes~$\eta_{\bf r}$ que le compact~$P_R$ de l'espace affinoïde~$\Gamma$-strict~$X_R$ satisfait les conditions suffisantes du lemme \ref{testpolytanal} ci-dessus (en ce qui concerne la condition ii), on peut prendre~$f_i=T_i$ pour tout~$i$). C'est donc un~$c$-polytope analytique de~$X_R$, et également de~$\gman$. 

\deux{defskel} Soit~$X$ un espace~$k$-analytique topologiquement séparé et~$\Gamma$-strict. Nous dirons qu'une partie localement fermée~$\Sigma$ de~$X$ est un {\em~$c$-squelette} de~$X$ si l'ensemble des~$c$-polytopes analytiques de $X$ contenus dans~$\Sigma$ est un atlas~$c$-polytopal. Tout~$c$-squelette hérite par définition d'une structure d'espace~$c$-linéaire par morceaux. 

\trois{polegalsquel} Tout sous-espace~$c$-linéaire par morceaux d'un~$c$-polytope analytique~$P$ de~$X$ est un squelette de~$X$ : cela résulte de \ref{souspol}.

\trois{soussquel} Soit~$\Sigma$ un~$c$-squelette de~$X$ et soit~$\got A$ l'ensemble des~$c$-polytopes analytiques de~$X$ contenus dans~$\Sigma$. Soit~$\Sigma'$ un sous-ensemble localement fermé de~$\Sigma$. 

Supposons que~$\Sigma'$ soit un sous-espace~$c$-linéaire par morceaux de~$\Sigma$. La famille~$(P\cap \Sigma')_{P\in \got A}$ est un G-recouvrement de~$\Sigma'$ ; et pour tout~$P\in\got A$, l'intersection~$P\cap \Sigma'$ est un sous-espace~$c$-linéaire par morceaux de~$P$, et est donc G-recouverte par des éléments de~$\got A$, d'après \ref{polegalsquel}. Il s'ensuit que~$\Sigma'$ est un~$c$-squelette de~$X$. 

Réciproquement, supposons que~$\Sigma'$ soit un~$c$-squelette de~$X$. Il est alors G-recouvert par des éléments de~$\got A$, et est donc un sous-espace~$c$-linéaire par morceaux de~$\Sigma$. 

\trois{domsquel} Soit~$Y$ un domaine analytique~$\Gamma$-strict de~$X$. Une partie localement fermée de~$Y$ est un~$c$-squelette de~$Y$ si et seulement si c'est un~$c$-squelette de~$X$ : c'est une conséquence immédiate des définitions et de \ref{poldomgrand}. 

\trois{fonctsquel} Soit~$\Sigma$ un~$c$-squelette de $X$, soit~$Y$ un domaine analytique~$\Gamma$-strict de~$X$ et soit~$f$ une fonction inversible sur~$Y$. L'intersection~$Y\cap \Sigma$ est un~$c$-squelette de~$X$, ou encore, ce qui revient au même ({\em cf}. \ref{soussquel}) un sous-espace~$c$-linéaire par morceaux de~$\Sigma$ ; et la fonction~$|f|_{|\Sigma\cap Y}$ est~$c$-linéaire par morceaux. 

En effet, soit~$y\in Y\cap \Sigma$. Le point~$y$ possède un voisinage dans~$\Sigma$ qui est de la forme~$\bigcup P_i$ pour une famille finie~$(P_i)$ d'éléments de~$\got A$, tels que~$y\in \bigcap P_i$. Pour tout~$i$, l'intersection~$Y\cap P_i$ est un sous-espace~$c$-linéaire par morceaux de~$P_i$. Par conséquent,~$y$ possède dans~$Y\cap P_i$ un voisinage de la forme~$\bigcup P_{ij}$, ou~$(P_{ij})_j$ est une famille finie de~$c$-polytopes de~$P_i$ qui contiennent~$y$. 

La réunion~$\bigcup\limits_{i,j}P_{ij}$ est alors un voisinage de~$y$ dans~$Y\cap \Sigma$, et chacun des~$P_{ij}$ est un élément de~$\got A$ contenant~$y$. Il s'ensuit que~$Y\cap \Sigma$ est un sous-espace~$c$-linéaire par morceaux de~$Y$. 

Comme chaque~$P_{ij}$ appartient à~$\got A$, la fonction~$|f|_{|P_{ij}}$ est~$c$-linéaire par morceaux pour tout~$(ij)$. Par conséquent,~$|f|_{|Y\cap \Sigma}$ est~$c$-linéaire par morceaux au voisinage de~$y$, ce qui achève de prouver l'assertion requise. 

\deux{exsquel} Soit~$n\in \NN$. Nous allons montrer que~$S_n$ est un~$c$-squelette de~$\gman$ (rappelons que~$S_n$ est égal à~$\{\eta_{\bf r}\}_{{\bf r}\in \rst n}$). Soit~$\got A$ l'ensemble des~$c$-polytopes analytiques de~$\gman$ contenus dans~$S$. 

\medskip
Les~$c$-polytopes analytiques~$P_R$ définis au \ref{exempolytanal} pour~$R\in ]1;+\infty[\cap \sqrt{|k^*|\cdot \Gamma}$ appartiennent tous à~$\got A$ et sont contenus dans~$S$ ; de plus, leurs intérieurs (dans~$S$) recouvrent~$S$. Par conséquent,~$S$ est G-recouvert par les éléments de~$\got A$. 

\medskip
Par ailleurs, soient~$Q$ et~$Q'$ deux éléments de~$\got A$. Comme ils sont compacts, il sont tous deux contenus dans~$P_R$ pour un certain~$R\in  ]1;+\infty[\cap \sqrt{|k^*|\cdot \Gamma}$. Ce sont alors deux  ~$c$-polytopes analytiques de~$\gman$ contenus dans le ~$c$-polytope analytique~$P_R$ ; il s'ensuit que ce sont deux ~$c$--polytopes de~$P_R$, et leur intersection est dès lors un ~$c$-polytope de~$P_R$, et partant un ~$c$-polytope analytique ; c'est par conséquent un ~$c$-polytope de~$Q$ aussi bien que de~$Q'$, et les deux structures ~$c$-polytopales dont il hérite ainsi coïncident en vertu de \ref{souspol}.

\medskip
Ainsi,~$\got A$ est un atlas~$c$-polytopal sur~$S$,  et celui-ci est un~$c$-squelette de~$\gman$.

\deux{testsquel} {\bf Proposition.} {\em Soit~$X$ un espace~$k$-analytique~$\Gamma$-strict topologiquement séparé et soit~$(X_i)$ un G-recouvrement de~$X$ par des domaines analytiques~$\Gamma$-stricts. Soit~$\Sigma$ une partie localement fermée de~$X$. Les propositions suivantes sont équivalentes : 

\medskip
i)~$\Sigma$ est un~$c$-squelette de~$X$ ; 

ii) pour tout~$i$, l'intersection~$\Sigma\cap X_i$ est un~$c$-squelette de~$X_i$. }

\medskip
{\em Démonstration.} On sait déjà que i)$\Rightarrow$ii) (\ref{fonctsquel}). Réciproquement, supposons que ii) soit vraie. Soit~$\got A$ l'ensemble des~$c$-polytopes analytiques de~$X$ contenus dans~$\Sigma$. Soit~$x\in \Sigma$. Par définition d'un G-recouvrement, il existe un ensemble fini~$I$ d'indices tel que~$x\in \bigcap\limits_{i\in I}X_i$ et tel que ~$\bigcup\limits_{i\in I}X_i$ soit un voisinage de~$x$ dans~$X$. Fixons~$i$. Comme~$\Sigma\cap X_i$ est un squelette de~$X_i$, il existe une famille finie~$(P_{ij})_j$ de~$c$-polytopes analytiques de~$X_i$ (ce sont aussi des~$c$-polytopes analytiques de~$X$) qui contiennent~$x$, et tels que~$\bigcup P_{ij}$ soit un voisinage de~$x$ dans~$\Sigma\cap X_i$. 

\medskip
La réunion des~$P_{ij}$ pour~$i\in I$ et~$j$ variable est alors un voisinage de~$x$ dans~$\Sigma$. Comme chacun des~$P_{ij}$ est un élément de~$\got A$ qui contient~$x$, on voit que~$\got A$ constitue un G-recouvrement de~$\Sigma$. 

\medskip
Soient maintenant~$P$ et~$Q$ deux éléments de~$\got A$. Nous allons montrer que~$P\cap Q$ est un~$c$-polytope analytique, ce qui permettra de conclure. 

Comme~$P$ est un~$c$-polytope analytique,~$P\cap X_i$ est pour tout~$i$ un sous-espace~$c$-linéaire par morceaux de~$P$, et est donc G-recouvert par les~$c$-polytopes analytiques qu'il contient. Ainsi,~$P$ est G-recouvert par ses~$c$-polytopes qui sont contenus dans l'un des~$X_i$. Par compacité, il s'écrit comme une union finie~$\bigcup P_j$, où chaque~$P_j$ est un polytope analytique de~$X_{i(j)}$ pour un certain indice~$i(j)$. 

On écrit de même~$Q$ comme une union finie~$\bigcup Q_\ell$, où chaque~$Q_\ell$ est un~$c$-polytope analytique de~$X_{i(\ell)}$ pour un certain indice~$i(\ell)$. 

\medskip
Les~$c$-polytopes analytiques de~$X$ contenus dans~$P$ sont exactement les~$c$-polytopes de~$P$, et une union finie de~$c$-polytopes analytiques de~$X$ contenue dans~$P$ est donc encore un~$c$-polytope analytique de~$X$. Pour montrer que~$P\cap Q$ est un~$c$-polytope analytique, il suffit dès lors de montrer que~$P_j\cap Q_\ell$ est un~$c$-polytope analytique pour tout~$(j,\ell)$. On se ramène ainsi au cas où il existe deux indices~$i_1$ et~$i_2$ tel que~$P\subset X_{i_1}$ et~$Q\subset X_{i_2}$. 

\medskip
L'intersection de~$P$ avec~$X_{i_1}\cap X_{i_2}$ est un sous-espace~$c$-linéaire par morceaux de~$P$, et est donc un~$c$-squelette de~$X_{i_1}\cap X_{i_2}$ ; ce~$c$-squelette est contenu dans l'intersection~$\Sigma\cap X_{i_1}\cap X_{i_2}$, qui est lui-même un~$c$-squelette de~$X_{i_1}\cap X_{i_2}$ puisque~$\Sigma \cap X_{i_1}$ est un~$c$-squelette de~$X_{i_1}$. Par conséquent, ~$P\cap X_{i_1}\cap X_{i_2}$ est un sous-espace~$c$-linéaire par morceaux de~$\Sigma\cap X_{i_1}\cap X_{i_2}$. Il en va de même de~$Q\cap X_{i_1}\cap X_{i_2}$. Il en résulte que~$$P\cap Q=(P\cap (X_{i_1}\cap X_{i_2}))\cap (Q\cap (X_{i_1}\cap X_{i_2}))$$ est un sous-espace~$c$-linéaire par morceaux du squelette ~$\Sigma\cap X_{i_1}\cap X_{i_2}$. Étant de sucroît compact,~$P\cap Q$ est réunion finie de~$c$-polytopes analytiques ; comme il est inclus dans le~$c$-polytope analytique~$P$, il est lui-même un~$c$-polytope analytique, ce qui achève la preuve.~$\Box$ 

\section{Images réciproques du squelette standard de~$\gman$}

Cette section est entièrement consacrée à la démonstration du théorème suivant ; on rappelle que~$S_n$ désigne le sous-ensemble~$\{\eta_{\bf r}\}_{{\bf r}\in (\RR^*_+)^n}$ de~$\gman$. Si~$F$ est une extension complète de~$k$, on notera~$S_{n,F}$ le sous-ensemble~$\{\eta_{F,\bf r}\}_{{\bf r}\in (\RR^*_+)^n}$ de~${\mathbb G}_{m,F}^{n,\rm an}$. C'est un~$(\QQ,\sqrt {|F^*|\cdot \Gamma})$-squelette de~${\mathbb G}_{m,F}^{n,\rm an}$ (\ref{exsquel}). 

\deux{imrecsquels} {\bf Théorème.} {\em Soit~$n$ un entier et soit~$X$ un espace~$k$-analytique topologiquement séparé et de dimension~$\leq n$. Soit~$(\phi_1,\ldots, \phi_m)$ une famille finie de morphismes de~$X$ vers~$\gman$. Pour toute extension complète~$F$ de~$k$, on désigne par~$c_F$ le couple~$(\QQ,\sqrt{|F^*|\cdot \Gamma})$ ; on écrira~$c$ au lieu de~$c_k$. 

\medskip
1) La réunion des~$\phi_j^{-1}(S_n)$ est un~$c$-squelette de~$X$, vide si~$\dim {} X<n$. 

2) Pour tout~$j$, l'application~$\phi_j^{-1}(S_n)\to S_n$ est G-localement une immersion d'espaces~$c$-linéaires par morceaux (notons que~$\phi_j^{-1}(S_n)$ est un~$c$-squelette de~$X$, d'après l'assertion 1) appliquée au cas d'un seul morphisme). 

3) Pour toute extension complète~$F$ de~$k$, l'application~$$\bigcup \phi_{j,F}^{-1}(S_{n,F}) \to \bigcup \phi_{j,F}^{-1}(S_n)$$ est surjective, et c'est G-localement une immersion d'espaces~$c_F$-linéaires par morceaux. 

4) Supposons de plus que~$X$ est compact. Il existe alors une extension finie séparable~$F_0$ de~$k$ telle que pour toute extension complète~$F$ de~$F_0$ l'immersion d'espaces~$c_F$-linéaires par morceaux~$$\bigcup \phi_{j,F}^{-1}(S
_{n,F}) \to \bigcup \phi_{j,F_0}^{-1}(S_{n,F_0})$$ soit un isomorphisme.}

\deux{dkxegaln} Commençons par quelques remarques. Soit~$x\in S_n$, soit~$j\in\{1,\ldots,m\}$ et soit~$y\in \phi_j\inv(x)$. Comme~$d_k(x)=n$, on a~$d_k(y)\geq n$ pour tout antécédent~$y$ de~$x$ par~$\phi_j$. La dimension de~$X$ étant majorée par~$n$, il vient~$d_k(y)=n$, et partant~$d_{\hres(x)}(y)=0$. En conséquence : 

\medskip
$\bullet$ si~$\dim{}X<n$ alors~$\phi_j\inv(S_n)=\emptyset$ ; 

$\bullet$ les fibres de~$\phi_j$ en les points de~$S_n$ sont purement de dimension nulle. 

\deux{planpreuve} Le plan de la preuve est le suivant. Nous allons tout d'abord établir le théorème dans le cas où~$m=1$, c'est-à-dire dans le cas où l'on ne considère qu'un seul morphisme de~$X$ vers~$\gman$. Nous montrerons ensuite l'assertion 1) dans le cas général, ce qui suffira à conclure : la combinaison de 1) dans le ças général et de 3) et 4) dans le cas où~$m=1$ entraîne aussitôt 3) et 4) dans le cas général. 

\subsection*{Preuve dans le cas où~$m=1$}

Pour alléger les notations, nous écrirons~$\phi$ au lieu de~$\phi_1$. Nous allons commencer par quelques réductions. 

\deux{corpsetred} Les énoncés requis peuvent se tester après extension des scalaires au complété~$\widehat{k_{\rm parf}}$ de la clôture parfaite~$k_{\rm parf}$ de~$k$. C'est immédiat, sauf peut-être à propos de 4) : concernant cette dernière assertion, il faut utiliser le fait que~$F\mapsto F\otimes_k\widehat{k_{\rm parf}}$ établit une équivalence entre la catégorie des extensions finies séparables de~$k$ et celle des extensions finies séparables de~$k_{\rm parf}$. Cela résulte de l'invariance topologique du site étale, appliquée à~$k\hookrightarrow k_{\rm parf}$, et du lemme de Krasner, pour passer à la complétion.  

\medskip   
Ils sont par ailleurs insensibles à la présence de nilpotents : on peut donc supposer~$k$ parfait et~$X$ réduit. 

\deux{bloc} En vertu de \ref{fonctsquel} et de la proposition \ref{testsquel}, on peut raisonner G-localement sur~$X$. Cela permet tout d'abord de se ramener au cas où~$X$ est affinoïde (et toujours~$\Gamma$-strict). Puis il suffit, par compacité, d'exhiber pour tout~$x\in X$ un voisinage affinoïde~$\Gamma$-strict~$V$ de~$x$ dans~$X$ tel que~$(V,\phi_{|V})$ satisfasse le théorème. 

\medskip
Soit~$x\in X$. Si~$x$ n'appartient pas au compact~$\phi\inv (S_n)$, il existe un voisinage affinoïde~$\Gamma$-strict~$V$ de~$x$ dans~$X$ qui ne rencontre pas~$\phi\inv (S_n)$ ; le théorème est trivialement vérifié par~$(V,\phi_{|V})$. 

\medskip
Supposons maintenant que~$x\in \phi\inv(S_n)$. On a alors~$d_k(x)=n$ d'après \ref{dkxegaln} ; il s'ensuit que~$\dim x X=n$, et que l'anneau local~$\sch O_{X,x}$ est artinien (\cite{flatn}, cor.~1.12). L'espace~$X$ étant réduit, $\sch O_{X,x}$ est un corps. Il est en particulier régulier, ce qui implique,~$k$ étant parfait, que~$X$ est quasi-lisse en~$x$. 

\medskip
 Le lemme \ref{algquasil} assure alors qu'il existe une~$k$-variété affine et lisse~$\sch X$, de dimension~$n$, et un voisinage affinoïde~$\Gamma$-strict~$V$ de~$x$ dans~$X$ qui s'identifie à un domaine affinoïde de~$\sch X\an$. On peut toujours supposer~$V$ connexe, puis~$\sch X$ connexe, et donc intègre. 
 
 \medskip
 Comme~$d_k(x)=n$, le point~$x$ est situé au-dessus du point générique~$\eta$ de~$\sch X$. Le morphisme~$\phi$ est donné par~$n$ fonctions inversibles~$f_1,\ldots, f_n$ sur~$V$. Le corps~$\kappa(\eta)=k(\sch X)$ est dense dans~$\hres(x)$. Il existe donc des fonctions rationnelles~$g_1,\ldots,g_n$ sur~$\sch X$ telles que~$|g_i(x)-f_i(x)|<|f_i(x)|$ pour tout~$i$. On peut restreindre~$V$ de sorte que les~$g_i$ soient définies sur~$V$, et de sorte que l'on ait~$|f_i-g_i|<|f_i|$ en tout point de~$V$. Soit~$\psi : V\to \gman$ le morphisme induit par les~$g_i$. 
 
 \medskip
Soit~$F$ une extension complète du corps~$k$ et soit~$y$ un point de~$V_F$. On a~$|f_i(y)-g_i(y)|<|f_i(y)|$ ; autrement dit,~$\widetilde{f_i(y)}=\widetilde{g_i(y)}$ dans le corps gradué résiduel~$\widetilde{\hres(y)}$. La famille~$(\widetilde{f_i(y)})_i$ est donc algébriquement indépendante sur~$\red F$ si et seulement si c'est le cas de~$(\widetilde{g_i(y)})_i$. Cette équivalence peut, en vertu de \ref{transeteta}, se reformuler ainsi : le point~$y$ appartient à~$\phi\inv(S_{n,F})$ si et seulement si il appartient à~$\psi\inv(S_{n,F})$. On peut donc, en remplaçant~$f_i$ par la restriction de~$g_i$ à~$V$ pour tout~$i$, supposer que~$\phi$ est induit par un morphisme~$\sch U\to {\mathbb G}_{m,k}^n$, où~$ \sch U$ est un  ouvert de Zariski de~$\sch X$ tel que~$V\subset \sch U\an$ (par abus, on notera encore~$\phi$ ce morphisme, et~$(f_1,\ldots, f_n)$ la famille de fonctions inversible sur~$\sch U$ qui le définit). 

Comme~$\phi(x)\in S_n$, il est situé au-dessus du point générique de ~${\mathbb G}_{m,k}^n$. Par conséquent,~$\phi(\eta)$ est le point générique de~${\mathbb G}_{m,k}^n$, et le corps~$L:=k(\sch X)$ apparaît ainsi comme une extension finie de~$k(T_1,\ldots, T_n)$. 


\medskip
Le théorème \ref{propvalg} assure l'existence d'une famille finie~$(g_1,\ldots, g_r)$ d'éléments de~$L$ telle que les~$g_i$ séparent les prolongements de~$\eta_{\bf r}$ à~$L$ pour tout~${\bf r}\in (\RR^*_+)^n$ ; comme~$L$ est un corps, on peut supposer que les~$g_i$ sont inversibles. Le théorème~\ref{propvalg} assure également qu'il existe une extension finie séparable~$F_0$ de~$k$ et une famille finie~$p_1,\ldots, p_s$ d'éléments de~$F_0\otimes_kL$ telles que~$p_1,\ldots, p_s$ séparent les extensions de~$\eta_{F_0,{\bf r}}$ à~$F\otimes_kL$ pour toute extension complète~$F$ de~$F_0$.

 \medskip
 Il existe un voisinage affinoïde~$\Gamma$-strict~$W$ de~$x$ dans~$\sch U\an$ sur lequel les~$g_i$ sont définies et inversibles et tel que les~$p_j$ soient définies sur~$W_{F_0}$ (pour ce dernier point, on utilise le fait que~$\sch U\an_{k'}\to \sch U\an$ est fini, et en particulier topologiquement propre). 
 
 \medskip
 Le point~$x$ étant situé au-dessus du point générique de~$\sch X$, et~$\phi : \sch X\to {\mathbb G}^n_{m,k}$ étant génériquement finie,~$\phi$ est finie en~$x$. Par ailleurs, le point~$\phi(x)$ est situé sur~$S_n$, et~$\sch O_{\gman, \phi(x)}$ est dès lors un corps ; en conséquence,~$\phi$ est fini et plat en~$x$. 
 
 On peut dès lors restreindre~$W$ de sorte que~$\phi$ induise un morphisme fini et plat~$W\to W'$, où~$W'$ est un voisinage affinoïde~$\Gamma$-strict de~$\phi(x)$ dans~$\gman$. Un morphisme fini et plat est ouvert ; il en résulte que~$\phi(W)$ contient un voisinage de~$\phi(x)$ dans~$S_n$, et en particulier un pavé~$P$ de la forme~$\{\eta_{(r_1,\ldots, r_n}\}_{s_i<r_i<t_i}$, où les~$s_i$ et les~$r_i$ appartiennent à~$\sqrt{|k^*|\cdot \Gamma}$ pour tout~$i$, et où~$s_i<|T_i(x)|<t_i$ pour tout~$i$. Soit~$W'_0$ l'intersection de~$W$ avec le domaine affinoïde~$\Gamma$-strict de~$\gman$ défini par la condition~$|{\bf T}|\in P$ ; c'est un voisinage de~$\phi(x)$ dans~$W'$. Quitte à remplacer~$W$ par~$W\times_{W'}W'_0$,  on peut supposer que ~$\phi(W)\cap S_n$ est égal au pavé~$P$. 
 
 On peut remplacer~$V$ par~$V\cap W$, et donc supposer que~$V$ est un domaine affinoïde~$\Gamma$-strict de~$W$. Il résulte du \ref{fonctsquel} que pour que 1), 2), 3) et 4) soient vraies pour~$(V,\phi)$, il suffit qu'elles le soient pour~$(W,\phi)$.

 \trois{preuves}{ \em Preuve de 1, 2), et 3) pour~$(W,\phi)$.} Pour prouver 1), nous allons établir que~$\phi_{|W}\inv(S_n)$ est un~$c$-polytope analytique, en utilisant les conditions suffisantes exhibées par le lemme \ref{testpolytanal}. Commençons par une remarque : si~$(h_1,\ldots, h_m)$ est une famille de fonctions analytiques inversibles sur~$W$, et si l'image par~$|{\bf h}|$ de~$W$ est un~$c$-polytope de~$(\RR^*_+)^m$, son image par~$(|h_1|,\ldots, |h_{\ell}|)$ est un~$c$-polytope de~$(\RR^*_+)^\ell$ pour tout~$\ell\leq m$ (on utilise la projection de~$(\RR^*_+)^n$ sur les~$\ell$ premières coordonnées). 
 
En conséquence, pour que les conditions suffisantes du lemme \ref{testpolytanal} sont satisfaites, il suffit que pour toute famille~$(h_1,\ldots, h_\ell)$ de fonctions inversibles sur~$W$, l'application~$(|f_1|,\ldots,|f_n|, |g_1|,\ldots, |g_r|,|h_1|,\ldots |h_\ell|)$ induise un homéomorphisme entre~$Q:=\phi_{|W}\inv(S_n)$ et un~$c$-polytope de~$(\RR^*_+)^{n+r+\ell}$. 
 
 \medskip
 Soit donc~$(h_1,\ldots, h_\ell)$ une famille de fonctions analytiques inversibles sur~$W$, et soit~$\lambda$ l'application~$(|f_1|,\ldots,|f_n|, |g_1|,\ldots, |g_r|,|h_1|,\ldots |h_\ell|)_{|W}$. Le théorème \ref{tropicglob} assure que~$\lambda(W)$ est un~$c$-polytope de~$(\RR^*_+)^{n+r+\ell}$. Comme les fonctions~$g_i$ séparent les antécédents de~$\eta_{\bf r}$ sur le corps~$L$ pour tout~$\bf r$, la restriction de~$\lambda$ à~$Q$ est injective, et induit en conséquence, par un argument de compacité, un homéomorphisme entre~$Q$ et son image~$\lambda(Q)$ dans~$\lambda(W)$. Soit~$\pi$ la projection sur les~$n$ premières coordonnées, de~$(\RR_+^*)^{n+r+\ell}$ sur~$(\RR^*_+)^n$. 
 
 \medskip
 Il reste à s'assure que~$\lambda(Q)$ est un~$c$-polytope. Pour cela, on choisit une triangulation~$\sch T$ de~$\lambda(W)$ en~$c$-polytopes convexes (on ne demande pas qu'il s'agisse de simplexes) de sorte que la projection~$\pi$ sur les~$n$ premières coordonnées soit affine sur chacune des cellules fermées de la triangulation. On note~$Q'$ la réunion des cellules fermées et de dimension~$n$ de~$\sch T$ auxquelles la restriction de~$\pi$ est injective ; nous allons montrer que~$\lambda(Q)=Q'$, par double inclusion. 
 
 \medskip
 {\em On a~$\lambda(Q)\subset Q'$.} Soit~$\xi$ un point de~$\lambda(W)-Q'$. Le point~$\xi$ possède un voisinage ouvert~$\Omega$ dans~$\lambda(W)$ dont l'image par~$\pi$ est contenue dans un polytope compact de dimension~$n-1$. 
  
 L'image réciproque~$\lambda\inv(\Omega)$ est un ouvert de~$W$. Si~$\lambda\inv(\Omega)$ rencontrait~$\phi\inv(S_n)$, alors en vertu du caractère ouvert des morphismes finis et plats, son image par~$\phi$ contiendrait un~$c$-polytope analytique~$P'$ de~$\gman$, contenu dans~$S_n$ et de dimension~$n$ ; on aurait par construction~$P'\subset \pi(\Omega)$, en contradiction avec le choix de~$\Omega$. Ainsi,~$\xi\notin \lambda(Q)$. 
 
 \medskip
 {\em On a~$Q'\subset \lambda(Q)$.} Soit~$\xi$ un point de~$Q'$. Supposons que~$\lambda\inv(\xi)$ ne rencontre pas~$\phi\inv(S_n)$. Par propreté topologique de~$\lambda$, il existe un voisinage polytopal compact~$R$ de~$\xi$ dans~$\lambda(W)$ tel que~$\lambda\inv(R)$ ne rencontre pas~$\phi\inv(S_n)$. 
 
 Soit~$y\in \lambda\inv(R)$. Comme il n'appartient pas à~$\phi\inv(S_n)$, la famille des~$\widetilde{f_i(y)}$ est algébriquement liée sur le corps gradué~$\red k$ (\ref{transeteta}) ;   le théorème \ref{theotroploc} assure alors
 que~$y$ possède un voisinage analytique compact dans~$\lambda\inv(R)$ dont l'image par~$(|f_1|, \ldots, |f_n|)$ est un polytope de dimension~$\leq n-1$. 
 
 Par compacité, l'image de ~$\lambda\inv(R)$ par~$(|f_1|, \ldots, |f_n|)$ est un polytope de dimension~$\leq n-1$. Mais cette image est égale à~$\pi(R)$. Or~$R$ est un voisinage polytopal compact de~$\xi$, lequel appartient à~$Q'$ ; par conséquent,~$\pi(R)$ contient un polytope de dimension~$n$, et on aboutit ainsi à une contradiction. 
 
 Il s'ensuit que 1) est satisfaite par~$(W,\phi)$, et donc vraie en général dès que~$m=1$. 
 
 \medskip
 Montrons maintenant 2). Il résulte de 1) que l'application~$(|f_1|,\ldots, |f_n|)_{|Q}$ est~$c$-linéaire par morceaux. En conséquence,~$Q\to S_n$ est~$c$-linéaire par morceaux. Par ailleurs, on a vu au \ref{dkxegaln} que les fibres de~$\phi$ en les points de~$S_n$ sont purement de dimension nulle. Par compacité, les fibres de~$\phi_{|W}$ en les points de~$S_n$ sont finies ; ainsi,~$Q\to S$ est une application~$c$-linéaire par morceaux à fibres finies, et c'est donc une immersion G-locale d'espaces~$c$-linéaires par morceaux. Il en résulte que 2) est vraie pour~$(W,\phi)$, et donc en général lorsque~$m=1$. 
 
 \medskip
 Il reste à établir 3). Soit~$F$ une extension complète de~$k$. Considérons le diagramme commutatif~$$\diagram (\phi_{|W})_F\inv(S_{n,F}) \dto \rto &S_{n,F}\dto \\ \phi_{|W}\inv (S_n)\rto &S_n\enddiagram.$$ Les flèches horizontales sont, en vertu de l'assertion 2) déjà prouvée lorsque~$m$ est égal à~$1$, des immersions G-locales d'espaces~$c_F$-linéaires par morceaux, et la flèche verticale de droite est un isomorphisme d'espaces~$c_F$-linéaires par morceaux, d'après la description directe de~$S_n$, de~$S_{n,F}$ et des structures linéaires par morceaux sur ces derniers. Il s'ensuit que la flèche verticale de gauche est une immersion G-locale d'espaces~$c$-linéaires par morceaux. Elle est par ailleurs surjective. En effet, soit~$y\in \phi\inv(S_n)$ et soit~$z$ le point de~$S_{n,F}$ qui correspond à~$\phi(y)$ {\em via} l'homéomorphisme~$S_{n,F}\to S_n$. Tout idéal maximal de la~$\hres(z)$-algèbre finie et non nulle~$\hres(z)\otimes_{\hres(\phi(y)}\hres(y))$ définit alors un antécédent de~$z$ sur~$W_F$ situé au-dessus de~$y$. Ainsi, 3) est vraie pour~$(W,\phi)$, et partant en général lorsque~$m=1$. 
 
 \trois{preuve4wphi} {\em Preuve de 4) pour~$(W,\phi)$}. Soit~$F$ une extension complète de~$F_0$. On dispose d'un diagramme commutatif 
 $$\diagram \phi_F\inv(S_{n,F})\dto\rto&S_{n,F}\dto\\
 \phi_{F_0}\inv(S_{n,F_0})\rto&S_{n,F_0}\enddiagram$$
 dans lequel la flèche verticale de droite est un homéomorphisme. Par construction, les fonctions~$p_i$ sont définies sur~$W_{F_0}$, 
 et les fonctions~$\val{p_i}$ séparent les points des fibres
 de chacune des deux flèches horizontales du diagramme. Il s'ensuit immédiatement que la flèche verticale de gauche est injective. 
 Comme c'est par ailleurs, en vertu
 de l'assertion 3) déjà prouvée pour~$m=1$, une surjection,
 et une immersion G-locale d'espaces~$c_F$-linéaires par morceaux, c'est un isomorphisme d'espaces~$c_F$-linéaires par
 morceaux. Ceci achève de montrer~4) pour~$(W,\phi)$, et donc en général lorsque~$m=1$.

\deux{astucesection} {\em Preuve de 1) dans le cas général.} On peut là encore, en vertu du lemme \ref{testsquel}, raisonner G-localement sur~$X$ et donc supposer~$X$ compact. Pour tout~$j$, notons ~$f_{1,j},\ldots, f_{n,j}$ les fonctions inversibles qui définissent~$\phi_j$. Soit~$Y$ l'espace affine relatif sur~$X$ de dimension (relative) égale à~$nm$ et soit~$(T_{i,j})$ la famille des fonctions coordonnées sur~$Y$. Pour tout~$i$, posons~$g_i=\sum_j f_{i,j} T_{i,j}$. 

Pour tout~$x\in X$, désignons par~$\sigma(x)$ le point~$\eta_{r_{i,j}(x)}$ de~$Y_x\simeq \Aff^{nm,{\rm an}}_{\hres(x)}$, où~$r_{i,j}(x)$ est égal pour tout~$(i,j)$ à~$1/|f_{i,j}(x)|$. L'application~$\sigma$ est une section de~$Y\to X$. Elle est continue : pour le voir, on se ramène immédiatement au cas où~$X$ est affinoïde, d'algèbre associée~$\mathscr A$ ; et l'on remarque alors que pour tout polynôme~$P=\sum a_I {\bf T}^I \in \mathscr A[{\bf T}]$ l'application qui envoie le point~$x$ sur~$|P(\sigma(x))|=\max |a_I(x)| {\bf r}(x)^I$ est continue. 

\medskip
Soit~$\psi$ le morphisme~$Y\to \gm^{nm+n,\rm an}$ défini par les~$T_{i,j}$ et les~$g_i$, et soit~$\Theta$ le fermé~$\psi^{-1}(S_{m+n})$ de~$Y$. En vertu du cas~$m=1$ déjà traité,~$\Theta$ est un~$c$-squelette de~$Y$.

\trois{sigmaxintheta} Soit~$x\in \bigcup \phi_j\inv(S_n)$. Le point~$y:=\sigma(x)$ de~$Y$ appartient alors à~$\Theta$. En effet, soit~$j_0$ tel que~$x\in \phi_{j_0}\inv(S_n)$. Il s'agit de montrer que~$$\left(\;(\widetilde{T_{i,j}(y)})_{i,j},\;(\widetilde{g_i(y)})_i\right)$$ est une famille algébriquement indépendante sur le corps gradué~$\red k$ (\ref{transeteta}). Or il résulte de notre construction que :

\medskip
\begin{itemize} 
\item les~$\widetilde{T_{i,j}(y)}$ sont algébriquement indépendants sur~$\widetilde{\hres(x)}~$, et \emph{a fortiori} sur~$\red k( \widetilde{f_{i,j}(x)})_{i,j}$ ; 

\item les~$\widetilde{f_{i,j_0}(x)}$ sont algébriquement indépendants sur~$\red k$. 

\end{itemize} 
\medskip
L'assertion requise résulte dès lors du lemme~\ref{indepchiant} ci-dessous. 

\medskip
\trois{conclusquelm} Le cas~$m=1$ déjà traité assure que~$\phi_j\inv(S_n)$ est un~$c$-squelette de~$X$ pour tout~$j$ ; il s'écrit donc comme une union finie~$(P_{j,\ell})_{\ell}$ de~$c$-polytopes analytiques de~$X$. En conséquence, le compact~$\bigcup \phi_j\inv(S_n)$ est une union finie de~$c$-polytopes analytiques de~$X$. 

\medskip
Soient~$P$ et~$Q$ deux~$c$-polytopes analytiques de~$X$ contenus dans~$\bigcup \phi_j\inv(S_n)$. Nous allons démontrer que leur intersection est un $c$-polytope analytique de~$X$, ce qui achèvera de prouver que~$\bigcup\phi_j\inv(S_n)$ est un squelette. 

\medskip
Comme~$P$ est un~$c$-polytope analytique de~$X$, il existe un domaine analytique~$\Gamma$-strict~$Z$ de~$X$ contenant~$P$, et une famille~$(h_1,\ldots, h_r)$ de fonctions inversibles sur~$Z$ telles que~$(|h_1|,\ldots, |h_r|)$ identifie~$P$ à un~$c$-polytope~$\Pi$ de~$(\RR^*_+)^r$. L'image réciproque de~$Z$ sur~$Y$ est un domaine analytique~$\Gamma$-strict de~$Y$. Son intersection avec~$\Theta$ est donc un sous-espace~$c$-linéaire par morceaux de ce dernier, et~$|h_s|_{|(Y\times_XZ)\cap \Theta}$ est~$c$-linéaire par morceaux pour tout~$s$. En vertu de \ref{sigmaxintheta}, le compact~$\sigma(P)$ est contenu dans~$(Y\times_XZ)\cap \Theta$. Comme~$P\simeq \sigma(P)$, la famille~$(|h_1|,\ldots, |h_r|)_{|(Y\times_XZ)\cap \Theta}$ identifie~$\sigma(P)$ au~$c$-polytope~$\Pi$. 

\medskip
Les fonctions~$|h_s|_{|(Y\times_XZ)\cap \Theta}$ étant linéaires par morceaux, il en résulte que~$\sigma(P)$ est un~$c$-polytope du squelette~$\Theta$, et que sa structure~$c$-polytopale est induite par l'homéomorphisme~$$(|h_1|,\ldots, |h_r|)_{|\sigma(P)}: \sigma(P)\simeq \Pi.$$ Autrement dit,~$\sigma$ induit un isomorphisme~$c$-linéaire par morceaux entre~$P$ et~$\sigma(P).$ De même, le compact~$\sigma(Q)$ est un~$c$-polytope de~$\Theta$ (et~$\sigma$ induit un isomorphisme~$c$-linéaire par morceaux entre~$Q$ et~$\sigma(Q)$.)

Comme~$\sigma(P)$ et~$\sigma(Q)$ sont deux~$c$-polytopes de~$\Theta$, leur intersection est un~$c$-polytope de~$\Theta$, et est donc un~$c$-polytope de~$\sigma(P)$. Il s'ensuit que~$P\cap Q$ est un~$c$-polytope de~$P$ et donc un $c$-polytope analytique de $X$, ce qui achève la démonstration du théorème.~$\Box$

\deux{indepchiant} {\bf Lemme.} {\em Soit~$K$ un corps gradué et soit~$L$ une extension 
graduée 
de~$K$. Soient~$(f_{i,j})$ et~$(T_{i,j})$ 
deux familles finies d'éléments homogènes de~$L$ indexées
par les couples~$(i,j)$ d'entiers tels que~$1\leq i\leq n$ et~$1\leq 
j\leq m$.
Pour tout~$(i,j)$ on note~$s_{i,j}$ le degré de~$T_{i,j}$ ;
on suppose que pour tout~$i$, les éléments homogènes 
$T_{i,j}f_{i,j}$ ont tous même degré~$r_i$ lorsque~$j$ varie de~$1$ 
à~$m$,
et l'on pose ~$g_i=\sum_j f_{i,j}T_{i,j}$.
On fait les hypothèses suivantes : 

\medskip
$a)$ il existe~$j_0$ tel que la famille~$(f_{i,j_0})_i$ 
soit algébriquement indépendante sur~$K$ ; 

$b)$ la famille~$(T_{i,j})_{i,j}$ est algébriquement indépendante 
sur~$K(f_{i,j})_{i,j}$. 

\medskip
La famille obtenue en concaténant les familles 
$(T_{i,j})_{i,j}$ 
et la famille~$(g_i)_i$ est alors
algébriquement indépendante sur 
$K$. 
}

\medskip
{\em Démonstration.} Soit~$\rho>0$ et soit~$P$ un élément homogène de degré~$\rho$ de
$K[(X_{i,j}/s_{i,j})_{i,j},Y_1/r_1,\ldots, Y_n/r_n]$ tel
que~$P((T_{i,j})_{i,j},g_1,\ldots,g_n)=0$ ; on veut montrer que
$P=0$, et l'on raisonne pour ce faire par l'absurde. On suppose
donc que~$P$ est non nul ; écrivons 
\[ P=\sum a_{(e_{i,j}), (e'_i)} \prod
X_{i,j}^{e_{i,j}}\prod Y_i^{e'_i}.\]
 L'ensemble~$\mathscr E$ des
multi-exposants~$((e_{i,j}), (e'_i))$ tels que ~$a_{(e_{i,j}),
(e'_i)}\neq 0$ est non vide ; soit~$N$ le plus grand des entiers
\[ \sum_i e_{i,j_0}+\sum_i e'_i\]
 pour ~$((e_{i,j}), (e'_i))$
parcourant~$\mathscr E$, et soit~$\mathscr E'$ le sous-ensemble de
$\mathscr E$ formé des exposants en lesquels ce maximum est atteint.

On a \[ P((T_{i,j})_{i,j},g_1,\ldots,g_n)=\sum a_{(e_{i,j}),
(e'_i)} \prod_{i,j} T_{i,j}^{e_{i,j}}\prod_i \left(\sum _j
f_{i,j}T_{i,j}\right)^{e'_i}.\]

Lorsqu'on développe cette expression, on obtient une somme de termes
de la forme~$A_{(e_{i,j})}\prod T_{i,j}^{e_{i,j}}$, où~$\sum
e_{i,j_0}\leq N$, et où~$A_{(e_{i,j})}$ est somme de monômes en les
$f_{i,j}$ qui sont tous homogènes de degré~$\rho\prod 
s_{i,j}^{-e_{i,j}}$
; l'hypothèse b) assure alors que chacun des~$A_{(e_{i,j})}$ est
nul.

Soit~$((e_{i,j}), (e'_i))\in \mathscr E$. Il résulte de la définition
de~$N$ que 
$$a_{(e_{i,j}), (e'_i)}\prod T_{i,j}^{e_{i,j}}\prod
\left(\sum _j f_{i,j} T_{i,j}\right)^{e'_i}$$ 
 s'écrit comme la somme de
$$a_{(e_{i,j}), (e'_i)} \prod f_{i,j_0}^{e'_i} \prod 
T_{i,j}^{e_{i,j}}\prod
T_{i,j_0}^{e'_i}$$ et de monômes en les~$T_{i,j}$ dont le degré
total (monomial) en les variables~$T_{i,j_0}$ est strictement
inférieur à~$N$. Quant au degré monomial total en les variables
$T_{i,j_0}$ du terme~$$a_{(e_{i,j}), (e'_i)}\prod 
f_{i,j_0}^{e'_i}\prod
T_{i,j}^{e_{i,j}}\prod T_{i,j_0}^{e'_i},$$ il est majoré par~$N$
avec égalité si et seulement si~$((e_{i,j}), (e'_i))\in \mathscr
E'$.

Fixons~$((e_{i,j}), (e'_i))\in \mathscr E'$. Pour tout~$(i,j)$,
l'on pose~$e'_{i,j}=e_{i,j}$ si~$j\neq j_0$ et 
$e'_{i,j_0}=e_{i,j_0}+e'_i$.
Il découle de ce qui précède que~$A_{(e'_{i,j})}=\sum 
a_{(\epsilon_{i,j}),
(\epsilon'_i)} \prod f_{i,j_0}^{\epsilon'_i}$, où~$((\epsilon_{i,j}),
(\epsilon'_i))$ parcourt la famille~$\mathscr F$ des multi-exposants
tels que l'on ait pour tout~$(i,j)$ les égalités 
$\epsilon_{i,j}=e'_{i,j}$
si~$j\neq j_0$ et~$\epsilon_{i,j_0}+\epsilon '_i=e'_{i,j_0}$ ;
notons que la donnée de~$(\epsilon'_i)$ détermine entièrement un
multi-exposant de~$\mathscr F$, et que~$((e_{i,j}),(e'_i))\in \mathscr 
F$.

Ainsi,~$A_{(e'_{i,j})}$ apparaît comme un polynôme en les~$f_{i,j_0}$. On a vu 
plus haut que ~$A_{(e'_{i,j})}=0$, ce qui entraîne que 
chacun de ses coefficients est nul en
vertu de l'hypothèse~a) ; 
mais l'un de ses coefficients est égal à~$a_{(e_{i,j}),
(e'_i)}$ qui est non nul, ce qui débouche sur une contradiction et
achève la démonstration.~$\Box$


{\footnotesize \begin{thebibliography}{bb} 

\bibitem{brk1} {\sc V. Berkovich}, {\em Spectral theory and analytic
    geometry over non-archimedean fields}, Mathematical Surveys and
    Monographs {\bf 33}, AMS, Providence, RI, 1990.
\bibitem{brk2} {\sc V. Berkovich}, {\em \'Etale cohomology for
    non-archimedean analytic spaces}, Inst. Hautes Etudes
    Sci. Publ. Math. {\bf 78} (1993) 5-161. 
    \bibitem{loc2} {\sc V. Berkovich}, {\em Smooth~$p$-adic spaces are locally contractible II}, in {\em Geometric Aspects of Dwork Theory}, Walter de Gruyter \& Co., Berlin, 2004, 293-370.
    \bibitem{bgrov} {\sc R. Bieri and J.R.J. Groves}, {\em The geometry of the set of characters induced by valuations}, J. Reine Angew. Math. {\bf 347} (1984), 168-195. 
    \bibitem{bgr}  {\sc S. Bosch, S. G\"untzer and U. Remmert}, {\em Non-Archimedean analysis. A systematic approach to rigid analytic geometry}, Grundlehren der Mathematischen Wissenschaften  {\bf 261}, Springer-Verlag, Berlin, 1984.
\bibitem{fdceb} {\sc A. Chambert-Loir et A. Ducros}, {\em Formes différentielles et courants sur les espaces de Berkovich}, travail en cours. 
\bibitem{cotmk}{\sc B. Conrad and M. Temkin}, {\em Non-Archimedean analytification of algebraic spaces}, preprint.
\bibitem{moder}{\sc A. Ducros}, {\em Les espaces de Berkovich sont modérés}, exposé {\bf 1056} du séminaire Bourbaki. 
\bibitem{imrsk} {\sc A. Ducros}, {\em Image réciproque du squelette par un morphisme entre espaces de Berkovich de même dimension}, Bull. Soc. Math. France {\bf 131} (2003), no. 4, 483--506. 
\bibitem{semialg} {\sc A. Ducros}, {\em Parties semi-alg\'ebriques d'une vari\'et\'e alg\'ebrique~$p$-adique}, Manuscripta Math. {\bf 111} no. 4 (2003), 513-528.
\bibitem{dim} {\sc A. Ducros}, {\em Variation de la dimension relative en g\'eom\'etrie analytique~$p$-adique}, Compositio. Math. {\bf 143} (2007),1511--1532. 
\bibitem{exc} {\sc A. Ducros}, {\em Les espaces de Berkovich sont excellents}, Ann. Inst. Fourier {\bf 59} (2009), no. 4, 1407-1516.
\bibitem{form} {\sc A. Ducros}, {\em Toute forme modérément ramifiée d'un polydisque ouvert est triviale}, à paraître dans Math. Z. 
\bibitem{flatn} {\sc A. Ducros}, \emph{Flatness in non-Archimedean analytic geometry}, preprint. 
\bibitem{hhmelim} {\sc D. Haskell, E. Hrushovski, D. Macpherson}, {\em Definable sets in algebraically closed valued fields: elimination of imaginaries}, J. reine angew. Math. {\bf 597} (2006), 175 -- 236. 
\bibitem{hl} {\sc E. Hrushovski and F. Loeser}, {\em Non-archimedean tame topology and stably dominated types}, preprint. 
\bibitem{angel} {\sc J. Poineau}, {\em Les espaces de Berkovich sont angéliques}, prépublication. 
\bibitem{tmk1} {\sc M. Temkin}, {\em On local properties of non-Archimedean analytic spaces}, Math. Annalen  {\bf 318}, (2000), 585-607. 
\bibitem{tmk2} {\sc M. Temkin}, {\em On local properties of non-Archimedean analytic spaces. II.},  Israel J. Math.  {\bf 140}  (2004), 1-27.
\bibitem{tmk3} {\sc M. Temkin}, {\em A new proof of the Gerritzen-Grauert theorem}, Math. Ann. {\bf 333} (2005), 261-269.
\end{thebibliography} }
\end{document}